\documentclass[11pt,a4paper,leqno]{amsart}

\usepackage{latexsym,enumitem}
\usepackage{amssymb}
\usepackage[cp850]{inputenc}
\usepackage{epsfig}
\usepackage{psfrag}
\usepackage{amsthm}
\usepackage{amscd}
\usepackage{amsmath}
\usepackage{amsfonts}
\usepackage{graphics,caption}
\usepackage[all]{xy}
\usepackage{etoolbox}
\patchcmd{\quote}{\rightmargin}{\leftmargin 2em \rightmargin}{}{}
 \usepackage[sort,numbers]{natbib} 
 \usepackage{mathtools}
\usepackage[normalem]{ulem}
\usepackage{array}

\usepackage[english]{babel}
\usepackage{comment}
\usepackage{color}
\usepackage[colorlinks]{hyperref}

\DeclareMathOperator{\tr}{tr} 

\let\phi\varphi

\DeclareMathOperator{\eqdef}{\coloneqq} 
\newcommand{\f}[2]{\frac{#1}{#2}} 

\let\epsilon\varepsilon
\let\subset\subseteq
\newcommand{\be}{\begin{equation*}}
 \newcommand{\ee}{\end{equation*}}
\newcommand{\bpf}{\begin{dimo}}
\newcommand{\epf}{\end{dimo}}
\newcommand{\bdefi}{\begin{defin}}
\newcommand{\edefi}{\end{defin}}
\newcommand{\bthm}{\begin{thm}}
\newcommand{\ethm}{\end{thm}}
\newcommand{\blem}{\begin{lem}}
\newcommand{\elem}{\end{lem}}
\newcommand{\bcor}{\begin{cor}}
\newcommand{\ecor}{\end{cor}}

\newcommand{\bprop}{\begin{prop}}
\newcommand{\eprop}{\end{prop}}
\newcommand{\bese}{\begin{ese}}
\newcommand{\eese}{\end{ese}}
\newcommand{\brem}{\begin{rem}}
\newcommand{\erem}{\end{rem}}
\newcommand{\bpfc}{\begin{dimoclaim}}
\newcommand{\epfc}{\end{dimoclaim}}

\newcommand{\rar}{\rightarrow} 

\newcommand{\diff}{\mathop{}\!\mathrm{d}} 
	
\newcommand{\abs}[1]{\left\lvert#1\right\rvert}						
\newcommand{\norm}[1]{\left\lVert#1\right\rVert}					
\newcommand{\set}[1]{\left\{#1\right\}}					
	
\newcommand{\quotient}[2]{\left.\raisebox{.1em}{$#1\!$}\middle/\raisebox{-.1em}{$#2$}\right.}

\DeclareMathOperator{\R}{\mathbb R}			
\DeclareMathOperator{\C}{\mathbb C}			
\DeclareMathOperator{\Z}{\mathbb Z}			

\DeclareMathOperator{\im}{Im} 

\DeclareMathOperator{\id}{id} 

\newcommand{\B}{\mathscr B}

\newcommand{\defrac}[2]{\frac{\partial #1}{\partial #2}} 
\newcommand{\difrac}[2]{\frac{\diff#1}{\diff#2}} 

\newenvironment{quot}
{
	\vspace{-0.2cm}
	\vspace{0.2cm}
}

\theoremstyle{definition}
\newtheorem{d1}{Definition}[section] 

\newenvironment{defin}
{
	\begin{quot}
		\begin{d1}
		}
		{\end{d1}
	\end{quot}

}

\theoremstyle{definition}
\newtheorem{r1}[d1]{Remark}

\newenvironment{rem}
{
	\begin{quot}
		\begin{r1}
		}
		{\end{r1}
	\end{quot}
}

\theoremstyle{definition}
\newtheorem{e1}[d1]{Exercise}

\theoremstyle{definition}
\newtheorem{ese1}[d1]{Example}

\newenvironment{ese}
{
	\begin{quot}
		\begin{ese1}
	}
	{	
		\end{ese1}
	\end{quot}
}

\theoremstyle{definition}
\newtheorem{f1}[d1]{Formulation}

\theoremstyle{definition}
\newtheorem{f2}[d1]{Fact}

\theoremstyle{definition}

\theoremstyle{definition}

\theoremstyle{definition}
\newtheorem{t1}[d1]{Theorem}

\newenvironment{thm}
{
	\begin{quot}
		\begin{t1}}
		{\end{t1}
	\end{quot}
}

\theoremstyle{definition}
\newtheorem*{T1*}{Theorem}

\newenvironment{teor*}
{
	\begin{quot}
		\begin{T1*}}
		{\end{T1*}
	\end{quot}
}

\newenvironment{dimo}
{\begin{proof}[Proof]
	}
	{\end{proof}}

\newenvironment{dimoclaim}{\emph{Proof of Claim:}\;}{\hfill$\square$}

	\theoremstyle{definition}
	\newtheorem{l1}[d1]{Lemma}
	
	\newenvironment{lem}
	{
		\begin{quot}
			\begin{l1}}
			{\end{l1}
		\end{quot}
	}
	\theoremstyle{definition}
	\newtheorem{p1}[d1]{Proposition}
	
	\newenvironment{prop}
	{
		\begin{quot}
			\begin{p1}}
			{\end{p1}
		\end{quot}
	}
	
	\theoremstyle{definition}
	\newtheorem{c1}[d1]{Corollary}
	
	\newenvironment{cor}
	{
		\begin{quot}
			\begin{c1}}
			{\end{c1}
		\end{quot}
	}

		\renewenvironment{abstract}
	{\list{}{\rightmargin\leftmargin}%
		\item[\textbf{Abstract:}]\relax}
	{\endlist}

\newtheorem{innercustomthm}{Theorem}
\newenvironment{customthm}[1]
  {\renewcommand\theinnercustomthm{#1}\innercustomthm}
  {\endinnercustomthm}

 \newtheorem*{Theorem*}{Theorem}
 \newtheorem*{Proposition*}{Proposition}
 \newtheorem*{Lemma*}{Lemma}

\newtheorem{innercustomdef}{Definition}
\newenvironment{customdef}[1]
  {\renewcommand\theinnercustomdef{#1}\innercustomdef}
  {\endinnercustomdef}

\usepackage[margin=70pt,marginparwidth=60pt]{geometry}
\usepackage{xcolor}

\renewcommand{\bar}{\overline}

\newcommand{\A}{ d\widetilde{\mathcal A}}

\newcommand{\HH}{{\mathbb H}}
\newcommand{\re}{\mathrm{Re}}

\newcommand{\hess}{\mathrm{Hess}}
\newcommand{\II}{I\hspace{-0.1cm}I}

\renewcommand{\B}[2]{B\left( #1, e^{2#2}#1\right)}

\renewcommand{\bar}{\overline}

\newcommand{\argsinh}{{\mathrm{arsinh}}}
\newcommand{\argtg}{{\mathrm{artg}}}

\newtheorem{theorem}{\rm\bf Theorem}[section]

\newcounter{notes}%

\newcounter{TCcomments}

\newcounter{VGcomments}

\begin{document}

\title[Behaviour of the Schwarzian derivative on long complex projective tubes]{Behaviour of the Schwarzian derivative on long complex projective tubes}

\author{Tommaso Cremaschi}
\address{Tommaso Cremaschi:
Trinity College Dublin, School of Mathematics,
17 Westland Row, Trinity College Dublin, Dublin 2, Ireland}
\email{cremasct@tcd.ie}

\author{Viola Giovannini}
	\address{Viola Giovannini:
	ETH Z\"urich, Department of Mathematics, 
	101 R\"amistrasse
	8092, Switzerland}
\email{viola.giovannini@math.ethz.ch}

\maketitle
\begin{abstract}
	The Schwarzian derivative parametrizes the fibers of the space of complex projective structures on a surface as a vector bundle over its Teichm\"uller space. We study its behaviour on long complex projective tubes, and get estimates for its pairing with infinitesimal earthquakes and graftings. As the real part of the Schwarzian derivative coincides with the differential of the renormalized volume, we obtain bounds for the variation of renormalized volume under complex earthquake paths, and its asymptotic behaviour under pinching a compressible curve. 
\end{abstract}

\section{Introduction and results}

Complex projective structures and Riemann surfaces appear in multiple areas of math, in particular, they arise in the study of convex co-compact (or, more generally, geometrically finite) hyperbolic $3$-manifolds as naturally induced structures on their boundary at infinity. This leads to a parametrization of the deformation space of convex co-compact metrics on a hyperbolic $3$-manifold via the so-called conformal boundary at infinity \cite{Ma2016,MT1998}.

The deformation space of complex projective structures on a closed surface $S$ of genus $g\geq 2$ forms a vector bundle over the moduli of underlying Riemann surfaces, whose fibers are parameterized by the Schwarzian derivative of the developing map $f$ comparing a complex projective structure to its \textit{Fuchsian} conformal one. Depending on $f$ being univalent or not, its Schwarzian derivative $S(f)$ behaves quite differently. In this work, we investigate the case in which the restriction of $f$ to a `banana' neighborhood of a lift of a simple closed curve is a covering map onto an annulus. In particular, $f$ is not univalent. A case falling under this setting is when $f:\mathbb H^2\rar \Omega\subsetneq\mathbb{CP}^1$ is a covering map onto a \textit{non-simply-connected hyperbolic} domain. This case arises naturally in the context of convex co-compact hyperbolic $3$-manifolds with compressible boundary, and their renormalized volume. This is, in fact, the main application of the present work.

Recall that two Riemann surfaces are isomorphic if there exists a biholomorphic diffeomorphism between the two. Two complex projective surfaces are isomorphic if there exists a diffeomorphism between the two pulling back the atlas of projective charts of the second one into the one of the first one. In both cases, when the diffeomorphism is isotopic to the identity, we say that the structures are \textit{marked isomorphic}. Riemann structures on an orientable compact surface $S$ admit a large deformation space, which, when considered up to marked isomorphisms, is called the \textit{Teichm\"uller space} of $S$, and it is denoted by $\mathcal{T}(S)$. Complex projective structures on $S$, up to marked isomorphism, admit an even larger deformation space which is identified, through the forgetting projection $\pi$, with a vector bundle $\mathcal{P}(S)$ over $\mathcal{T}(S)$ (see \cite{dumas-survey}). 

The fibers of $\mathcal{P}(S)$ are parameterized by the \textit{Schwarzian derivative} of the \textit{developing map} of the complex projective structure.
The Schwarzian derivative of a locally injective holomorphic map $f\colon D\rightarrow \mathbb{CP}^1$, for an open domain $D\subseteq \mathbb{C}$, is the central object of this paper, and it is the holomorphic quadratic differential defined as
\[S(f)=\biggl(\biggl(\dfrac{f^{''}}{f^{'}}\biggr)^{'}-\dfrac{1}{2}\biggl(\dfrac{f^{''}}{f^{'}}\biggr)^{2}\biggr)dz^2~,\]
where $z$ denotes the complex coordinate of $\mathbb{C}$, and $dz^2=dz\otimes dz$ the associated holomorphic symmetric $(0,2)$-tensor.
 A projective structure $Z$ on a surface $S$ can be lifted to the universal cover $\widetilde{S}$, and the \textit{developing map} of $Z$ is an immersion $f\colon \widetilde{S}\rightarrow \mathbb{CP}^1$ such that it restricts to projective charts on small enough opens of $\widetilde{S}$.  Also, recall that the group of \textit{M\"obius transformations} $\mathrm{PSL}(2,\mathbb C)$ is the group of biholomorphisms of $\mathbb{CP}^1$ and that the developing map is unique up to post-composing with a M\"obius transformation.

The Schwarzian derivative $S(f)$ is equal to zero if and only if $f$ is a M\"obius transformation, and it therefore measures how far a projective structure is from the \textit{Fuchsian} one in $\mathcal{T}(S)$, identified with the zero of the fiber (see Section \ref{hypmanifolds}). A Fuchsian complex projective structure $Z_0$ on a closed surface $S$ identifies the universal cover $\widetilde{S}$ with $\mathbb{H}^2$, and the fundamental group $\pi_1(S)$ with a discrete subgroup $\Gamma< \mathrm{PSL}(2,\mathbb R)$, for any point of its own fiber. The developing maps of the points $Z$ in $\mathcal{P}(S)$ in the fiber of $Z_0$ are then ($\Gamma$, $\Gamma'$)-equivariant, with $\Gamma'< \mathbb{P}SL(2,\mathbb{C})$ the image of the complex projective holonomy $\rho$ of $Z$. The pair $(f,\rho)$ determines $Z$ uniquely. We simply call the Schwarzian derivative of the developing map $f\colon \mathbb{H}^2\rightarrow \mathbb{CP}^1$ of a complex projective structure $Z$ on $S$, the \textit{Schwarzian of $Z$}.

Given $Z$, by pulling back to the universal cover $\widetilde{S}$ the hyperbolic metric conformal to $\pi(Z)$, we can define the \textit{injectivity radius} of $Z$ to be $\delta$ if any developing map $f\colon \widetilde{S}\rightarrow \mathbb{CP}^1$ is univalent (i.e. holomorphic injective) on every hyperbolic ball of radius $\delta$, and $\delta$ is the largest positive number with this property. When $f$ is defined on the upper half-plane $(\mathbb{H}^2, \rho\abs{dz}^2)\subseteq \mathbb{C}$, with $\rho\abs{dz}^2$ the unique complete hyperbolic metric, we define the \textit{infinity norm} of $S(f)=qdz^2$ as
  \[\norm{S(f)}_\infty=\sup_{z\in\mathbb{H}^2} |q(z)|/\rho(z)~.\] 
If $f\colon \mathbb{H}^2\rightarrow \mathbb{CP}^1$ is univalent, the injectivity radius is $+\infty$, its image is a disk and, by the Kraus-Nehari bound \cite{Ne1949}, the infinity norm of $S(f)$ is bounded:
\[\norm{S(f)}_\infty\leq 3/2~.\] 
If, instead, $f$ is not univalent, the infinity norm is bounded from above and below by functionsmeglio  of the injectivity radius.

More precisely, by Kra-Maskit \cite[Lemma $5.1$]{KM}:
\[\norm{S(f)}_\infty\leq \f32\coth^2(\delta/2)~,\]
and
 \[\norm{S(f)}_\infty\geq \f12\coth^2(\delta/2)~.\]
  In particular, when $\delta\sim 0$, then
  \[\norm{S(f)}_\infty\gtrsim 2/\delta^2~.\]
Thus, it is interesting to further investigate the behaviour of the Schwarzian of a complex projective structure $Z$ on a closed surface $S$ when its injectivity radius $\delta$ is small (hence, $f$ is not univalent), and, consequently, the norm of the Schwarzian is large. Motivated by Kleinian groups, with an application to renormalized volume and hyperbolic 3-manifolds, we are interested in bounding the Schwarzian along hyperbolic geodesics that have \emph{wide} collar neighborhoods. We now introduce our current setting and our main Definition \ref{defilongprojtube}. 

For a hyperbolic surface $X$, simple closed geodesics of length $\ell_X(\gamma)\leq \epsilon_0=2\sinh^{-1}(1)$ have disjoint collars of positive width, whose widths tend to infinity as $\ell_X(\gamma)\rar 0$, see \cite[Theorem $4.1.6$ (Collar Lemma)]{Bu1992}. The large collars are often referred to as \emph{long tubes}, while $\gamma$ is called \textit{short}. We say that an $L$-neighborhood\footnote{Here we are taking a closed neighborhood} of a lift $\widetilde{\gamma}\subseteq\widetilde{X}$ of a short curve $\gamma$ is \textit{wide}, if $L$ is the one given by the Collar Lemma (see Definition \ref{margulitubedefi}). Identifying $\widetilde{S}=\mathbb{H}^2$, this is a banana neighborhood of width $L$. Before introducing our main object of study we also need the following definition. Given $f\colon \mathbb{H}^2\rightarrow \mathbb{CP}^1$ locally univalent, a topological disk $D\subseteq\mathbb{H}^2$ is called \textit{$f$-round} if the restriction $f\vert_{D}$ is univalent and $f(D)$ is a round disk in $\mathbb{CP}^1$. The set of $f$-round disks of $\mathbb{H}^2$ is partially ordered by inclusion. For each point $z\in\mathbb H^2$, there is a distinguished choice of a maximal element in the collection of $f$-round disks containing $z$, which we call \textit{the maximal disk $D_z$} for $f$ at $z$ (see Section \ref{ssc:thmetric} for the definition of $D_z$). 
	
	We now introduce a notion of long tube in the setting of complex projective structures.

\begin{customdef}{\ref{defilongprojtube}}
	
Let $Z=(f,\rho)$ be a complex projective structure on $S$, and $X=\pi(Z)$ be the induced conformal hyperbolic surface. Let $\gamma\subset X$ be a short simple closed geodesic. We say that $Z$ contains a nontrivial \textit{long tube} with core $\gamma$ if, for any lift $\widetilde\gamma\subset \widetilde X$, up to post-composition by a M\"obius map, the restriction of the developing map $f$ to the wide neighborhood $\widetilde{\mathcal{A}}$ of $\widetilde{\gamma}$ is a covering onto an essential annulus $\mathcal{A}\subseteq\mathbb{C}^*\subseteq \mathbb{CP}^1$ with deck group $\langle\gamma\rangle$, and, $\forall z\in \widetilde{\mathcal A}$, the maximal disk $D_z$ satisfies $f(D_z)\subset\C^*$. We call $\mathcal{A}$ a \textit{long complex projective tube}.
  \end{customdef}

 	We point out that by $\langle\gamma\rangle$ we denote the group generated by the loxodromic isometry of $\mathbb{H}^2$ translating along $\widetilde{\gamma}$ associated to $\gamma$ by identifying the fundamental group of $S$ with the discrete subgroup of $\mathbb{P}SL(2, \mathbb{R})$ realizing $X$. Moreover, up to M\"obius transformation, any proper annulus $\mathcal{A}\subseteq \mathbb{CP}^1$ can be assumed to be essential in $\mathbb{C}^*$. Hence, what one is really asking in Definition \ref{defilongprojtube} is that $\cup_{z\in\widetilde{\mathcal A}}f(D_z)$ does not saturate a complementary component of $\mathcal A$. Note that whenever the deck group of a covering $f$ onto a long complex projective tube $\mathcal{A}\subseteq \mathbb{CP}^1$ acts by isometries of $\widetilde{\mathcal{A}}$ (i.e. it is a subgroup of isometries fixing $\tilde\gamma$), we can equip $\mathcal{A}$ with a hyperbolic metric by pushing forward the one on $\widetilde{\mathcal{A}}\subseteq\mathbb{H}^2$ via $f$. On a long tube, the (restriction of the) injectivity radius coincides with half of the hyperbolic length of the geodesic representative $\gamma$ of the unique non-trivial simple closed curve in $\mathcal{A}$, which is called the \textit{core} of $\mathcal{A}$.
 	
 	In the next remark we discuss some of the conditions in Definition \ref{defilongprojtube}.
	
	 \brem Note that assuming the restriction of $f$ being a covering onto an annulus, with deck group $D$ acting by isometries on ${\A}$, does not imply $D=\langle \gamma\rangle$ nor that the complex projective holonomy of $\gamma$ is trivial. However, one has two distinct cases:
 \begin{itemize}
 \item [(i)] $\mathcal A$ is a round annulus in $\mathbb C^*$;
 \item[(ii)] $\mathcal A$ is not round.
 \end{itemize} 
 Since $\rho(\gamma)(\mathcal A)=\mathcal A$ and $\mathcal{A}$ is proper, the holonomy of $\gamma$ in $\mathbb{P}SL(2,\mathbb{C})$ is forced to be an elliptic M\"obius transformation.
Moreover, if $\mathcal A$ is not round, then it has to be finite-order.  This implies that, in the latter case, the deck group $D\cong\langle\eta\rangle$ is rational with respect to $\gamma$, that is, $\eta^n=\gamma^m$ for some positive integers $m,n$. Thus, up to taking common finite index covers we can always return to the more general case of Definition \ref{defilongprojtube} and, since $S( \tilde f)=p^*(S( f))+S(p)$, for $p$ the covering map and $\tilde f$ the developing map of the covering annulus, a bound on the infinity norm of $S(\tilde f)$ and $S(p)$ descends to a bound on the one of $S(f)$. Then, we can bound the infinity norm of $S(\tilde f)$ either by using our methods (in the case of core curves of length $\leq \epsilon_0$) or Kraus-Nehari while, in general, the infinite norm of $S(p)$ depends on the degree of the cover. 

The round annulus case is instead similar to the Toy Model discussed in Section \ref{sec:symmsch}. 
Thus, up to changing some of the numbers in the statements of Theorems \ref{Schwarzian} and \ref{pairing} the same kind of bounds hold. However, as we are mainly interested in the asymptotic behavior, which is useful in the case of Kleinian groups, we restrict ourselves to this more restrictive setting for ease of notation and readability.  One could go through the paper and replace any asymptotic $\ell(\gamma_t)\rar 0$ with $\ell(\eta_t)\rar 0$, and have the error terms depending on $\ell(\eta_t)$ and $L_t=L(\ell(\eta_t))$, where the latter essentially depends on the rational relation between $\eta_t$ and $\gamma_t$ (essentially the degree of the needed covers).\erem

  In the following theorem we present the behaviour of the Schwarzian on long complex projective tubes.

 \begin{theorem}\label{Schwarzian} Let $\mathcal{A}$ be a long tube of a complex projective surface $Z$, and let $\gamma$ be its core of length $\ell\leq \varepsilon_0$. Let also $f\colon \mathbb{H}^2\rightarrow \mathbb{CP}^1$ be the developing map of $Z$, and let $S(f)$ be its Schwarzian derivative. Then, in a tubular neighborhood of the lift $\widetilde{\gamma}=i\mathbb{R}_{+}$, the Schwarzian $S(f)$ behaves as follows 
 
 	 \[S(f)=\f{1}{2z^2}\left(1+\f{4\pi^2}{\ell^2}\right)dz^2+q_\ell(z)dz^2~,\]
	where $z$ is the complex coordinate of $\mathbb{H}^2$, and $\norm{q_\ell(z)dz^2}_{\infty}$ is bounded and in $O\left(\f{e^{-\pi^2/(2\ell)}}{\ell^2}\right)$.
 
 \end{theorem}
 
Note that the same result holds for any other lift $\widetilde{\gamma}$, after a M\"obius change of coordinates sending $i\mathbb{R}_{+}$ to $\widetilde{\gamma}$.
 
 As the space of holomorphic quadratic differentials on a Riemann surface $X$ is identified with the cotangent space $T^{*}_{X}\mathcal{T}(S)$ of the Teichm\"uller space at $X$, the Schwarzians of the complex projective structures in the fiber of $X$ can be paired with vectors in the tangent space $T_{X}\mathcal{T}(S)$. Vectors in $T_{X}\mathcal{T}(S)$ can be expressed as \textit{harmonic Beltrami differentials}, i.e. as objects of the type $(\overline{f}/\rho) \partial z\otimes d\bar{z} $, with $f$ a holomorphic function on $X$ and $\rho$ such that $\rho |dz|^2$ is the unique hyperbolic representative in $X$. The pairing is then realized by integration on the surface $X$ (see \cite{Hubbard2016}, \cite{QTT}).
 
  The real part of the Schwarzian derivative, as a holomorphic quadratic differential, is of particular interest as it coincides with the differential of the \textit{renormalized volume} $V_R$, see \cite[Lemma $8.3$]{S2008} or \cite[Corollary $3.11$]{compare}. In particular, the variation of renormalized volume along paths in the deformation spaces of convex co-compact manifolds with incompressible boundary, i.e. in the Teichm\"uller space of their boundary, has been extensively studied in the case of flow-lines of the gradient flow, see \cite{BBB2018}, \cite{bridgeman-brock-bromberg:gradient}, \cite{bridgeman-bromberg-pallete:convergence}. For compressible boundary the layout of this problem changes: in general $V_R$ does not converge under gradient flow as the norm of its derivative is unbounded, see for example Theorem \ref{mainVR}, and \cite{averages}.

  \textit{Earthquaking} and \textit{grafting} are two natural paths inside the Teichm\"uller space, and by infinitesimal we mean the tangent vector at time zero of their respective deformation paths (see \cite{McMcomplexearth}, or Section \ref{complexearth}). In brief, an earthquake on a simple closed curve $\gamma\subseteq X$ consists in cutting along $\gamma$, twisting by a certain parameter $t\in \mathbb{R}$ the left-hand side of the surface, and gluing back. Again in short, grafting on $\gamma$ means to cut along $\gamma$ and attach a Euclidean tube of some height $s\in\mathbb{R}^{+}$, this gives a new well-defined conformal structure on $X$. Infinitesimal earthquake and infinitesimal grafting on the same (multi)curve at $X$ are orthogonal with respect to the complex structure on $\mathcal{T}(S)$ \cite[Theorem $2.10$]{McMcomplexearth}.
  
In the following theorem we show that the real part of the pairing of the Schwarzian of a complex projective surface containing a long tube is almost trivial along \textit{infinitesimal earthquakes} on the core, while, asymptotically, it grows as $-\frac{\pi^2}\ell$ in the \textit{infinitesimal grafting} direction. In the following result $F_e(\ell)$ and $F_{gr}(\ell)$ are to be interpreted as real-valued functions estimating the error along an earthquake and grafting path in Teichm\"uller space, respectively. A version of Theorem \ref{pairing} with explicit bounds is stated and proved at the end of Section \ref{sec:graft}.

\begin{theorem}\label{pairing}
	Let $Z$ be a complex projective surface containing a long tube of core $\gamma$, let $S(f)$ be its Schwarzian, and let $X=\pi(Z)$ its underlying Riemann surface. Let also $\mu$ and $\nu$ be, respectively, the infinitesimal earthquake and grafting on the simple closed curve $\gamma\subseteq X$ of hyperbolic length $\ell\leq\varepsilon_0$. Then 
	\[\abs{\re\left\langle S(f), \mu\right\rangle}\leq F_{e}(\ell) \] 
	and 
	\[\abs{\re\left\langle S(f), \nu\right\rangle+\f{\pi^2}{\ell}}\leq  \f14\ell+ F_{gr}(\ell) \]
	with $F_{e}(\ell)$ and $F_{gr}(\ell)$ two explicit functions such that \[\abs{F_{e}(\ell)},\ \abs{F_{gr}(\ell)}\leq C\f{e^{-\pi^2/(2\ell)}}{\ell}~,\]
	for some constant $C>0$.
\end{theorem}
A heuristic of why one should believe Theorem \ref{pairing} is presented in Section \ref{sec:symmsch}, where we show the result for the \emph{toy model} case of a \textit{symmetric} complex projective tube.

\textbf{Connection to Hyperbolic $3$-manifolds.} A  hyperbolic $3$-manifold $(M, g)$ is \textit{convex co-compact} if it is homotopy equivalent to the interior of a compact convex subset $C\subseteq M$ with boundary components of genus $g\geq 2$. The Riemannian metric $g$ does not extend to the boundary $\partial\bar{M}$ as it diverges in its collar, in particular, $M$ has infinite hyperbolic volume. However, $M$ can be conformally compactified, i.e. $g$ induces a well-defined complex structure on $\partial\bar{M}$, which is called the \textit{conformal boundary at infinity} of $M$. As a consequence of $\mathbb{CP}^1$ being the conformal boundary of the closure $\overline{\mathbb{H}^3}$ of the $3$-hyperbolic space $\mathbb{H}^3$, and of the fact that the orientation preserving isometry group of $\mathbb{H}^3$ coincides with $\mathbb{P}SL(2, \mathbb{C})$, the conformal boundary at infinity is naturally equipped with a complex projective structure. This is more explicitly obtained as the quotient of the discontinuity domain $\Omega(\Gamma)\subseteq\mathbb{CP}^1$ of the Kleinian group $\Gamma<\mathbb{P}SL(2,\mathbb{C})$ realizing $M$, by $\Gamma$ itself. In this case, the associated developing map is then a covering onto its image $\Omega(\Gamma)$. Up to quotienting $\mathcal{T}(\partial\bar{M})$ by the group generated by the Dehn twists on \textit{compressible} curves, the conformal boundary at infinity realizes a homeomorphism between the deformation space of convex co-compact hyperbolic metrics on $M$ up to homotopy, denoted by $CC(M)$, and $\mathcal{T}(\partial\bar{M})$ \cite[Thm 5.1.3.]{Ma2016}. An \emph{essential} curve in $\partial\bar{M}$ is compressible if it bounds a disk in $M$, the boundary $\partial\bar{M}$ is \textit{incompressible} if there are no compressible curves. The \textit{renormalized volume} $V_R(\cdot)$ is a real-analytic function on $CC(M)$, or, equivalently, on $\mathcal{T}(\partial\bar{M})$ (see \cite{S2008}). Depending on whether $\partial\bar{M}$ is incompressible or not, the renormalized volume behaves very differently. If $\partial\bar{M}$ is incompressible  $V_R$ is always non-negative (see \cite{VP2016}, \cite{BBB2018}). If $\partial\bar{M}$ is \textit{compressible}, the work of Canary and Bridgeman in \cite{BC2017}, and of Schlenker and Witten in \cite{averages}, show that the infimum of the renormalized volume is $-\infty$. The condition for $M$ of having compressible boundary translates into having a non-univalent developing map of the natural complex projective structure induced on the boundary at infinity $\partial \overline{M}$, with image a non-simply-connected hyperbolic domain of $\mathbb{CP}^1$. We can then exploit Theorems \ref{Schwarzian} and \ref{pairing} to get information on the behaviour of the renormalized volume function in the compressible boundary case. Note, indeed, that the developing map of any complex projective structure naturally arising on the boundary at infinity of a hyperbolic $3$-manifold with compressible boundary satisfies the hypothesis of Definition \ref{defilongprojtube}, provided that the compressible loop is short enough.

As already pointed out, the differential of the renormalized volume $(d V_R)_X(\cdot)$ at a point $X$ coincides with the real part of the Schwarzian derivative of the developing map $f_M$ of the projective structure associated to the corresponding convex co-compact $3$-manifold $M$. Theorem \ref{Schwarzian} and, thanks to our key Lemma \ref{intbypart} (see also \cite[Lemma $5.1$]{CGS2024}),  Theorem \ref{pairing} can  then be rephrased in terms of $d V_R$.

In \cite{McMcomplexearth} McMullen introduced the notion of $\lambda$-\textit{complex earthquake} with $\lambda=t+is\in\mathbb{H}^2$, which consists in moving inside the Teichm\"uller space with first a parameter $t$ earthquake, and then a parameter $s$ grafting, on the same simple closed curve (actually, more generally, on the same measured lamination). With the following two theorems, we bound the change of renormalized volume along a complex earthquake path on a short simple closed geodesic. More specifically, in Theorem \ref{mainVR} we bound the change of renormalized volume along an earthquake of parameter $t$, while in Theorem \ref{asymptotic} we calculate the change along a grafting path of parameter $s$. The next is a corollary of the first part of Theorem \ref{pairing}, and the key Lemma \ref{intbypart}.

\begin{theorem}\label{mainVR}
		Let $M$ be a convex co-compact hyperbolic $3$-manifold. Let $X_0\in \mathcal{T}(\partial\overline{M})$, and let $X_t\in \mathcal{T}(\partial\overline{M})$ be the Riemann surface obtained by a parameter $t\in\mathbb{R_+}$ earthquake on $X_0$ along a simple closed compressible geodesic $\gamma$ in $\partial\bar{M}$, of length $\ell\leq\varepsilon_0$. Then, the following estimate for the renormalized volume of the associated convex co-compact manifolds $M_0$ and $M_t$ holds:
	\[\abs{V_R(M_t)-V_R(M_0)}\leq 
	F(\ell)t~,\]
	with $F(\ell)$ an explicit function of $\ell$ such that \[|F(\ell)|\leq C\f{e^{-\pi^2/(2\ell)}}{\ell}\] for some explicit constant $C>0$.

\end{theorem}

When $\ell$ is small enough, Theorem \ref{mainVR} gives a sharper bound, for the change of renormalized volume under earthquake, than the one of Theorem $1.4$ in \cite{CGS2024}. 

As another application of our estimate for the Schwarzian, we obtain the asymptotic behavior of the renormalized volume under \textit{pinching} a compressible curve in the boundary at infinity of a convex co-compact hyperbolic $3$-manifold, recovering in particular Theorem $A.15$ in \cite{averages}. To pinch a simple closed curve in a Riemann surface means to send to $0$ the length of its geodesic representative with respect to the hyperbolic metric. Pinching on $\gamma$ can be realized by grafting on it with a parameter tending to infinity. This gives a bound for the change of renormalized volume under grafting a short enough curve.

\begin{theorem}\label{asymptotic}
	Let $M_0$ be a convex co-compact hyperbolic $3$-manifold, and let $\gamma\subseteq\partial\bar{M}$ be a simple closed compressible curve of length $\ell_0(\gamma)\leq\varepsilon_0$ in its boundary. The composition of the renormalized volume with the grafting path $(M_s)_{s\in[0,\infty)}$ satisfies
		 \[V_R(M_s)-V_R(M_0)= -\f{\pi^3}{\ell_{s}(\gamma)}+\f{\pi^3}{\ell_0(\gamma)}+\f{\pi}{4}(\ell_s(\gamma)-\ell_0(\gamma))+C(s)~,\]
	with $C(s)$ bounded and such that $C(s)\rar 0$ as $s\rar 0$. In particular, when $s\rightarrow \infty$ the renormalized volume diverges as $ -\f{\pi^3}{\ell_s(\gamma)}$.
\end{theorem}

Moreover, the tools developed can give a sharper estimate of the error function $C(s)$ in Theorem \ref{asymptotic}, see Remark \ref{sharperboundschwarzian}.

 \subsection{Outline of the paper} 

Section \ref{Background} contains the main objects and notions that will be needed, and the key Lemma \ref{intbypart}. 

In Section \ref{sec:symmsch} we give a motivating example by explicitly calculating the Schwarzian derivative of \textit{symmetric} complex projective tubes. This is done in Proposition \ref{symmetric}, which is the analogue of Theorem \ref{mainthm} in this setting. 

Section \ref{sc:estimates} is dedicated to the study of the Schwarzian of long tubes in the general case. In Subsection \ref{Schwarziantensor} we prove Theorem \ref{Schwarzian}. In Subsection  \ref{ssc:estimatesearth}, we give bounds for the real part of the pairing between Schwarzian and infinitesimal earthquakes, while in Subsection \ref{sec:graft} there are the analogous bounds for the pairing with infinitesimal grafting. These are obtained by the tool of the Osgood-Stowe tensor, as explained in Section \ref{Schwarziantensor}. The proof of Theorem \ref{mainthm} is at the end of Section \ref{sec:graft}. Section \ref{ssc:fourier} contains the Fourier analysis leading to the key estimates exploited in the previous subsections. These are used accurately in Remarks \ref{Mell} and \ref{utube}, which can also be applied to obtain sharper versions of Theorem \ref{asymptotic}, see Remark \ref{sharperboundschwarzian}.

Finally, Section \ref{sc:appl} contains the applications of the Schwarzian bounds to the renormalized volume for convex co-compact hyperbolic manifolds with compressible boundary, namely Theorem \ref{mainVR} (stated in a more explicit form), and Theorem \ref{asymptotic}.

 In Appendix \ref{standardmetrics}, we recall relations between the Osgood-Stowe tensor \cite{OS1992} (and a non-traceless extension) and the Schwarzian derivative, which can also be found in \cite{Epstein-wvolume-OS}.

\subsection*{Acknowledgements. }The authors would like to thank Jean-Marc Schlenker for stimulating conversation and for allowing us to include Appendix A, which he wrote. Tommaso Cremaschi was partially supported by MSCA grant 101107744--DefHyp. Viola Giovannini was supported by FNR AFR grant 15777, and the Swiss State Secretariat for Education, Research and Innovation (SERI): MB25.00004.

\tableofcontents

\section{Notation and background}\label{Background}
In this section we recall the main objects and tools that we will use in this work.

\subsection{Riemann, Complex Projective, and Hyperbolic Surfaces}\label{hypmanifolds}
A \textit{Riemann surface} (of genus $g\geq 2$) is a smooth $2$-manifold with an atlas of charts into opens (of the upper half-plane $\mathbb{H}^2$) of $\mathbb{C}$, whose transition maps are biholomorphisms. A \textit{complex projective surface} is a smooth $2$-manifold with an atlas consisting of charts into opens of $\mathbb{CP}^1$ whose transition maps are restriction of M\"obius transformations $\phi\in\mathbb{P}SL(2, \mathbb{C})$, i.e. restriction of biholomorphisms of $\mathbb{CP}^1$. This can be rephrased by saying that a complex projective surface is a $(\mathbb{P}SL(2, \mathbb{C}), \mathbb{CP}^1) $ structure. Since M\"obius transformations are holomorphic, a complex projective surface $Z$ also defines a Riemann surface $X=\pi(Z)$. By Riemann's Uniformization Theorem, a complex surface of genus $g\geq 2$ has a $(\mathbb{P}SL(2, \mathbb{R}), \mathbb{H}^2)$ structure.  In particular, since $\mathbb{P}SL(2, \mathbb{R})< \mathbb{P}SL(2, \mathbb{C})$, a Riemann surface is also a complex projective surface, which is called \textit{Fuchsian}.

The $2$-dimensional hyperbolic space is the unique simply connected complete Riemannian manifold of constant (sectional) curvature $-1$. Here, we will use its \textit{upper half-space} model \[\left(\mathbb{H}^2, \rho(z)|dz|^2\right)~,\] with $\mathbb{H}^2\subseteq\mathbb{C}$ the set of complex vectors of positive imaginary part, $z=x+iy$ the complex coordinate of $\mathbb{C}$, $|dz|^2=dx^2+dy^2$ and $\rho(z)=1/y^2$ (see \cite[Chapter $2$]{Mar16}). Every smooth surface of negative Euler characteristic can be equipped with a smooth complete hyperbolic Riemannian metric, and expressed as a quotient \[\mathbb{H}^2/\Gamma\]
with $\Gamma< \mathbb{P}SL(2, \mathbb{R})=\text{Isom}^{+}(\mathbb{H}^2)$ a discrete and torsion-free subgroup of the group of orientation preserving isometries of $\mathbb{H}^2$.

By Riemann's Uniformization Theorem \cite[Theorem 1.1.1]{Hubbard2016}, every Riemann structure on $S$ of genus $g\geq 2$ can be regarded as a conformal class of metrics on $S$ where \[g\sim h\  \text{  iff  }\  h=e^{2u}g\]
for a smooth function $u\colon S\rightarrow \mathbb{R}$, in which there exists a unique hyperbolic representative. In this paper, we will use all three viewpoints, and move from one viewpoint to another depending on which one is more suited.

We now describe the setup we are interested in. We are concerned with points $Z$ in the vector bundle $\pi\colon\mathcal{P}(S)\rightarrow \mathcal{T}(S)$ of complex projective structures on $S$, and the Schwarzian derivative of their developing maps (see \cite{dumas-survey}). Consider a complex projective structure $Z\in\mathcal{P}(S)$ over a hyperbolic surface $X=\pi(Z)$, such that $Z$ has holonomy representation $\rho\colon\pi_1(S)\rar  \mathbb{P}SL(2, \mathbb{C})$, and let $\rho_0\colon \pi_1(S)\rightarrow \mathbb{P}SL(2, \mathbb{R})$ be the Fuchsian representation of a fixed uniformizing map $\widetilde{S}\rightarrow\mathbb{H}^2$ for $X$. We assume that $\rho$ has non-trivial kernel $K$, and that, for every simple closed curve $\gamma$ whose class is in $K$, the restriction of the developing map to the wide neighborhoods $\widetilde{\mathcal{A}}$ of any lift of $\gamma$ is a covering onto an annulus, with deck group $\langle\gamma\rangle$. This yields a non-univalent, $(\rho_0, \rho)$-equivariant developing map $f\colon\mathbb H^2\rar \mathbb {CP}^1$. If we further assume that the image of each maximal disk $D_z$ (see Section \ref{ssc:thmetric}), for $z\in\widetilde{\mathcal{A}}$, is contained in $\mathbb{C}^*$, then the image of $f$ contains a long complex projective tube. The restriction of $f$ to each wide neighborhood descends to a univalent map $\hat f$ on the image annulus $\mathcal{A}$. Since the restriction of $f$ is a covering map with deck group $\langle\gamma\rangle$ acting by isometries of $\mathbb{H}^2$, one can push-forward the metric of $\widetilde{\mathcal{A}}/\langle\gamma\rangle$ on $\mathcal{A}$. In particular, we can talk about lengths of curves, and the injectivity radius of points in $\mathcal{A}$. 

This setting naturally arises when considering the conformal boundary $\partial_{\infty}N$ of a convex co-compact hyperbolic $3$-manifold $N$ with boundary $\partial N$ being compressible, as explained in Section \ref{hypback}. In this case, the developing map is a covering whose image is a non-simply-connected hyperbolic domain $\Omega\subsetneq\mathbb{CP}^1$, and the complex projective holonomy will have non-trivial kernel.

\subsubsection{Teichm\"uller space and Weil-Petersson pairing } We define here the deformation space of Riemann structures over a closed oriented surface $S$ of genus $g\geq 2$, by using the hyperbolic geometry interpretation (the others are analogous). For a complete reference we refer to \cite[Chapter 6-7]{Hubbard2016} and \cite[Section $10$]{FM2011}.
The \textit{Teichm\"uller space of $S$} can be defined as:
 \[\mathcal{T}(S)= \{h\ |\ h\text{ is a complete hyperbolic metric on $S$}\}/\text{Diffeo}_0(S)~,\]
where $\text{Diffeo}_0(S)$ is the group of diffeomorphisms of $S$ isotopic to the identity, and it acts by pull-back.

The tangent and the cotangent space of the Teichm\"uller space of $S$ at the Riemann surface $X$ are identified (see \cite[Section $6.6$]{Hubbard2016}, or \cite{QTT}), respectively, with the space of \textit{harmonic Beltrami differentials} and \textit{holomorphic quadratic differentials} on $X$. Given a Riemann surface $X$, we denote by $z$ and $\bar{z}$ its local holomorphic and anti-holomorphic coordinate. A Beltrami differential on $X$ is a $(1,-1)$-tensor $\mu$ that can be expressed in local coordinate as 
\[\mu(z)= \eta(z)\partial z\otimes d\bar{z}~,\]
with $\eta$ a measurable complex-valued function. A holomorphic quadratic differential $Q$ is a symmetric $(0,2)$-tensor that can be locally written as \[Q(z)=q(z)dz^2\] with $q(z)$ a holomorphic function. A Beltrami differential $\mu$ is \textit{harmonic} if there exists a holomorphic quadratic differential $qdz^2$ such that \[\mu=(\bar{q}/\rho)\partial z\otimes d\bar{z}\]
where $\rho |dz|^2$ is the unique hyperbolic metric representative in $X$. We denote the space of harmonic Beltrami and holomorphic quadratic differentials on $X$, respectively, with $B(X)$ and $Q(X)$. For any $\mu\in B(X)$ and any $Q\in Q(X)$, with $\mu=\eta \partial_z\otimes d\bar{z}$ and $Q=qdz^2$, it is possible to define the natural non-degenerate pairing \[\langle Q, \mu\rangle = \int_{X} (q\eta) dx\wedge dy\]
with $z=x+iy$ the local complex coordinate of $X$. We will be mostly concerned with its real part \[\re\langle Q, \mu\rangle = \re\int_{X} (q\eta) dx\wedge dy~.\]

\subsubsection{Margulis tubes}\label{ssc:margtubes}
As already pointed out in the introduction, our focus will be on long tubes. We have already defined there the notion of long tube in the complex projective setting. We give here the definition of a particularly useful long tube in the hyperbolic viewpoint. 

\bdefi \label{margulitubedefi}
A \emph{thin tube} in a hyperbolic surface $X$, is the set of points $\mathbb T(\ell)$ around a simple closed geodesic $\gamma$ of length $\ell\leq \epsilon_0=2\argsinh(1)$ that are at a distance at most $L\eqdef\argsinh\left(\f1{\sinh\left(\f\ell 2\right)}\right)$.
\edefi
The hyperbolic metric on $\mathbb T(\ell)$ can be written as 
\[d\rho^2+\left(\f{\ell}{2\pi}\right)^2\cosh^2(\rho)d\theta^2~,\] 
with $\theta\in[0,2\pi]$ and $\rho\in[-L,L]$, for details see \cite[Thm 4.1.1]{Bu1992}. Moreover, in $\mathbb T(\ell)$, the injectivity radius is bounded as 

\[\f{\ell}2 \leq \operatorname{inj}(p)=\argsinh\left(\sinh(\ell/2)\cosh(L-d)\right),\qquad d=d(p,\partial \mathbb T(\ell))~,\] 
and its maximum is achieved on $\partial \mathbb T(\ell)$, see \cite[Thm 4.1.6]{Bu1992}.
In particular we get that each component of the boundary of such a tube $\mathbb T(\ell)$ has length 
\[\ell \cosh\left(\argsinh \left(\f 1 {\sinh(\ell/2)}\right) \right)~,\] and the injectivity radius on points $p\in\partial \mathbb{T}(\ell)$ is given by $\argsinh\left(\sinh(\ell/2)\cosh(L)\right)$.

\subsubsection{Complex Earthquakes along simple closed geodesics}\label{complexearth}
We restrict here the treatment to the case of twisting or earthquaking on a simple closed geodesic, but both operations can be extended to the space of measured laminations (see \cite{McMcomplexearth}).

Let $\gamma$ be a simple closed geodesic on a Riemann surface $X\in\mathcal{T}(S)$. We fix an identification of the universal cover $\widetilde{X}$ with $\mathbb{H}^2$ such that the geodesic $\gamma$ lifts to the vertical geodesic $\tilde{\gamma}$ of the upper half-space $\mathbb{H}^2$ through $0$ and $\infty$. Given a hyperbolic surface $X$ and a simple closed geodesic $\gamma$ the parameter $t$  \textit{earthquake} $\text{eq}_{t, \gamma}(X)$, $t\in\mathbb{R}$, is defined in the lift of a neighborhood of $\gamma$ as 
\[ z \longrightarrow \begin{cases}
	e^{t}z &\text{ if } z\in\mathbb{H}^2_{s} \\
	z &\text{ if } z\in\mathbb{H}^2_{d}
\end{cases}~,\]
with $\mathbb{H}^2_{d}=\{ z\in\mathbb{H}^2\ |\  \text{Re}z >0\}$ and $\mathbb{H}^2_{s}=\{ z\in\mathbb{H}^2\ |\  -\text{Re}z >0\}$, while it is the identity outside. This defines a new hyperbolic metric on $S$, therefore, a new point in $\mathcal{T}(S)$, and also a new Fuchsian projective structure.

Given a hyperbolic surface $X$ and a simple closed geodesic $\gamma$ the parameter $s$ \textit{grafting} $\text{gr}_{s, \gamma}(X)$, $s\in\mathbb{R}^+$, is obtained by cutting along $\gamma$ and inserting a Euclidean tube of height $s$ and core equal to the hyperbolic length $\ell_{\gamma}(X)$ of $\gamma$. The gluing of the Euclidean and the hyperbolic metric, which remain unchanged in $X\setminus\gamma$, produces a new conformal structure, and therefore a new point in Teichm\"uller space by uniformization, yielding a new hyperbolic metric. Let now $\widetilde{\mathcal{A}}(s)$ be an $s$-angular sector in $\mathbb{H}^2$. We equip the inserted Euclidean tube with the complex projective structure 
\[\widetilde{\mathcal{A}}(s)/\langle z\rightarrow e^{\ell_{\gamma}(X)}z\rangle~.\]
This, together with the Fuchsian structure on $X\setminus\gamma$, produces a complex projective structure for $\text{gr}_{s \gamma}(X)$, which is denoted by $\text{Gr}_{s \gamma}(X)$. The grafting operation $\text{Gr}_{s \gamma}(X)$ lifts to the universal cover by cutting along $\tilde{\gamma}$ with the map 
\[ z \longrightarrow \begin{cases}
	e^{is}z &\text{ if } z\in\mathbb{H}^2_{s} \\
	z &\text{ if } z\in\mathbb{H}^2_{d}
\end{cases}~,\]
and inserting $\widetilde{\mathcal{A}}(s)$.

Finally, given a parameter $\lambda=t+is\in\C$, $s>0$, we define the \textit{projection} of the $\lambda$-\textit{complex earthquake} on $\gamma$ in $X$ as the path in Teichm\"uller space given by the composition of the two operations just defined: 
\[\pi\circ\text{Eq}_{\lambda \gamma}(X)= \pi\circ\text{Gr}_{s\gamma}\circ \text{eq}_{t\gamma}(X)=\text{gr}_{s \gamma}(X)\circ\text{eq}_{t\gamma}(X) ~,\]
where $\pi\colon\mathcal{P}(S)\rightarrow \mathcal{T}(S)$ is the forgetting projection.

We call, respectively, the \textit{infinitesimal earthquake} and the \textit{infinitesimal grafting} at $X$ along $\gamma$ the following two Beltrami differentials $\mu, \nu \in T_{X}\mathcal{T}(S)$
\[\mu = \partial_t|_{t=0}\left(\text{eq}_{t \gamma}(X)\right), \qquad \nu = \partial_s|_{s=0}\left(\text{gr}_{s\gamma}(X)\right)~.\]
By McMullen's Theorem $2.10$ in \cite{McMcomplexearth}, the complex earthquake map is holomorphic in $\lambda$ and at any $X$: 
\[\nu=i\mu~.\]

\subsubsection{A Stokes type lemma}
In this subsection we state a key lemma in our work. This is an application of the Stokes Theorem in the setting of the pairing between a holomorphic quadratic differential and the Beltrami differentials $\mu$ and $\nu$ corresponding to, respectively, infinitesimal earthquake and infinitesimal grafting. The first part of the statement is Lemma $5.1$ in \cite{CGS2024}. The proof of the second part follows, by applying the same proof, from the fact that $\nu= i\mu$.

 \blem\label{intbypart}
 Let $X$ be a Riemann surface structure on $S$, with $z$ its complex coordinate. Let $\gamma\colon \mathbb{R}/\ell\mathbb{Z}\rightarrow X$ be a unit-speed parameterization of a simple closed geodesic with respect to the unique hyperbolic representative in $X$, of length $\ell$. Let also $\mu\in T_{X}\mathcal{T}(S)$ be the infinitesimal earthquake along $\gamma$, $\nu\in T_{X}\mathcal{T}(S)$ be the infinitesimal grafting along $\gamma$, and $Q(z)=q(z)dz^2\in T^{*}_{X}\mathcal{T}(S)$ be a holomorphic quadratic differential. Then the following holds:

\[ \re\langle Q, \mu\rangle=\f 12\int_{\R/\ell\Z}\re(Q(i\dot\gamma(t),\dot\gamma(t))) dt~,\] 
and

\[\re \langle Q, \nu\rangle=-\f 12\int_{\R/\ell\Z}\re(Q(\dot\gamma(t),\dot\gamma(t))) dt~.\]\elem

\bpf
We just need to prove the second equality. Let $w$ be the vector field realizing the infinitesimal grafting $\nu$ along the image of $\gamma$, that is $\bar\partial w=\nu$. It is enough to notice that, if $v$ is the vector field realizing the infinitesimal earthquake along the same $\gamma$, then, since the complex earthquake map $Eq_{\lambda\gamma}(\cdot)$ is holomorphic (see Theorem 2.10 in \cite{McMcomplexearth}) and also $gr_{s\gamma}(X)=\pi\circ Eq_{is\gamma}(X)$ and $tw_{t\gamma}(X)=Eq_{t\gamma}(X)$ \[w=iv~.\]  
The proof of the version for infinitesimal earthquakes in \cite{CGS2024} then applies straightforwardly after substituting $v$ with $w=iv$.
\epf

In this paper, we will apply Lemma \ref{intbypart} with $q(z)dz^2=S(f)$ the Schwarzian of a complex projective surface.

\subsection{Hyperbolic 3-manifolds} \label{hypback}

 We now recall some facts about hyperbolic $3$-manifolds that we will use in the present work. For simplicity we assume that all manifolds are convex co-compact.

By the Tameness Theorem \cite{AG2004,CG2006}, a hyperbolic $3$-manifold with finitely generated fundamental group is homeomorphic to the interior of a compact manifold with boundary. A Riemannian $3$-manifold is hyperbolic if it is locally isometric to $\mathbb H^3$, where $\mathbb{H}^3$ is the unique simply connected $3$-manifold of sectional curvature $-1$, the three-dimensional hyperbolic space. 
 Analogously to hyperbolic surfaces (and any hyperbolic $n$-manifold), a hyperbolic 3-manifold $M$ is isometric to $\mathbb H^3/\Gamma$ for $\Gamma<\mathbb{P}SL(2, \mathbb{C})$ a discrete, torsion-free subgroup of the orientation preserving isometries of $\mathbb H^3$.

The action of $\Gamma$ on $\mathbb{H}^3$ can be naturally extended to $\partial\overline{\mathbb{H}}^3=\mathbb{CP}^1$, with $\mathbb{CP}^1$ the complex projective line, but it does not remain proper discontinuous, i.e., the closure of the orbit of a general point $x\in \mathbb{H}^3$ has non-empty \textit{limit set} $\Lambda(\Gamma)$ in $\overline{\mathbb{H}^3}$. The complement of the limit set in $\mathbb{CP}^1$ is called the \textit{domain of discontinuity} and it is denoted by $
\Omega(\Gamma)$. The action of $\Gamma$ on $\Omega(\Gamma)$ is properly discontinuous, so that one can define the \textit{conformal boundary at infinity of $M=\mathbb{H}^3/\Gamma$} as the surface 
\[\partial_{\infty}M=\Omega(\Gamma)/\Gamma~.\] 
Since $\Omega(\Gamma)$ is a proper open subset of $\mathbb{CP}^1$ and $\mathbb{P}SL(2,\mathbb{C})$ is the group of biholomorphisms of $\mathbb{CP}^1$, when $\partial_{\infty}M$ is homeomorphic to a closed surface, it is naturally equipped with a complex projective structure, and, therefore, also a Riemann surface structure. Moreover, each component of $\Omega(\Gamma)$ is a hyperbolic domain. We will denote by $f_M$ the (product of the) developing map(s) of the complex projective structure of (the components of) $\partial_{\infty}M$, where the universal cover of each connected component is identified with $\mathbb{H}^2$.

The smallest non-empty convex subset $C(M)$ of $M$ whose inclusion is a homotopy equivalence is called the \textit{convex core} of $M$. When $C(M)$ is also compact, $M$ is called \textit{convex co-compact}. If $M$ is convex co-compact, its conformal boundary at infinity is homeomorphic to a closed surface $S$, therefore \[[\partial_{\infty}M]\in \mathcal{T}(S)~.\]

Connected components of $\overline{M\setminus C(M)}$, or, more generally, of the complement of a geodesically convex compact subset homotopy equivalent to $M$, are called \textit{ends}. By Tameness an end is homeomorphic to $S_i\times [0, +\infty)$, with $S_i$ a component of $S=\partial_{\infty}M$, and it has infinite hyperbolic volume. The \textit{conformal compactification} $\bar{M}=M\cup \partial_{\infty}M$ of $M$ yields a nice parametrization of the deformation space $CC(M)$ of convex co-compact hyperbolic structures on $M$ up to homotopy, see \cite[Thm 5.1.3.]{Ma2016} and \cite[Thm 5.27]{MT1998}
\[ CC(M)=\quotient{\mathcal T(\partial\overline M)}{T_0(D)},\]
where $T_0(D)\subset MCG(\partial\overline M)$ is the subgroup generated by Dehn twists along compressible curves of $\partial\overline M$ and $\mathcal T(\partial\overline M)$ is the product of the Teichm\"uller spaces of the connected components of $\partial\overline M$.

\subsubsection{The Thurston metric}\label{ssc:thmetric} This paragraph will only be needed in the last part of Section \ref{sc:estimates} on Fourier analysis, which is rather technical and can be ignored at first reading. We will need to compare the Euclidean metric on $\mathbb{C}^*$ and the hyperbolic metric on an annular domain $\mathcal{A} \subset\mathbb {CP}^1$. Due to the different curvatures and the asymmetry of $\mathcal{A}$, these metrics behave quite differently. However, there exists a metric, the \emph{Thurston metric}, that can be compared with both. The Thurston metric is also known as the \textit{projective metric}, as it can be defined for every locally univalent map $f\colon\mathbb{H}^2\rightarrow\mathbb{CP}^1$, and so, for any developing map.

\bdefi
Let $f\colon\mathbb{H}^2\rightarrow\mathbb{CP}^1$ be a locally univalent map. The \textit{Thurston metric} on $\mathbb{H}^2$ associated to $f$ is defined as: \[h_{Th}(z)=\inf_{D}h_D(z)~,\]
where the infimum is taken over all the topological disks $D\subseteq\mathbb{H}^2$ with image $f(D)$ a round disk, and $h_D$ is the hyperbolic metric on $D$. 
\edefi

Note that, if $f$ is the developing map of a complex projective surface $Z$, with $\pi(Z)=X$ and $X\cong\mathbb{H}^2/\Gamma_0$, since for each $\gamma\in\Gamma_0$ we have $f\circ\gamma=\rho(\gamma)\circ f$ for some $\rho(\gamma)\in\mathbb{P}SL(2,\mathbb{C})$, and since M\"obius maps preserve the class of round disks, the Thurston metric descends to a conformal metric on $X$.

Any point $p\in \mathbb{H}^3$ defines a conformal metric on $\partial\bar{\mathbb{H}}^3=\mathbb{CP}^1$ given by the push-forward of the spherical metric on the unit tangent bundle $T^1_p\mathbb{H}^3$ at $p$ through the exponential map at $p$, which is a homeomorphism from $T^1_p\mathbb{H}^3$ to $\mathbb{CP}^1$. This metric is called the \textit{visual metric} from $p$. Given a point $z\in \partial \bar{\mathbb{H}}^3$ the visual metric $h_p(z)$ from $p$ at $z$ coincides with $h_{D(p,z)}(z)$, where $h_{D(p,z)}(z)$ is the hyperbolic metric on the disc ${D(p,z)}$ detected on $\partial \bar{\mathbb{H}}^3=\mathbb{CP}^1$ by the unique geodesic hyperplane $\mathbb{H}^2\subset \mathbb{H}^3$ containing $p$ and perpendicular to the geodesic ray connecting $p$ to $z$. More explicitly, considering the upper-half space model for the hyperbolic 3-space $\mathbb{H}^3=\{(z,x)\in \mathbb{C}\times \mathbb{R}^{+}\}$, given $p=(w,x)\in\mathbb{H}^3$, one can verify that the visual metric from $p$ at $z$ is given by
\[h_p(z)=\dfrac{4x^2}{(\abs{w-z}^2+x^2)^2}|dz|^2.\]
From this formula follows that, for any fixed $z\in \partial \bar{\mathbb{H}}^3$, the set of points $p$ in $\mathbb{H}^3$ such that $h_p(z)$ is constant forms a horosphere centered at $z$. By the work of Thurston (see, for example, \cite[Section $2.2$]{BBB2018}), \cite[Section $4.2$]{dumas-survey}, and Epstein \cite{epsteinSurf}, to any $f$
can be associated a locally convex pleated plane in $\mathbb{H}^3$. Briefly, the construction is as follows. For each $z\in\mathbb{H}^2$, restrict $f$ to a neighborhood where it is injective, and consider the push-forward $f_{*}h_{Th}$ of the relative Thurston metric. Consider then the horosphere pointed at $f(z)$ given by the set of points whose visual metric coincides with $f_{*}h_{Th}(z)$. Finally, take the tangent point between this horosphere and the hyperplane in $\mathbb{H}^3$ projecting on $f(D_z)$, with $D_z$ the $f$-round disk realizing the Thurston metric at $z$, which exists, is unique, and it is called \textit{the maximal disk} for $f$ at $z$ \cite[Lemma $1.2.7$]{KamishimaTan1992}. The disk $D_z$ is indeed a disk containing $z$, with round image, such that the restriction of $f$ to $D_z$ is injective, and maximal with respect to the inclusion. This defines a map from $\mathbb{H}^2$ to $\mathbb{H}^3$, whose image is a locally convex pleated plane.

A case of particular interest is when $f$ is a covering onto $\Omega$ the domain of discontinuity of $\Gamma< \mathbb{P}SL(2, \mathbb{C})$, with $M=\mathbb{H}^3/\Gamma$ a convex co-compact hyperbolic $3$-manifold. In this case the Thurston metric descends to a conformal metric on $X=[\Omega/\Gamma]=[\partial_{\infty}M]\in \mathcal{T}(\partial \overline{M})$ the boundary at infinity of $M$. Moreover, by \cite[Prop 2.13]{BBB2018}, the horosphere obtained by imposing the equality between the visual metric and the push-forward by $f$ of the Thurston metric at $z\in \Omega(\Gamma)\subseteq \partial\bar{\mathbb{H}}^3$ is tangent to the lift to $\HH^3$ of the convex core of $M=\mathbb{H}^3/\Gamma$, and so \[h_{Th}(z)=\dfrac{1}{r_z^2}|dz|^2,\]
where $r_z$ is the radius of the horosphere centered at $z$ and tangent to $C(M)$.
Following the proof of Proposition $2.13$ in \cite{BBB2018}, one can check that, anytime the push-forward of the Thurston metric through (the restriction of) $f$ is well-defined on some sub-domain $U\subseteq\mathbb{H}^2$, and the maximal disks satisfy $f(\cup_{z\in U}(D_z))\subseteq\mathbb{C}^*$, the same argument applies for the smallest closed convex subset $CH(\mathbb{CP}^1\setminus f(\cup_{z\in U}(D_z)))$ of $\mathbb{H}^3$ whose closure in $\mathbb{CP}^1$ is $\mathbb{CP}^1\setminus f(\cup_{z\in U}(D_z))$.

\subsubsection{The renormalized volume} \label{ssc:renormvol}

For a complete introduction to the renormalized volume in the setting of convex co-compact hyperbolic $3$-manifolds, we refer to \cite{S2008} and \cite{compare}. 

The ends of a convex co-compact hyperbolic $3$-manifold $M$ have infinite hyperbolic volume. They can be nevertheless renormalized through equidistant foliations determined by a metric representative $g$ in the conformal boundary at infinity. Considering a geodesically convex compact subset $N\subseteq M$, homotopy equivalent to $M$ and with smooth boundary, we call \textit{W-volume} of $N$ the quantity
\[W(N)=\text{Vol}(N)-\f12\int_{\partial N} H dA\]
with $\text{Vol}(\cdot)$ the hyperbolic volume, $H$ the mean curvature of the boundary $\partial N$, and $dA$ its induced area form. It turns out (see \cite[Lemma 3.6]{compare}) that the $W$-volumes of the geodesically convex compacts $N_r$, defined as the set of points at distance less than or equal to $r$ in $M$ from $N$ satisfy
\[W(N_r)+\pi r\chi(\partial N)=W(N)~,\]
where $\chi(\cdot)$ is the Euler characteristic, and we note that topologically $\partial N \cong \partial \bar{M}$. For any $g\in[\partial_{\infty}M]$, through an explicit construction due to Epstein \cite{epsteinSurf}, it is possible to associate bijectively an exhaustion $\{N_r(g)\}_r$ of the same type as above, and therefore also to $g$ its $W$-volume \[W(M,g)=W(N_r(g))+\pi r \chi(\partial M)~.\] Finally, the renormalized volume is defined as 
 \[V_R(M)=W(M, h)~,\]
with $h\in[\partial_{\infty}M]$ the unique hyperbolic representative. Thanks to the parameterization of the deformation space of convex co-compact hyperbolic structures on $M$ through $\mathcal{T}(\partial \bar{M})$, the renormalized volume is a real-valued function over the Teichm\"uller space.

The renormalized volume satisfies the following differential formula (see \cite{S2008}), which will be crucial for our applications of Section \ref{sc:appl}.

\bthm\label{dVR}
Let $M$ be a convex co-compact hyperbolic $3$-manifold, $S(f_M)$ be the holomorphic quadratic differential given by the Schwarzian derivative of the developing map of $\partial_{\infty}M$, and $\mu\in T_{[\partial_{\infty}M]}\mathcal{T}(\partial\bar{M})$, then 
\[d V_R(\mu) =\text{Re}\langle S(f_M), \mu\rangle=\text{Re} \int_{\partial_{\infty}M} S(f_M) \mu\  .\] 
\ethm

We can then apply Lemma \ref{intbypart} to get a simpler formula for the variation of the renormalized volume in the special case of $\mu$ equal to an infinitesimal earthquake or an infinitesimal grafting.

\subsection{Asymptotic Notation} We will use $f(x)\in O(g(x))$ if either for large $x$ or small $x$, which will be clear from the context, one has that $f(x)\leq C g(x)$ for some positive constant $C$.

\section{Toy model: the symmetric tube case} \label{sec:symmsch}

In this section we outline the idea behind Theorem \ref{mainthm}. The main observation is that the result is almost trivial in the case of a \textit{symmetric} complex projective tube.
Here a complex projective tube is symmetric if its image under its developing map in $\mathbb {CP}^1$ is bounded by two disjoint round circles $C_1$ and $C_2$. The annular domain bounded by $C_1$ and $C_2$ is a \emph{symmetric tube}. Note that up to a M\"obius map this configuration can be made so that they are concentric, and with center at $0\in\mathbb{C}\subseteq\mathbb{CP}^1$. Moreover, composing by a M\"obius map does not change the Schwarzian derivative.

 Then, the bulk of the remaining work in proving Theorem \ref{mainthm} is to show that the Schwarzian of a general long complex projective tube behaves like the one in the symmetric case, up to a term that goes to zero in $\ell$, with $\ell$ the hyperbolic length of the core curve.
 
\vspace{0.3cm}

Consider the half-space model for the $2$-hyperbolic space with coordinate $z=\rho e^{i\theta}$, $\rho>0$, $\theta \in [0, \pi]$, and $f_{\ell}\colon \mathbb{H}^2\rightarrow \mathbb{C}$ the map uniformizing a symmetric tube centered at $0$ in $\mathbb{C}^*$ with a hyperbolic metric of core length $\ell$ defined as
\[f_{\ell}(z)=z^{\f{2\pi i}{\ell}}=e^{-\f{2\pi\theta}{\ell}+i\f{2\pi}{\ell}\log(\rho)}.\] 
We denote by $\widetilde{\mathcal{A}}(L)$ the $L$-neighborhood in $\mathbb{H}^2$ of the vertical geodesic between $0$ and $\infty$.
Since $\widetilde{\mathcal{A}}(L)$ corresponds to the annular sector $\theta\in[\theta(L), \pi-\theta(L)]$ with (see \cite[Lemma $5.2.7$]{Mar16}) \[\f{1}{\cosh(L)}=\sin(\theta(L))~,\] the restriction of $f_{\ell}$ to $\widetilde{\mathcal{A}}(L)$ is the developing map of the annulus 
\[\mathcal{A}_{\ell}(L)=\left\{r e^{i\theta}\in\mathbb{C}\ \ |\ \ e^{-\f{2\pi^2}{\ell}+\f{2\pi}{\ell}\arcsin\left(\f{1}{\cosh(L)}\right)} \leq r \leq e^{-\f{2\pi}{\ell}\arcsin\left(\f{1}{\cosh(L)}\right)} \right\}~.\]
Note that this is indeed \textit{symmetric}.
\brem\label{L}
For any symmetric complex projective tube $\mathcal{A}$, up to composing by a M\"obius transformation, there exist $\ell>0$ and $L>0$ such that $\mathcal{A}$ is equal to $\mathcal{A}_{\ell}(L)$. We also recall that composing by M\"obius transformations leaves the Schwarzian derivative invariant.
\erem

Observe that, if we equip $\mathcal{A}_{\ell}(L)$ with the restriction of the flat metric $h_0= \f 1{r^2}|dz|^2$ on the infinite tube $\C^*$, we get a well-truncated Euclidean tube with core length $2\pi$ and of some height $2m$, with $2m$ such that $(\mathcal{A}_{\ell}(L), h_0)$ is conformal to a truncated hyperbolic tube with core length $\ell$ and height $2L$.

\brem\label{modulus}
	The conformal structure of a (truncated) tube is uniquely determined by its \textit{modulus}, which, in the case of a symmetric annulus in $\mathbb{C}$ with boundaries centered at the origin and of radii $r_1$ and $r_2$, is equal to $\f{1}{2\pi}\log\left(r_2/r_1\right)$. 
	Then, any (truncated) tube equipped with a Riemann surface structure is conformal to a symmetric one (but in general not projectively equivalent to it).  
\erem

By a direct computation we have:
\blem\label{symmlemma}
Let $\mathcal{A}$ be a symmetric complex projective tube. Then, there exist $\ell>0$ and $L>0$ such that the Schwarzian of $\mathcal{A}$ identifying the universal cover $\widetilde{\mathcal{A}}$ with $\widetilde{\mathcal{A}}(L)$ is equal to
\[S(f_{\ell})(z)=\f{1}{2z^2}\left(1+\f{4\pi^2}{\ell^2}\right)dz^2~.\]
   
\elem

\bpf
By remark \ref{L}, we can assume $\mathcal{A}$ to be equal to $\mathcal{A}_{\ell}(L)$, for some $\ell>0$ and $L>0$. Then, we just have to compute the Schwarzian derivative of its developing map $f_{\ell}$ on $\widetilde{\mathcal A}(L)$:
\begin{align*}
q=S(f_{\ell})&=\biggl(\biggl(\dfrac{(f_{\ell})^{''}}{(f_{\ell})^{'}}\biggr)^{'}-\dfrac{1}{2}\biggl(\dfrac{(f_{\ell}){''}}{(f_{\ell})^{'}}\biggr)^{2}\biggr)dz^2 \\
&=\left(\left(\left(\f{2\pi i}{\ell}-1\right)\f1z\right)^{'}-\f12  \left(\left(\f{2\pi i}{\ell}-1\right)\f1z\right)^{2}\right)dz^2 \\
&= \f{1}{2z^2}\left(1+\f{4\pi^2}{\ell^2}\right)dz^2 ~.
\end{align*}
This completes the proof. \epf 

 In Section \ref{sc:estimates}, the width $L$ will usually be the one given by the Collar Lemma (see Subsection \ref{ssc:margtubes}): $L=\operatorname{arcsinh}\left(\f{1}{\sinh(\ell/2)}\right)$, so that the neighborhoods $\widetilde{\mathcal{A}}$ considered are \emph{wide}. In this case, we denote the symmetric tube simply by $\mathcal{A}_{\ell}$, and we say that $(\mathcal{A}_{\ell}, f_{\ell})$ is the \textit{standard symmetric tube} of core length $\ell$.

\brem
The expression of $S(f_{\ell})$ in Lemma \ref{symmlemma} is in $\mathbb H^2$, however, if one wants to use the coordinates of $\mathcal{A}$ rather than viewing it as a neighborhood of the vertical axis in $\mathbb{H}^2$, the formula changes as follows.
We first need to recall that the Schwarzian derivative satisfies the following composition rule:  
\begin{equation}\label{composition}
	S(f\circ g)= g^{*}S(f)+S(g)~,
\end{equation}
whenever $f$ and $g$ are two locally injective holomorphic maps, whose composition $f\circ g$ is well-defined.
Denoting by $w$ the complex coordinate in $\mathbb{C}$, a local inverse $f_{\ell}^{-1}$ of $f_{\ell}$ sends $w$ to $z=w^{-\f{i\ell}{2\pi}}$, and then also $dz(w)^2=-\f{\ell^2}{4\pi^2}w^{2\left(-\f{i\ell}{2\pi}-1\right)}dw^2$

Therefore, by applying equation \eqref{composition} with $f=f_{\ell}$ and $g=f_{\ell}^{-1}$ we obtain: 
\[S\left(f_{\ell}^{-1}\right)=-\left(f_{\ell}^{-1}\right)^*(S(f_{\ell}))\]
and then
\[S\left(f_{\ell}^{-1}\right)= \f{1}{2w^2}\left(1+\f{\ell^2}{4\pi^2}\right)dw^2~.\]

\erem

Lemma \ref{symmlemma} applied to the standard symmetric tubes $(\mathcal{A}_{\ell}, f_{\ell})$ is the ``toy model" version of Theorem \ref{Schwarzian}. The next proposition is the ``toy model" version of Theorem \ref{pairing}. This is because, thanks to the key Lemma \ref{intbypart}, the two integral terms coincide with the pairing of the Schwarzian with, respectively, infinitesimal earthquakes and infinitesimal graftings.  Then, in Section \ref{sc:estimates} we show that the result in the general case is equal to the formulas in Lemma \ref{symmlemma} and Proposition \ref{symmetric} up to an error term decaying exponentially to zero in $\ell$.

\bprop\label{symmetric}
On the core $\gamma$ of the symmetric complex projective tube $\mathcal{A}_{\ell}$, the following equalities for the Schwarzian of $f_{\ell}^{-1}$ hold: \[ \re\left(S\left(f_{\ell}^{-1}\right)(i\dot\gamma, \dot\gamma)\right)= 0~,\qquad \re\left(S\left(f_{\ell}^{-1}\right)(\dot\gamma, \dot\gamma)\right)=\f{2\pi^2}{\ell^2}+\f12~.\]
Moreover
\[\f12\int_0^{\ell}\re(S\left(f_{\ell}^{-1}\right)(i\dot\gamma, \dot\gamma))dt= 0~,\qquad\f12 \int_0^{\ell}\re(S\left(f_{\ell}^{-1}\right)(\dot\gamma, \dot\gamma))dt=\f{\pi^2}{\ell}+\f{\ell}{4}~.\]
\eprop

\bpf
Up to lifting $\mathcal{A}_{\ell}$ with $f_{\ell}$, we have $\gamma(t)=ie^t$, $t\in [0, \ell]$, and $S(f_{\ell})$ as in Lemma \ref{symmlemma}. Then, $\dot\gamma(t)=e^t\f{\partial}{\partial y}$, $i\dot\gamma(t)=-e^t\f{\partial}{\partial x}$, and, since $dz=dx+idy$, on $\gamma$, we have 
\[\re \left(\f{1}{2z^2}\left(1+\f{4\pi^2}{\ell^2}\right)dz^2\right)(i\dot\gamma, \dot\gamma)=\re \left(-\f{1}{2e^{2t}}\left(1+\f{4\pi^2}{\ell^2}\right)ie^{2t}\right)=0~,\]
and 
\[\re \left(\f{1}{2z^2}\left(1+\f{4\pi^2}{\ell^2}\right)dz^2\right)(\dot\gamma, \dot\gamma)=\re \left(\f{1}{2e^{2t}}\left(1+\f{4\pi^2}{\ell^2}\right)e^{2t}\right)=\f12 \left(1+\f{4\pi^2}{\ell^2}\right)~,\]
completing the proof of the first part of the statement. The second one follows directly by integration. \epf

\section{Schwarzian derivative on long tubes}
\label{sc:estimates}

This section is dedicated to the study of the Schwarzian on long complex projective tubes, which we now recall:
\bdefi\label{defilongprojtube}
	Let $Z=(f,\rho)$ be a complex projective structure on $S$, and $X=\pi(Z)$ be the induced conformal hyperbolic surface. Let $\gamma\subset X$ be a short simple closed geodesic. We say that $Z$ contains a nontrivial \textit{long tube} with core $\gamma$ if, for any lift $\widetilde\gamma\subset \widetilde X$, up to post-composition by a M\"obius map, the restriction of the developing map $f$ to the wide neighborhood $\widetilde{\mathcal{A}}$ of $\widetilde{\gamma}$ is a covering onto an essential annulus $\mathcal{A}\subseteq\mathbb{C}^*\subseteq \mathbb{CP}^1$ with deck group $\langle\gamma\rangle$, and, $\forall z\in \widetilde{\mathcal A}$, the maximal disk $D_z$ satisfies $f(D_z)\subset\C^*$. We call $\mathcal{A}$ a \textit{long complex projective tube}.
\edefi

 In the first subsection, we start by introducing the central object used in our analysis: the \textit{Osgood-Stowe tensor of a complex projective structure}. We will denote it by $B(\rho, e^{2u}\rho)$, where $\rho$ stands for the hyperbolic metric of $\mathbb{H}^2$, and $u\colon\mathbb{H}^2\rightarrow \mathbb{R}$ is a smooth function which is the \textit{logarithmic conformal factor} between the two metrics $\rho$ and $e^{2u}\rho$, and depends on the projective structure. Right after, we prove Theorem \ref{Schwarzian} exploiting the estimates of Subsection \ref{ssc:fourier}.
 
We then define two suitable norms for $B(\rho, e^{2u}\rho)$, which coincide, respectively, with the real part of the pairing between the Schwarzian with infinitesimal earthquakes and infinitesimal graftings (see Definitions \ref{definB} and \ref{definBgr}, and Lemma \ref{intbypart}). Then, in Subsections \ref{ssc:estimatesearth} and \ref{sec:graft} we give bounds for the pairing along earthquakes and grafting respectively. Finally, we prove Theorem \ref{pairing} at the end of Subsection \ref{sec:graft}.

 The main idea behind these estimates is to express $u$ as a sum of three other logarithmic conformal factors of different metrics, such that one of them, denoted by $u_1$, satisfies $\Delta u_1 = 0$. We will then require the Fourier analysis of Subsection \ref{ssc:fourier}, where we bound the \textit{flat logarithmic conformal factor $u_1$} and its derivatives in terms of the hyperbolic length of the long tube.

\subsection{The Osgood-Stowe tensor}\label{Schwarziantensor}

We introduce here the main tool of our analysis: the \textit{Osgood-Stowe tensor} presented in \cite{OS1992} by the same authors, and whose relation with the $3$-dimensional hyperbolic geometry has been largely studied by Bridgeman and Bromberg in \cite{Epstein-wvolume-OS}. The Osgood-Stowe tensor $B$ associates to any pair of conformal metrics $g$ and $g'=e^{2u}g$ on a Riemann surface $X$, with $u\colon X\rightarrow \mathbb{R}$ a smooth function, a real $(0,2)$-tensor, defined through the logarithmic conformal factor $u$ between the two metrics (see Definition \ref{tensor}). We will denote it by $B(g, e^{2u}g)$.
The key point is then that the real part of the Schwarzian  $\re\left(S(f)\right)$ of a long tube in a complex projective structure $Z$ (i.e. $f\colon \widetilde{\mathcal{A}}\subseteq\mathbb{H}^2\rightarrow \mathbb{CP}^1$ is the restriction of the developing map of $Z$ on $S$ covering a long complex projective tube $\mathcal{A}$) is equal to $\B{\rho} u$, with $\rho$ the hyperbolic metric on $\mathbb{H}^2$, and $u$ such that \[f^{*}(|dz|^2)=e^{2u}\rho~,\] with $z$ the complex coordinate of $\mathbb{C}\subseteq \mathbb{CP}^1$ (see Theorem \ref{tensorschw}). We then say that $\B{\rho} u$ is the \textit{Osgood-Stowe tensor of the long complex projective tube} of $Z$.

From now on, we will then focus on what happens on the restriction of the developing map $f$ of a complex projective structure $Z$ on a closed surface $S$ covering a long complex projective tube, which with some abuse of notation we will still often denote by $f$. More specifically, this means that the image $\mathcal{A}$ of $f$ is an annulus and that its \textit{core} geodesic $\gamma$, in the induced hyperbolic metric $h_{\infty}$ pushed forward on $\mathcal{A}$ by $f$, has length $\ell\leq\varepsilon_0$, see Section \ref{hypmanifolds}. The underlying Riemann structure $X=\pi(Z)$, where, recall, $\pi\colon \mathcal{P}(S)\rightarrow\mathcal{T}(S)$ is the vector bundle of projective structures on a closed surface $S$, defines a natural identification of the universal cover $\widetilde{S}$ with $\mathbb{H}^2$. There are then two subgroups $\Gamma_0<\mathbb{P}SL(2, \mathbb{R})$ (isomorphic to $\pi_1(S)$) and $\Gamma'< \mathbb{P}SL(2, \mathbb{C})$ (isomorphic to $\pi_1(S)$/$K$, with $K$ the kernel of the holonomy $\rho$ of $Z$) such that $f$ is $(\Gamma_0, \Gamma')$-equivariant.
Up to composing by a M\"obius transformation, which does not affect the Schwarzian of $f$, we can assume $f^{-1}(\gamma)$ to be the geodesic through $0$ and $\infty$ in $\mathbb{H}^2$. With an abuse of notation, we will also denote by $\gamma$ its projection on $X$, and its lift $f^{-1}(\gamma)$ in $\mathbb{H}^2$. Equipping $X$ with its unique hyperbolic metric, we consider the \textit{thin tube} $\mathbb{T(\ell)}\subseteq X$, i.e. the $L$-neighborhood of $\gamma\subseteq X$, with $L=\argsinh\left(\f{1}{\sinh\left(\f \ell 2\right)}\right)$ (see Definition \ref{margulitubedefi}). This lifts to the wide $L$-neighborhood $\widetilde{\mathcal{A}}\eqdef \widetilde{\mathcal{A}}(L)\subseteq\mathbb{H}^2$ of $f^{-1}(\gamma)$. Then, by our hypothesis, its image $\mathcal{A}=f(\widetilde{\mathcal{A}})$ through $f$ is a long complex projective tube with core $\gamma$. Thus, we have the following diagram:
\[ \xymatrix{ \widetilde{ \mathcal A}\ar[r]^ \subset \ar[d]_{/\langle\gamma\rangle}& \widetilde{ S} \ar[d]^{/\langle\gamma\rangle} \\
	\widetilde{\mathcal{A}}/\langle\gamma\rangle\ar[r]^\subset \ar[d]_{/\Gamma}& \widetilde{S}/\langle\gamma\rangle\ar[d]^{/\Gamma}\\
	\mathbb T\ar[r]^\subset & X } \]
with $\Gamma$ isomorphic to $\pi_1(S)/\langle\langle\gamma\rangle\rangle$, where we recall we are always identifying $\widetilde{S}$ with $\mathbb{H}^2$. Note also that $\mathcal{A}$ equipped with the push-forward of the hyperbolic metric on $\widetilde{\mathcal{A}}$ via $f$ is isometric to $\widetilde{\mathcal{A}}/\langle\gamma\rangle$.

We now remind the reader that our main goal is to study $S(f)$ on $\widetilde{\mathcal{A}}$. To this aim, we will widely use the Osgood-Stowe tensor, which we now formally define (see \cite{OS1992}).

\bdefi \label{tensor}
Let $(X, g)$ and $(X', g')$ be two Riemannian surfaces, and $f\colon (X,g)\rightarrow (X',g')$ be a conformal local diffeomorphism. If $u\colon X\rightarrow \mathbb{R}$ is the logarithmic conformal factor such that $f^{*}g'=e^{2u}g$, then the Osgood-Stowe tensor $B$ of $u$ is defined as 
\[ B\left( g,  f^* g' \right) =\B{g} {u}=\text{Hess}(u)-du\otimes du-\dfrac{1}{2}(\Delta u-||\nabla u||^2)g\]
where the hessian, the laplacian, and the gradient are with respect to the metric $g$.
\edefi

The Osgood-Stowe tensor can be regarded as a generalization of the Schwarzian derivative (see Theorem \ref{tm:standard}) and satisfies analogous properties. In particular, given  any two conformal diffeomorphisms $f_1\colon (X,g)\rightarrow (X_1, g_1)$ and $f_2\colon (X_1, g_1)\rightarrow (X_2, g_2)$, property \eqref{composition} extends to additivity under composition 
	\begin{equation}\label{comprule}B\left(g, (f_2\circ f_1)^*g_2\right)=B\left(g, f_1^*g_1\right)+f_1^*B\left(g_1, f_2^*g_2\right)~,\end{equation}
 while the vanishing of the Schwarzian of M\"obius maps translates to the Osgood-Stowe tensor being zero when $f$ is a local isometry, see Appendix \ref{standardmetrics}, \cite{OS1992} and \cite{Epstein-wvolume-OS}.

 We define $h_{\infty}=f_{*}\rho$, with $\rho$ the restriction of the hyperbolic metric of $\mathbb{H}^2$ to $\widetilde{\mathcal A}$, so that $f\colon \left( \widetilde{\mathcal{A}},\rho\right) \rightarrow \left( \mathcal{A}, h_\infty\right)$ is a local isometry. Let us fix $z=x+iy$ and $|dz|^2=dx^2+dy^2$, respectively, the standard complex coordinate and the standard metric on $\mathbb{C}$ contained in $\mathbb{CP}^1$ through stereographic projection.

\bthm \label{tensorschw}
Let $f\colon (\widetilde{\mathcal{A}}, \rho)\rightarrow (\mathcal{A}, |dz|^2)$ be the developing map in the notation above, and $u\colon \widetilde{\mathcal A}\rightarrow \mathbb{R}$ be the logarithmic conformal factor such that $f^{*}|dz|^2=e^{2u}\rho$, then 
\[\re(S(f))=\B{\rho} u~.\]
\ethm

\bpf
This follows from Theorem \ref{tm:standard} in the Appendix, or Section $4$ of \cite{Epstein-wvolume-OS}. 
\epf

As we are not able to compute $\B {\rho} u$ directly, we split the Osgood-Stowe tensor in a sum of terms, and use the additivity of the Osgood-Stowe tensor \eqref{comprule} to compute $f_*\B {\rho} u$ on $\mathcal{A}$. We will use the following notation for the auxiliary metrics:
\begin{itemize}
\item $h_{\infty,0} = e^{2u_0} h_\infty$ where $h_{\infty,0}$ is the conformal flat metric on $\mathcal A$ which has the same rotational invariance as $h_{\infty}$, i.e. such that there exist $(\rho, \theta)\in [-L, L]\times \mathbb{S}^1$ and flat coordinates $(r, \theta)\in [-m, m]\times \mathbb{S}^1$ on $\mathcal{A}$ such that 
\[ dr=\dfrac{2\pi}{\ell\cosh{\rho}}d\rho~, \qquad m=\int_0^{L}\f{2\pi}{\ell\cosh(\rho)}d\rho~,\] 
while $h_\infty$ is written as
\[ h_{\infty}=d\rho^2+\left(\dfrac{\ell}{2\pi}\right)^2\cosh^2(\rho)d\theta^2~, \]
and \[h_{\infty, 0}=dr^2+d\theta^2~;\] 
in particular, in these coordinates\[u_0(\rho, \theta)= \log\left(\f{2\pi}{\ell\cosh(\rho)}\right)~.\]
\textbf{Convention.} 
	In both of the above metrics by the \emph{core} we will mean the geodesic representative of the non-trivial loop, and by its \emph{width} we mean $2L$ or $2m$.
		
\item $h_{0} = e^{2u_1} h_{\infty,0}$ is another conformal flat metric on $\mathcal A$, given by the restriction of the flat Euclidean metric $h_0=\f{\abs{dz}^2}{r^2}$ on $\mathbb{C}^*$, where $r^2=z\bar z$.\end{itemize}

Thus, we consider the following locally conformal maps:
 \[\xymatrix{
\left(\A, \rho \right) \ar[rrrr]^{f\vert_{ \A}}_u \ar[dr]^{f\vert_{ \A}}_v& & & & \left( \mathcal A,\abs{dz}^2\right)\\
& \left( \mathcal A, h_\infty\right)\ar[r]^{\id_0}_{u_0} & \left(\mathcal A, h_{\infty, 0}\right)\ar[r]^{\id_1}_{u_1} &\left(\mathcal{A}, h_0\right) \ar[ur]^{\id_2}_{w}&
}\]

\begin{figure}[htb!]
\def\svgwidth{\textwidth}
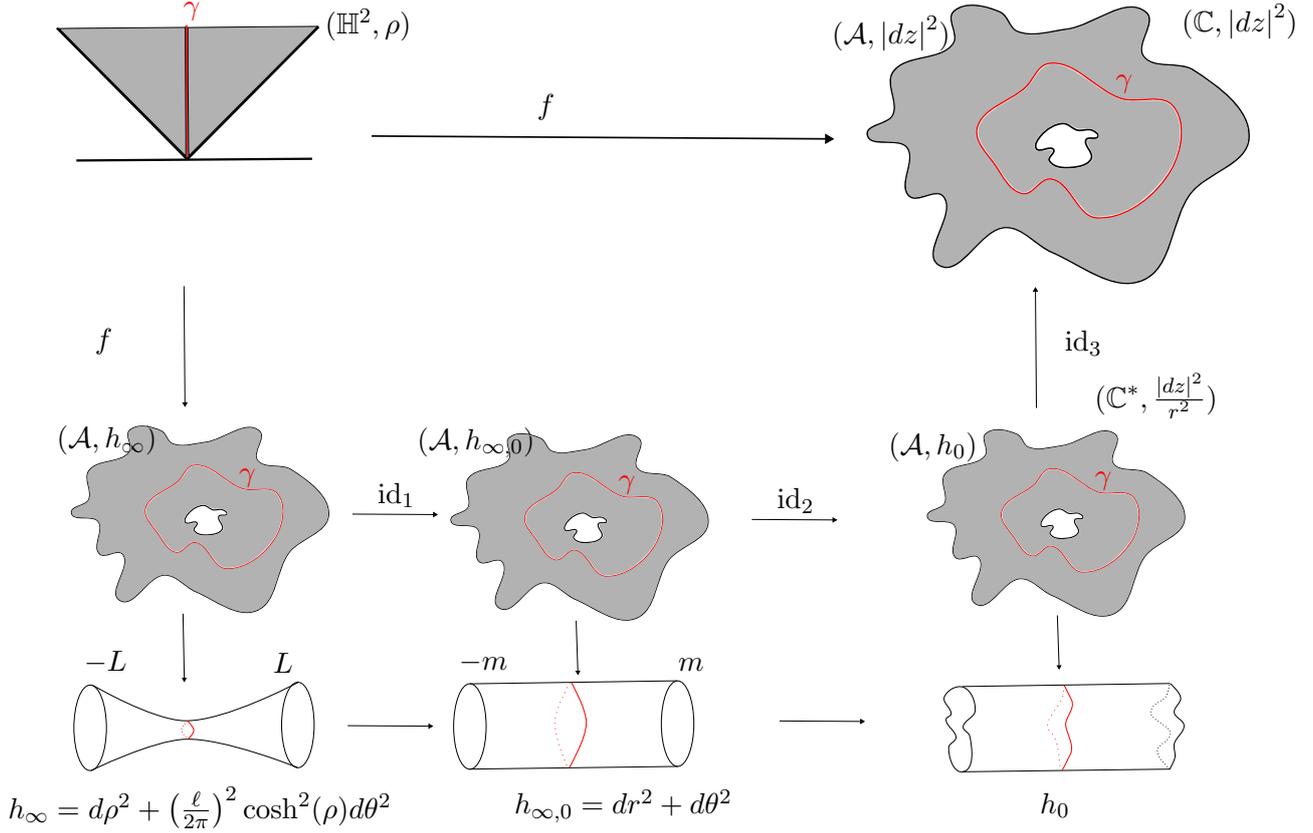
\caption{The long tube in the universal covers with the respective metrics.}\label{cylpic}\label{PicCylinders}
\end{figure}

See Figure \ref{cylpic} for sketches of the domains.

Clearly, we have that:
\[f\vert_{\A} =\id_2\circ \id_1\circ \id_0\circ f\vert_{\A}~.\]
By additivity of the Osgood-Stowe tensor 
\[ \B{\rho}{u\vert_{\A}}=\B\rho {v\vert_{\A}}+f^{*}\B{h_{\infty}}{u_0\vert_{\mathcal A}}+f^*\B{h_{\infty, 0}}{u_1\vert_{\mathcal A}}+f^*\B{h_0}{w\vert_{\mathcal A}},\]
where $v=0$ is the logarithmic conformal factor of a local isometry, so that its contribution is zero, thus
\begin{equation}
  \label{eq:sum}
f_*\B{\rho}{u\vert_{\A}}=\B{h_{\infty}}{u_0\vert_{\mathcal A}}+\B{h_{\infty, 0}}{u_1\vert_{\mathcal A}}+\B{h_0}{w\vert_{\mathcal A}}~,
\end{equation}
where the last term can be computed using that $h_0=\dfrac{\abs{dz}^2}{r^2}$, with $r$ the radial coordinate of $\mathbb{C}^*$. Moreover, thanks to the key Lemma \ref{intbypart}, to prove Theorem \ref{pairing}, we actually only need these terms on the core $\gamma$.

The following remark, together with the proof of Theorem \ref{Schwarzian} which follows afterwards, formalizes the heuristic of why the toy model of Section \ref{sec:symmsch} is a good model also for general long complex projective tubes.

\brem\label{RemSymm}
If $\mathcal{A}$ is symmetric (as in Section \ref{sec:symmsch}), the domains in Figure \ref{cylpic} all have round boundaries. Moreover, up to M\"obius map, $u_1=0$, while $u_0$ and $w$ remain unchanged, as also the respective Osgood-Stowe tensors.

We now discuss the case of a general long complex projective tube $\mathcal{A}$ of core length $\ell\leq\varepsilon_0$ and width $L=\operatorname{arcsinh}\left(\f{1}{\sinh(\ell/2)}\right)$ (with respect to the underlying hyperbolic structure). In the same notations and definitions of Section \ref{sec:symmsch}, where we denoted by $\mathcal{A}_{\ell}$ the standard symmetric tube of core length $\ell$, and by $f_{\ell}$ its developing map from $\widetilde{\mathcal{A}}$, there exists a conformal map \[F_{\ell}\colon \mathcal{A}_{\ell}\rightarrow \mathcal{A}\] such that \[f|_{\widetilde{\mathcal{A}}}=F_{\ell}\circ f_{\ell}\] and $F_{\ell}$ is a holomorphic isometry between $(\mathcal{A}_{\ell}, h_0)$ and $(\mathcal{A}, h_{\infty,0})$, i.e. \[F_{\ell}^*h_{\infty, 0}=h_0~.\]
The existence of a biholomorphism $F_{\ell}$ is guaranteed by the fact that, being the conformal structure of a (truncated) tube uniquely determined by its modulus, the underlying Riemann structure of a complex projective tube is conformal to the one of a symmetric tube with same core length and modulus, and, therefore, same length (see Remark \ref{modulus}). Then, the composition $F_{\ell}\circ f_{\ell}$ gives a developing map on $\widetilde{\mathcal{A}}$ for $\mathcal{A}$, which has to coincide, up to M\"obius transformation, with $f|_{\widetilde{\mathcal{A}}}$. We choose $F_{\ell}$ such that the equality between these maps holds. Finally, denoting by $\rho$ the hyperbolic metric on $\mathbb{H}^2$ \[F_{\ell}^{*}h_{\infty, 0}=e^{2u_0\circ F_{\ell}}F_{\ell}^{*}(f_{*}\rho)=e^{2u_0\circ F_{\ell}}F_{\ell}^*((F_{\ell}\circ f_{\ell})_*\rho)=e^{2u_0\circ F_{\ell}}(f_{\ell})_*\rho=h_0~.\]
\erem

 We now prove one of the main results of this paper:

 \begin{customthm}{\ref{Schwarzian}}Let $\mathcal{A}$ be a long tube of a complex projective surface $Z$, and let $\gamma$ be its core of length $\ell\leq \varepsilon_0$. Let also $f\colon \mathbb{H}^2\rightarrow \mathbb{CP}^1$ be the developing map of $Z$, and let $S(f)$ be its Schwarzian derivative. Then, in a tubular neighborhood of the lift $\widetilde{\gamma}=i\mathbb{R}_{+}$, the Schwarzian $S(f)$ behaves as follows

	 \[S(f)=\f{1}{2z^2}\left(1+\f{4\pi^2}{\ell^2}\right)dz^2+q_\ell(z)dz^2~,\]
	where $z$ is the complex coordinate of $\mathbb{H}^2$, and $\norm{q_\ell(z)dz^2}_{\infty}$ is bounded and in $O\left(\f{e^{-\pi^2/(2\ell)}}{\ell^2}\right)$.
 
\end{customthm}

\begin{proof}
	 We start by showing how the Osgood-Stowe tensor term $B(h_{\infty, 0}, e^{2u_1|_{\mathcal{A}}}h_{\infty, 0})$  is related to the real part of the Schwarzian derivative of $F_{\ell}$, with $F_{\ell}$ as in Remark \ref{RemSymm}. We show that they coincide up to pull-back and adding twice the Schwarzian uniformization for the infinite tube. By Theorem \ref{tm:standard} in the Appendix 
	 \[\re\left(S(F_{\ell})\right)=B\left(|dz|^2, F_{\ell}^*|dz|^2\right)~.\]
	By using the properties of the Osgood-Stowe tensor
	\begin{align*}
		B\left(|dz|^2, F_{\ell}^*|dz|^2\right)&=B\left(|dz|^2, h_0\right)+B\left(h_0, F_{\ell}^*h_0\right)+B\left(F_{\ell}^*h_0, F_{\ell}^*|dz|^2\right)\\
		&=B\left(|dz|^2, h_0\right)+F_{\ell}^*B\left(h_{\infty,0}, h_0\right)-F_{\ell}^*B\left(|dz|^2, h_0\right)
	\end{align*}
	where for the central term of the second equality we used that $(F_{\ell})^*h_{\infty, 0}=h_0$ (as seen in Remark \ref{RemSymm}). Let us also recall that $h_0=e^{2u_1}h_{\infty, 0}$. 
	Then, by the properties of the Schwarzian derivative and the relation just obtained: 
	\begin{align*} f_{*}S(f)&=f_{*}S(F_{\ell}\circ f_{\ell})\\
	&=f_*\left( f_\ell^*S(F_{\ell})+S(f_{\ell})\right)\\
	&=(F_\ell)_*S(F_{\ell})+f_*S(f_{\ell})\end{align*}
	with \begin{equation}\label{Fl}
		(F_{\ell})_{*}\re\left(S(F_{\ell})\right)=(F_{\ell})_{*}B\left(|dz|^2, h_0\right)+B\left(h_{\infty,0}, h_0\right)-B\left(|dz|^2, h_0\right)~.\end{equation}
		We are now going to estimate the terms in the right hand side of equation \eqref{Fl}. In what follows we will use the following abuse of notation. With $Q_2(f(\ell))$ we will mean a $(0,2)$-tensor (in the base $dx^2$, $dxdy$, $dy^2$, with $z=x+iy$) whose coefficients in absolute value are bounded and in $O(f(\ell))$ for small $\ell$. We will write various equalities involving quadratic differentials and, as we only care about their decay, we will write equations of the form $B(\cdot, \cdot)=Q_2(f(\ell))$, where the right hand side is to be interpreted as some quadratic differential in the class $Q_2(f(\ell))$. For a more precise version see Remark \ref{sharperboundschwarzian}.

		First, by using Definition \ref{tensor}, or equivalently by Remark \ref{B3remark}, we compute the term $B(\abs{dz}^2, h_0)$ with $h_0=\abs{dz}^2/r^2$:
		\[B(|dz|^2, h_0)=\f{1}{2r^2}dr^2-\f12d\theta^2=\re\left(\f{1}{2z^2}dz^2\right)~,\]
		where $z=re^{i\theta}$ is the complex coordinate of $\mathbb{C}$. Under the isometric change of coordinates $w=\log(z)$ from $(\mathbb{C}^{*}, h_0)$ to the horizontal infinite strip of height $2\pi$ in $(\mathbb{C}, w)$, the Osgood-Stowe tensor $B(\abs{dz}^2, h_0)$ becomes 
		\[B\left(e^{2\re (w)}\abs{dw}^2, h_0\right)=B\left(e^{2\re(w)}\abs{dw}^2, |dw|^2\right)=\re\left(\f{1}{2}dw^2\right)~.\]
		 Secondly, we note that, since $F_{\ell}$ is such that $(F_{\ell})^*h_{\infty, 0}=h_0$, with $h_0=e^{2u_1}h_{\infty, 0}$, the norm of its complex derivative is 
	 \[\norm{d F_{\ell}}=e^{u_1}~,\]
	 	where, by Lemma \ref{estcore} and Remark \ref{utube}, the function $e^{u_1}$ on a neighborhood of the core of $\mathcal{A}$ is such that
	\[e^{u_1(z)}= 1+ E_\ell(z) ~,\]
	where $E_\ell(z)$ is a bounded function in $O\left(e^{-\pi^2/(2\ell)}\right)$. Note also that $\norm{d F_{\ell}}$ in the $w$ coordinate is the complex norm of $F_{\ell}'$.
	Then, 
	\begin{align*} 
	(F_{\ell})_{*}B\left(e^{2\re(w)}\abs{dw}^2, |dw|^2\right)-B\left(e^{2\re(w)}\abs{dw}^2, |dw|^2\right)&=\re\left(\f{dw^2}{2\left(F_{\ell}^{'}(F_\ell^{-1}(w)\right)^2}-\f{dw^2}{2}\right)~,
	\end{align*}
	and we can write $F_{\ell}^{'}(F_\ell^{-1}(w))=\abs{F_{\ell}^{'}(F_\ell^{-1}(w))}e^{i\nu(w)}=e^{u_1(w)+i\nu(w)}$, where, by Cauchy-Riemann $|d\nu|=|du_1|\in O\left(e^{-\pi^2/(2\ell)}\right)$.
	
Choose the base point on the image core where \(\re (w)=\log|z|\) attains its maximum. Its tangent is vertical there, as is the tangent to the source core in its flat coordinates. Hence the complex derivative of \(F_\ell\), expressed in these flat coordinates, is real and nonzero at that point. Translating the target flat coordinate, we may take this point to be \(w=0\). Thus \(\nu(0)\in\pi\mathbb Z\), so \(e^{-2i\nu(0)}=1\). Then, the bound on \(d\nu\), together with the bound on \(u_1\), gives
\[
e^{-2u_1(w)-2i\nu(w)}-1
=O\!\left(e^{-\pi^2/(2\ell)}\right)
\]
on every fixed smaller tube.

	Therefore  
	\[(F_{\ell})_{*}B\left(e^{2\re(w)}\abs{dw}^2, |dw|^2\right)-B\left(e^{2\re(w)}\abs{dw}^2, |dw|^2\right)=\re\left(\f{1}{2}\left(e^{-2u_1(w)-2i\nu(w)}-1\right)dw^2\right)\]
	is in $Q_2\left(e^{-\pi^2/(2\ell)}\right)$, and then also
	\begin{equation}\label{Bz}
			(F_{\ell})_{*}B\left(|dz|^2, h_0\right)-B\left(|dz|^2, h_0\right)=Q_2\left(e^{-\pi^2/(2\ell)}\right) ~.
	\end{equation} 
	Concerning the $B\left(h_{\infty,0}, h_0\right)$ term in \eqref{Fl}, we prove the following claim.
	
	 \textbf{Claim.} On smaller tubular neighborhoods of $\gamma$ contained in $\mathcal{A}$ the Osgood-Stowe tensor satisfies $B\left(h_{\infty,0}, h_0\right)=Q_2\left(e^{-\pi^2/(2\ell)}\right) $.
	 
	\bpfc Using the isometric change of coordinates $(r, \theta)$ in which the tube $(\mathcal{A}, h_{\infty, 0})$ is expressed as $[-m, m]\times \mathbb{S}^1$, where $m$ is half of the length of $\mathcal{A}$ with respect to $h_{\infty, 0}$ (see Figure \ref{cylpic}), we can calculate the Osgood-Stowe tensor $B\left(h_{\infty,0}, h_0\right)=\B{h_{\infty, 0}} {u_1}$ with Definition \ref{tensor} as:
	
	\[ \B{h_{\infty, 0}} {u_1}= \begin{bmatrix} \f 12 (u_1)_{rr}-\f12 (u_1)_{\theta\theta}+\f 12 (u_1)_\theta^2-\f 12 (u_1)_r^2& (u_1)_{\theta r}-(u_1)_\theta (u_1)_r  \\ (u_1)_{\theta r}-(u_1)_\theta (u_1)_r & \f12 (u_1)_{\theta\theta}-\f12 (u_1)_{rr}+\f 12 (u_1)_r^2-\f12(u_1)_\theta^2  \end{bmatrix}~.\]
The claim now follows by the last block of estimates in Section \ref{est}, which bound the derivatives of $u_1$ in $\mathcal{A}$, together with the inequality 
\[e^{-m}\leq e^{-\pi^2/(2\ell)}~,\] 
proved in Remark \ref{mbound} as equation \eqref{mest}.\epfc

From equations \eqref{Fl}, \eqref{Bz} and the claim right above, we obtain $(F_{\ell})_{*}\re\left(S(F_{\ell})\right)=Q_2\left(e^{-\pi^2/(2\ell)}\right) $.
	Therefore, in a neighborhood of the core of $\mathcal{A}$ we have:
	\[f_{*}\left(\re(S(f))\right)=f_{*}\left(\re(S(f_{\ell}))\right)+Q_2\left(e^{-\pi^2/(2\ell)}\right)~.\]
	We are interested in the Schwarzian derivative on $\widetilde{\mathcal A}\subseteq \mathbb{H}^2$, and we first study its real part
	 \[\re\left(S(f)\right)=\re\left(S(f_{\ell})\right)+(f_{\ell})^*\re\left(S(F_{\ell}\right))~.\]
	Since $df=dF_\ell\circ df_{\ell}$, and, from Section \ref{sec:symmsch}, by a direct computation of the derivative of $f_{\ell}$, which is explicit, we can recover that the complex norm $\abs{df_{\ell}}$ is in $O(1/(\abs{z}\ell))$ as then also $\abs{df}$ is, we get \[(f_{\ell})^*\re\left(S(F_{\ell})\right)=f^*\left((F_{\ell})_{*}\re\left(S(F_{\ell})\right)\right)~=\f{1}{\abs{z}^2}Q_2\left(\f{e^{-\pi^2/(2\ell)}}{\ell^2}\right),\]
	where we used that $f_{\ell}^*=f^*\circ (F_{\ell})_{*}$, and therefore
	
	\[\re\left(S(f)\right)=\re\left(S(f_{\ell})\right)+\f{1}{\abs{z}^2}Q_2\left(\f{e^{-\pi^2/(2\ell)}}{\ell^2}\right)~.\]
	To get rid of the real part, given $z$ the coordinate of $\mathbb{H}^2$, we just need to notice that a holomorphic quadratic differential $Q(z)=q(z)dz^2$, where $q(z)=q_0(z)+iq_1(z)$ with $q_0$ and $q_1$ real valued functions, and $dz=dx+idy$, satisfies
	
	\[\re(S(f))=q_0(z)(dx^2-dy^2)-q_1(z)(dx\otimes dy + dy\otimes dx)\] 
	
	and 
	
	\[\im(S(f))=q_1(z)(dx^2-dy^2)+q_0(z)(dx\otimes dy + dy\otimes dx)~.\] 
	Then, the coefficients of its real part are $O(f(\ell))$ if and only if the ones of its imaginary part are as well. Finally, 
	
	\[S(f)=\re(S(f))+i\im(S(f))=S(f_{\ell})+\f{1}{\abs{z}^2}Q_2\left(\f{e^{-\pi^2/(2\ell)}}{\ell^2}\right)~.\]
	Noting that if $|z|^2|q(z)|$ is bounded and in $O\left(e^{-\pi^2/(2\ell)}/\ell^2\right)$ than $\norm{qdz^2}_{\infty}$ is bounded and in $O\left(e^{-\pi^2/(2\ell)}/\ell^2\right)$ on $\widetilde{\mathcal{A}}$, the equality right above concludes the proof thanks to Lemma \ref{symmlemma}.
	\end{proof}

We now define the two norms for the Osgood-Stowe tensor $B(g, e^{2u}g)$ that we are going to use and estimate in the following two subsections, which conclude with the proof of Theorem \ref{pairing}.
\bdefi\label{definB}
Given a smooth loop $\gamma$ we define:
\[ \norm{ \B g u}_\gamma\eqdef \abs{ \f12\int_\gamma \B g u(i\dot\gamma,\dot\gamma)} .\]
\edefi
The definition above is so that if $\mu$ is the infinitesimal earthquake along a unit speed simple closed geodesic $\gamma$ of a complex projective surface, we obtain, by Lemma \ref{intbypart} and Theorem \ref{tensorschw}, noting that $\left(\re(Q)\right)(\cdot, \cdot)=\re\left(Q(\cdot, \cdot)\right)$ for any quadratic differential $Q(z)$,
\[\vert\re \langle S(f), \mu\rangle\vert =\norm{ \B g u}_\gamma~.\]

Analogously, 

 \bdefi\label{definBgr}
Given a smooth loop $\gamma$ we define:
\[ \norm{ \B g u}^{gr}_\gamma\eqdef \abs{\f12 \int_\gamma \B g u(\dot\gamma,\dot\gamma)} .\]
\edefi
The definition above is so that if $\nu$ is the infinitesimal grafting along a unit speed simple closed geodesic $\gamma$, we obtain, by Lemma \ref{intbypart} and Theorem \ref{tensorschw} 
\[\vert\re \langle S(f), \nu\rangle\vert =\norm{ \B g u}^{gr}_\gamma~.\]

\subsection{Estimates on the Schwarzian along earthquakes} \label{ssc:estimatesearth}

In this section we are interested in studying $B(g, e^{2u}g)$ on the simple closed geodesic $\gamma$.

\brem\label{loopparam}
For Lemma \ref{intbypart}, in Definitions \ref{definB} and \ref{definBgr}, we need to take $\gamma$ to be the core geodesic, which in the metric $h_\infty$ and the coordinates $(\rho, \theta)\in [-L, L]\times \mathbb{S}^1$ is parameterized as:
\[\gamma(t)=\left(0,\f {2\pi}\ell t\right) ,\qquad t\in[0,\ell]~,\]
so that $\norm{\dot\gamma}_{h_\infty}=1$. We want to estimate $\norm{ \B \rho u}_\gamma$, with $u$ the logarithmic conformal factor between $h_{\infty}$ and $\abs{dz}^2$, through the following equality:
\[ f_*\B{\rho}{u\vert_{\A}}=\B{h_{\infty}}{u_0\vert_{\mathcal A}}+\B{h_{\infty, 0}}{u_1\vert_{\mathcal A}}+\B{h_0}{w\vert_{\mathcal A}}.\]
Since all the maps are the identity, the parameterisation of the curve $\gamma$ never changes. However, the curve is not necessarily geodesic in the other metrics. In particular, even if $\gamma$ has geodesic support for $h_{\infty,0}$ it is not in unit speed parametrisation as we have:
\[ \norm{\dot\gamma}_{h_{\infty,0}}= \f{2\pi}\ell~.\]\erem

We can finally state our estimate on the norm of $\B\rho u$. We first define the following auxiliary functions of $W$ and $m$, with $W$ a constant such that $0< W\leq 3.7$, which bounds $\abs{u_1}$ on the boundary $\partial\mathcal{A}$ (see Lemma \ref{lem1}), and $m$ half of the length of $(\mathcal{A}, h_{\infty, 0})$:
\[ G(\ell)\eqdef \min \set{e^{3.3},   1+4\pi W(3e^{3.3}+2)\f{e^m}{(e^m-1)^2}  +  96 W^2e^{3.3}\pi^2\left(  \f{e^m}{(e^m-1)^2}\right)^2}.\]
 and 
 \[\bar G(\ell)= \left( W  G(\ell) \left(\f {16}{15}\right)^2\right)^2\leq 17.7223G(\ell)^2~,\]
  see Remark \ref{Mell}. 

The main result we will now work towards is the following estimate:
\bthm\label{bgestimate} 
In the notations above, we have the following bound
\[\norm{\B {\rho} u }_{\gamma}\leq 144 \pi^4 \f{\bar G(\ell)}\ell  e^{-\pi^2/\ell}    +148\pi^2 \f{e^{-\pi^2/(2\ell)}}\ell~,\]
where $\ell$ is the length of the core $\gamma$ with respect to $h_{\infty}$.
\ethm

Theorem \ref{bgestimate} will follow from the equality between the Osgood-Stowe differentials
\[
f_*\B{\rho}{u\vert_{\A}}=\B{h_{\infty}}{u_0\vert_{\mathcal A}}+\B{h_{\infty, 0}}{u_1\vert_{\mathcal A}}+\B{h_0}{w\vert_{\mathcal A}}~,
\]
bounding, or computing, each term on the right-hand side, in order, in Propositions \ref{bounB1}, \ref{bounB2}, and \ref{bounB3}.

\bprop\label{bounB1}
For the Osgood-Stowe tensor $\B{h_{\infty}} {u_0}$: $\norm{\B{h_{\infty}} {u_0}}_{\gamma}=0$.\eprop

\bpf
We already noted that there exist coordinates $(\rho, \theta)$ in which $u_0(\rho, \theta)= \log\left(\f{2\pi}{\ell\cosh(\rho)}\right)$. Then, we can directly compute the Osgood-Stowe tensor $\B{h_{\infty}} {u_0}$. The only non-trivial Christoffel symbols for $h_\infty$ are 

\[\Gamma_{\rho\theta}^{\theta}=\Gamma_{\theta\rho}^{\theta}=\tanh(\rho),\qquad \Gamma_{\theta\theta}^{\rho}=-\left(\frac{\ell}{2\pi}\right)^2\sinh (\rho)\cosh (\rho)~,\]

so 

\[\hess(u_0)=-\f{1}{\cosh^2(\rho)}d\rho^2-\left(\frac{\ell}{2\pi}\right)^2\sinh^2(\rho)d\theta^2 ,\qquad d u_0\otimes d u_0=\tanh^2(\rho) d\rho^2~, \]
and their traces 

\[\Delta u_0= -1,\qquad \norm{\nabla u_0}_{h_\infty}^2= \tanh^2(\rho)~,\]
therefore, at the core $\rho=0$ 
\[\B{h_{\infty}} {u_0}=-\f 12 d\rho^2 +\f 12\left(\f{\ell}{2\pi}\right)^2 d\theta^2~.\]
The result now follows from Definition \ref{definB} and the fact that the unit speed parametrization with respect to $h_{\infty}$ of the core $\gamma$ has derivative $\dot\gamma(t) =  \f{2\pi}{\ell} \f{\partial}{\partial\theta}$, $t\in [0, \ell]$, and $i\dot\gamma(t) =  -\f{\partial}{\partial\rho}$.\epf

\bprop\label{bounB2}The norm of the Osgood-Stowe tensor $\B{h_{\infty, 0} }{u_1}$ is bounded as:
\[\norm{\B{h_{\infty, 0}}{u_1} }_{\gamma}\leq148\pi^2\frac{e^{-m}}{\ell},  \quad m=\f{2\pi\arctan(\sinh(L))}\ell,\quad  L=\argsinh\left(\f{1}{\sinh\left(\f \ell 2\right)}\right).\]
\eprop

\begin{proof} 

We recall that $2m$ is the length of $(\mathcal{A}, h_{\infty, 0})$ and that
 \[m=\int_0^{L} \f{2\pi}{\ell\cosh{\rho}} d\rho=\f{2\pi\arctan(\sinh(L))}\ell~,\] 
with $L=\argsinh\left(\f{1}{\sinh\left(\f \ell 2\right)}\right)$. Then, $e^{-m}$ is an increasing function of $\ell$, and, since $\ell\leq\epsilon_0=2\argsinh(1)$, we can bound $e^{-m}$ by:
 \begin{align*} e^{-m} &\leq \exp\left(-\f{2\pi}{\epsilon_0} \arctan\left(\sinh\left(\argsinh\left(\f{1}{\sinh\left(\f {\epsilon_0} 2\right)}\right)\right)\right)\right)\\
 	&\leq  \exp\left(-\f{2\pi}{\epsilon_0} \arctan(1)\right) \\
 	&\leq \f{1}{16}~.
 \end{align*}
 
The metrics $h_{\infty, 0}$ and $h_0=e^{2u_1}h_{\infty,0}$ are both Euclidean, and, by the conformal change of curvature:
\[K_0=e^{-2u_1}(K_{\infty,0}+\Delta u_1),\]
for $K_{\infty,0}$, $K_0$ the curvatures of $h_{\infty, 0}$ and $h_0$ respectively. Thus, $u_1(r,\theta)$ satisfies 
\[\Delta u_1 = 0~.\] 
Then, we can apply Fourier analysis, see Subsection \ref{ssc:fourier}, on the round tube $\mathcal{A} \cong  [-m,m]\times \mathbb S^1 $ to get the following estimates of $u_1$ on the core $\gamma(\theta)=\left(0,\f{2\pi}\ell \theta\right)$, $\theta\in[0,\ell]$ which has geodesic image in $h_{\infty,0}=dr^2+d
\theta^2$ but is not in unit speed parametrisation, see Remark \ref{loopparam}.
 
 Then, by the estimates in Subsection \ref{est}, we have: 
\begin{itemize}
\item $\abs{(u_1)_r}\leq   8W \sum_{k>0 }k e^{-k m}=8W\f{e^{m}}{(e^{m}-1)^2}$;
\item $\abs{(u_1)_{\theta}}\leq 4W \sum_{k>0} k e^{-km}=4 W \f{e^m}{(e^m-1)^2}$;
\item $\abs{(u_1)_{\theta r} }\leq 8W \sum_{k>0} k^2  e^{-km}=8W\f{e^{m} (e^{m}+1)}{(e^{m}-1)^3}$;
\end{itemize}
for $0<W\leq 3.7$ as in Lemma \ref{lem1}.  Since, in the coordinates $(r, \theta)\in [-m,m]\times\mathbb{S}^1 $ we can express the Osgood-Stowe tensor as:
\[ \B{h_{\infty, 0}} {u_1}= \begin{bmatrix} \f 12 (u_1)_{rr}-\f12 (u_1)_{\theta\theta}+\f 12 (u_1)_\theta^2-\f 12 (u_1)_r^2& (u_1)_{\theta r}-(u_1)_\theta (u_1)_r  \\ (u_1)_{\theta r}-(u_1)_\theta (u_1)_r & \f12 (u_1)_{\theta\theta}-\f12 (u_1)_{rr}+\f 12 (u_1)_r^2-\f12(u_1)_\theta^2  \end{bmatrix},\]
and we need to compute it on $(i\dot\gamma,\dot\gamma)$ with $\gamma(\theta)=\left(0,\f{2\pi}{\ell}\theta\right)$ we have that:
\begin{align*}  \B{h_{\infty, 0}} {u_1}(i\dot\gamma,\dot\gamma) &=\B{h_{\infty, 0}} {u}\left(\left(-\f{2\pi}\ell,0\right),\left(0,\f{2\pi}\ell\right)\right)\\
&=-\left(\f{2\pi}\ell\right) ^2\left((u_1)_{\theta r}-(u_1)_\theta (u_1)_r\right)~,\end{align*}
whose norm is bounded by:
\begin{align*}
\norm{\B{h_{\infty,0}}{u_1}}_{\gamma}
&\leq
\frac \ell2 \left(\frac{2\pi}{\ell}\right)^2
\left(
8W\frac{e^m(e^m+1)}{(e^m-1)^3}
+
32W^2\left(\frac{e^m}{(e^m-1)^2}\right)^2
\right)
\\
&=
\frac{2\pi^2}{\ell}
\left(
8We^{-m}\frac{1+e^{-m}}{(1-e^{-m})^3}
+
32W^2e^{-2m}\frac{1}{(1-e^{-m})^4}
\right)
\\
&=
\frac{16\pi^2}{\ell}We^{-m}
\left(
\frac{1+e^{-m}}{(1-e^{-m})^3}
+
4We^{-m}\frac{1}{(1-e^{-m})^4}
\right)
\\
&\leq
\frac{16\pi^2}{\ell}We^{-m}
\left(
1.29+1.2
\right),
\qquad e^{-m}\leq \frac{1}{16},
\\
&\leq
148\pi^2\frac{e^{-m}}{\ell},
\qquad W\leq 3.7.
\end{align*}

which completes the proof.\end{proof}

\brem\label{mbound} In Proposition \ref{bounB2} we get: 
\[\norm{\B{h_{\infty, 0}}{u_1} }_{\gamma}\leq 148\pi^2\frac{e^{-m}}{\ell}~, \]
where
\[ m=\f{2\pi\arctan(\sinh(L))}\ell~,\quad L=\argsinh\left(\f 1 {\sinh(\ell/2)}\right)~. \]
Since $\f{\ell}{2\pi} m$ is a decreasing function in $\ell\in(0,\epsilon_0]$ its minimum is at  $\epsilon_0$ and it is $\f\pi 4$. Then
\[ m\geq \f{2\pi}\ell\f{\pi}{4}= \f{\pi^2}{2\ell},\] and therefore 
\begin{equation}\label{mest}
e^{-m}\leq e^{-\pi^2/(2\ell)}~.
\end{equation}
Thus, we can write:
\[\norm{\B{h_{\infty, 0}}{u_1} }_{\gamma}\leq 148\pi^2 \f{e^{-\pi^2/(2\ell)}}\ell\in O\left(\f{{e^{-\pi^2/(2\ell)}}}{\ell}\right).\]
Moreover, as we have $e^{-m}<\f 1{16}$ we also get that: \[(e^m-1)^{-1}<\f {16}{15}e^{-m}~.\]

\erem

Before being able to prove the estimate for the norm on $\gamma$ of $B\left(h_0, \abs{dz}^2\right)$, we show that for small $\ell$ the curve $\gamma$ is almost a round circle in $\C^*$, that is, $\gamma$ is also almost geodesic with respect to the metric $h_0$. In the next lemma we have that $W$ is as in Lemma \ref{lem1} and $D^0$ is the Levi-Civita connection with respect to $h_0$.

 \blem\label{almostround} Given $\gamma:[0,\ell]\rar \left(\mathcal A,h_{\infty,0}\right)$ the unit speed $h_\infty$-geodesic of the core then: 
 \[ \max_{s\in[0,\ell]} \norm{ D^0_{\dot\gamma}\dot\gamma(s)}_{h_0} \leq 12We^{3.3}\left( \f{2\pi}\ell\right)^2 \f{e^m}{(e^m-1)^2},\quad W\leq 3.7~.\]
Moreover, by Remark \ref{mbound} we have:
  \[ \max_{s\in[0,\ell]} \norm{ D^0_{\dot\gamma}\dot\gamma(s)}_{h_0} \leq 12We^{3.3} \left(\f{32\pi}{15}\right)^2 \f{e^{-\pi^2/(2\ell)}}{\ell^2}.\]
  \elem 
 \bpf We consider the coordinates $(r, \theta)$ such that $\mathcal A\cong [-m,m]\times\mathbb S^1 $, $h_{\infty,0}=dr^2+d\theta^2$, and $\gamma(s)=\left(0,\f{2\pi}\ell s\right)$ with $\dot\gamma(s)=\left(0,\f{2\pi}\ell\right)$. The metric $h_0$ is equal to $e^{2u_1}h_{\infty,0}$ and we denote by $D^{\infty,0}$ the Levi-Civita connection for $h_{\infty,0}$. In what follows we use $u=u_1$ for the sake of notation.
 
 The change of the Levi-Civita connection under conformal changes of metric, see \cite[1.159 a]{Besse}, is given by the following formula:
 \begin{align*}
 D^0_{\dot\gamma}\dot\gamma(s)&=D^{\infty,0}_{\dot\gamma}\dot\gamma(s)+2du(\dot\gamma(s))\dot\gamma(s)-\norm{\dot\gamma(s)}^2_{h_{\infty,0}}D^{\infty,0}u\\
 &=2\left(\f{2\pi}\ell \right)^2 u_\theta\defrac{}{\theta}-\left(\f{2\pi}\ell \right)^2\left(u_r\defrac{} r+u_\theta\defrac{}{\theta}\right)\\
 &=\left(\f{2\pi}\ell \right)^2\left(u_\theta\defrac{}{\theta}-u_r\defrac{} r\right)~.
 \end{align*}
 Taking norms with respect to $h_0$:
 \[ \norm{ D^0_{\dot\gamma}\dot\gamma(s)}_{h_0}\leq \left(\f{2\pi}\ell \right)^2\left(\abs{u_\theta}\norm{\defrac{}{\theta}}_{h_0}+\abs{u_r}\norm{\defrac{}{r}}_{h_0}\right)~. \]
Since $\norm{v}_{h_0}=e^{u}\norm{v}_{h_\infty,0}$ and $u\leq 3.3$ at the core by Remark \ref{cstcorebound}, we obtain:
 
 \[ \norm{ D^0_{\dot\gamma}\dot\gamma(s)}_{h_0}\leq e^{3.3}\left(\f{2\pi}\ell \right)^2\left(\abs{u_\theta}+\abs{u_r}\right)~. \]
Finally, by the Fourier analysis bounds of Subsection \ref{est} we have:
\[ \abs{u_r}\leq   8 W\f{e^{m}}{(e^{m}-1)^2}\quad \abs{u_\theta}\leq 4 W\f{e^{m}}{(e^{m}-1)^2}~, \]
where again $0<W\leq 3.7$. This yields:
\[ \norm{ D^0_{\dot\gamma}\dot\gamma(s)}_{h_0}\leq 12We^{3.3}\left(\f{2\pi}\ell \right)^2\f{e^{m}}{(e^{m}-1)^2}~. \]
Taking the maximum of the left-hand side over $s$ completes the first part of the statement, while the second part follows by the computations in Remark \ref{mbound}.  \epf

We now need to compute $\norm{\B{h_0} {\log(r)} }_{\gamma}$ which is the last term of our bound for $\norm{B(\rho,\abs{dz}^2)}_{\gamma}$. 

\brem \label{B3remark}
We remark that by Theorem \ref{tm:standard} (or Section $4$ of \cite{Epstein-wvolume-OS}), the Osgood-Stowe tensor $\B{h_0} {\log(r)}$ is equal to the real part of the Schwarzian derivative of the uniformization map from the infinite flat tube, with core $2\pi$, to $\C^*$, that is, the conformal covering $\exp\colon \mathbb{C}\rightarrow \mathbb{C}^*$. 
This is because a lift of $(\mathcal{A}, h_0)$ is isometrically contained in the vertical strip equipped with the restriction of the metric $|dz|^2$ in $\mathbb{C}$ (so the relative Osgood-Stowe tensor is zero), and the restriction of $\exp{(z)}$ going from the strip to $\mathbb{C}^{*}\subseteq \mathbb{C}$ has real part of the Schwarzian equal to the Osgood-Stowe tensor with respect to $|dz|^2$ on both sides, which is the standard one for $\mathbb{C}$.
\erem 
 \bprop \label{bounB3} We have the following bound:
 \[\norm{\B{h_0} {\log(r)} }_{\gamma}\leq 144 \pi^4 W^2 \f{G(\ell)^2}\ell \f{e^{2m}}{(e^m-1)^4}~,\]
for $G(\ell)$ as in Remark \ref{Mell} and $0<W\leq3.7$ as in Lemma \ref{lem1}, already defined above. Moreover, we have:

  \[ \norm{\B{h_0} {\log(r)} }_{\gamma}\leq 144\pi^4 W^2 \f{G(\ell)^2}\ell\left(\f {16}{15}\right)^4  e^{-\pi^2/\ell}~,\]
and $G(\ell)$ goes to $1$ as $\ell $ goes to zero.
 \eprop 
 
 \bpf Let $\gamma(t):[0,\ell]\rar (\mathcal A,h_0)$ be the parameterization of the loop $\gamma$ in the flat tube $\mathcal A$ lifted isometrically to the vertical infinite strip $\mathcal C\cong\set{ z\in\C\vert 0\leq\im z \leq 2\pi}\subseteq \mathbb{C}$, with coordinates $z=x+iy$. Note that in this setting $\gamma(t)=(x(t),y(t))$ is not geodesic for $h_0$. However, one sees that $\dot\gamma$ has mean $\left(0,\f{2\pi}\ell\right)$ i.e.:
 \[ \bar x'=\f 1 \ell\int_0^\ell \dot x(t)dt=0 \qquad \bar y'= \f 1 \ell \int_0^\ell \dot y(t)dt=\f{2\pi}\ell~.\]
 Moreover, by the Mean Value Theorem we have times $\xi_1$ and $\xi_2$ such that $\dot x(\xi_1)=\bar x $ and $\dot y(\xi_2)=\bar y$. Since the metric is flat Euclidean, the covariant derivative is:
 \[  D^0_{\dot\gamma}\dot\gamma(s)=(\ddot x(s),\ddot y(s))~,\]
 and by Lemma \ref{almostround} and Remark \ref{Mell}:
 \[\max_{s\in[0,\ell]} \norm{ D^0_{\dot\gamma}\dot\gamma(s)}_{h_0} \leq 12W  G(\ell)\left( \f{2\pi}\ell\right)^2 \f{e^m}{(e^m-1)^2}~,\]
 we obtain:
 \[ \forall s \in [0,\ell]:\abs{\ddot x(s)},\abs{\ddot y(s)} \leq 12W G(\ell) \left( \f{2\pi}\ell\right)^2 \f{e^m}{(e^m-1)^2}~.\]

 Note:
\[\norm{\dot \gamma(s)-\bar \gamma'}_{h_0} \leq \frac{1 }{\ell}\int_{0}^\ell\norm{\dot \gamma(s)- \dot\gamma(t)}_{h_0}dt\leq \frac \ell 2\max_t \norm{\ddot \gamma(t)}_{h_0}\]
   Thus,
   \begin{equation}\label{normbound} 
   \norm{\dot \gamma(s)-\bar \gamma'}_{h_0} \leq 6  \ell W G(\ell) \left( \f{2\pi}\ell\right)^2 \f{e^m}{(e^m-1)^2}=\f{  24 W G(\ell) \pi^2}\ell\f{e^m}{(e^m-1)^2}~.
 \end{equation}
 We will now consider the restriction of the uniformization map from the infinite flat strip $\mathcal C$ of core $2\pi$ to $\C^*$. The map is 
 \[ f(z)=\lambda \exp(z)~,\]
 where $z$ is the coordinate on the infinite strip $\mathcal C\cong\set{ z\in\C\vert 0\leq\im z \leq 2\pi}$ and $\lambda$ is a complex number. As $\f{\ddot f}{\dot f}=1$ we have:
 \[ S(f)(z)=-\f 12 dz^2~.\]
 Since $\B{h_0} {\log(r)}$ on the core is equivalent to the real part of the Schwarzian right above, see Remark \ref{B3remark}, we need to compute
 \[ \re\left(\int_0^\ell S(f)(i\dot\gamma(s),\dot\gamma(s))ds\right)=-\im\left(\int_0^\ell S(f)(\dot\gamma(s),\dot\gamma(s))ds\right).\]
 Substituting $S(f)$ and considering $\dot\gamma(s)$ as in $\C$ we have:
 \begin{align*}
 -\im\left(\int_0^\ell S(f)(\dot\gamma(s),\dot\gamma(s))ds\right)&=\f 12\int_0^\ell\im(\dot\gamma(s)^2)ds\\
 &=\f 12\im \int_0^\ell\dot\gamma(s)^2ds~.
 \end{align*}
 Since the mean of $\dot\gamma$ is $\bar\gamma'=\f{2\pi}\ell i$ we have:
  \begin{align*}
\f 12\im \int_0^\ell\dot\gamma(s)^2ds&=\f12\im \int_0^\ell (\bar\gamma')^2ds+\f12\im \int_0^\ell(\dot\gamma(s)-\bar\gamma')^2ds+\im \int_0^\ell\bar\gamma'(\dot\gamma(s)-\bar\gamma')ds\\
&=\f12\im \int_0^\ell \left(\f{2\pi}\ell i\right)^2ds+\f12\im \int_0^\ell(\dot\gamma(s)-\bar\gamma')^2ds\\
&=\f12\im \int_0^\ell(\dot\gamma(s)-\bar\gamma')^2ds~.
 \end{align*}
 Then,
 \begin{align*}
 \norm{\B{h_0} {\log(r)} }_{\gamma}&=\abs{\f12\int_0^\ell \B{h_0} {\log(r)}(i\dot\gamma,\dot\gamma)ds}\\
 &=\abs{\re\left(\f12\int_0^\ell S(f)(i\dot\gamma(s),\dot\gamma(s))ds\right)}\\
 &\leq \f14 \int_0^\ell \norm{\dot\gamma(s)-\bar\gamma'}^2_{h_0}ds\\
 &\overset{\eqref{normbound} }{\leq}  \f\ell 4 \left(\f{ 24 W G(\ell) \pi^2}\ell\f{e^m}{(e^m-1)^2}\right)^2\\
 &\leq 144 \pi^4 W^2 \f{(G(\ell))^2}\ell \f{e^{2m}}{(e^m-1)^4}~.
  \end{align*} The second part follows by the computations in Remark \ref{mbound} and \ref{Mell}.   \epf
 
  We can now prove the main estimate, and recall that $\bar G(\ell)= \left( W G(\ell) \left( \f {16}{15}\right)^2\right)^2$.
 \begin{customthm}{\ref{bgestimate}}
We have the following bound
\[\norm{\B {\rho} u }_{\gamma}\leq 144 \pi^4 \f{\bar G(\ell)}\ell   e^{-\pi^2/\ell}    +148\pi^2 \f{e^{-\pi^2/(2\ell)}}\ell~,\]
where $\ell$ is the length of the core $\gamma$ with respect to $h_{\infty}$.\end{customthm}

 \bpf We have that:
  \small
 \[ \norm{\B \rho u}_{\gamma}= \abs{\frac 1 2 \int_{\gamma} \left( \B {h_{\infty}} {u_0}+\B{h_{\infty, 0}} {u_1}+\B{h_0}{\log(r)}\right) (i\dot\gamma,\dot\gamma)}\]\normalsize
 which gives:
 \[ \norm{\B\rho u}_{\gamma}\leq \norm{ \B {h_{\infty}} {u_0}}_{\gamma}+\norm{ \B {h_{\infty,0}} {u_1}}_{\gamma}+\norm{ \B{h_0} {\log(r)}}_{\gamma},\]
 which by Proposition \ref{bounB1}, Proposition \ref{bounB2},  Remark \ref{mbound}, and Proposition \ref{bounB3} yields:
 \[\norm{\B\rho u}_{\gamma}\leq 144 \pi^4 W^2 \f{G(\ell)^2}\ell\left(\f {16}{15}\right)^4 e^{-\pi^2/\ell}+148\pi^2 \f{e^{-\pi^2/(2\ell)}}\ell~,\]
 completing the proof.  \epf

  \subsection{Estimates on the Schwarzian along 
 	grafting}\label{sec:graft}

 In this section, we repeat the computations done in Section \ref{ssc:estimatesearth}, but this time for the norm of Definition \ref{definBgr}.
 As before, we are going to estimate the norm of the Osgood-Stowe tensor $B$ at each step. More specifically, we want to prove the following theorem. 
 
 \bthm\label{bgestimategr} 
	 We have the following bound
 	\[\abs{
\frac12\int_\gamma
B(\rho,e^{2u}\rho)(\dot\gamma,\dot\gamma)
-\frac{\pi^2}{\ell}}\leq \frac{\ell}{4}+ 144 \pi^4 \f{\bar G(\ell)}\ell  e^{-\pi^2/\ell}    +127\pi^2 \f{e^{-\pi^2/(2\ell)}}\ell~,\]
 	where $\ell$ is the length of the core $\gamma$ with respect to $h_{\infty}$. In particular $ \norm{\B \rho u}_{\gamma}^{gr}\sim\frac{\pi^2}\ell$ as $\ell\rar 0$.
 \ethm

 Theorem \ref{bgestimategr} will directly follow, analogously to Theorem \ref{bgestimate}, from splitting $\norm{\B {\rho} u }^{gr}_{\gamma}$ into the same three terms and bounding each of them. This is done in Propositions \ref{bounB1gr}, \ref{bounB2gr}, and \ref{bounB3gr} respectively.
 
 \bprop\label{bounB1gr}
 For the Osgood-Stowe tensor $\B{h_{\infty}} {u_0}$: \[\norm{\B{h_{\infty}} {u_0}}^{gr}_{\gamma}=\f14 \ell~,\]
 where $\ell$ is the length of the core $\gamma$ with respect to $h_{\infty}$.\eprop
 
 \bpf
 We have already seen in the proof of Proposition \ref{bounB1} that, in the coordinates $(\rho, \theta)$ in which $u_0(\rho, \theta)= \log\left(\f{2\pi}{\ell\cosh(\rho)}\right)$, on the core $\gamma$ the Osgood-Stowe tensor is given by
 \[\B{h_{\infty}} {u_0}=-\f 12 d\rho^2 +\f 12\left(\f{\ell}{2\pi}\right)^2 d\theta^2~.\]
 As $\dot\gamma=\f{2\pi}{\ell}\f{\partial} {\partial\theta}$ by a direct computation 
 \[\norm{\B{h_{\infty}} {u_0}}^{gr}_{\gamma}=\abs{\f12\int_0^\ell \B{h_{\infty}} {u_0}(\dot\gamma,\dot\gamma)ds}=\abs{\int_0^\ell \f14 ds }=\f\ell4~.\]
 This concludes the proof.  \epf
 
 \bprop\label{bounB2gr}The norm of the Osgood-Stowe tensor $\B{h_{\infty, 0} }{u_1}$ is bounded as:
 \[\norm{\B{h_{\infty, 0}}{u_1} }^{gr}_{\gamma}\leq 127\pi^2\f{e^{-m}}{\ell},  \quad m=\f{2\pi\arctan(\sinh(L))}\ell,\quad  L=\argsinh\left(\f{1}{\sinh\left(\f \ell 2\right)}\right) .\]
 \eprop
 
 \bpf
 As seen in the proof of Lemma \ref{bounB2}, from Subsection \ref{est} we have the following estimates for the derivatives of $u_1$ at the core:
 
 \begin{itemize}
 	\item $\abs{(u_1)_{r}}\leq   8W \sum_{k>0} k e^{-k m}=8W\f{e^{m}}{(e^{m}-1)^2}$;
 	\item $\abs{(u_1)_{\theta}}\leq 4W \sum_{k>0} k e^{-km}=4 W \f{e^m}{(e^m-1)^2}$;
 	\item $\abs{(u_1)_{rr}}\leq 4 W \sum_{k>0} k^2  e^{-km}=4 W\f{e^{m} (e^{m}+1)}{(e^{m}-1)^3}$;
 	\item $\abs{(u_1)_{\theta r}}\leq 8 W \sum_{k>0} k^2  e^{-km}=8 W\f{e^{m} (e^{m}+1)}{(e^{m}-1)^3}$;
 	\item $\abs{(u_1)_{\theta\theta}}\leq 4 W \sum_{k>0} k^2  e^{-km}=4 W\f{e^{m} (e^{m}+1)}{(e^{m}-1)^3}$;
 \end{itemize}
 for $0<W\leq 3.7$ as in Lemma \ref{lem1}, and the Osgood-Stowe tensor in the coordinates $(r, \theta)\in [-m,m]\times\mathbb{S}^1 $ is
 \[ \B{h_{\infty, 0}} {u}= \begin{bmatrix} \f 12 (u_1)_{rr}-\f12 (u_1)_{\theta\theta}+\f 12 (u_1)_\theta^2-\f 12 (u_1)_r^2& (u_1)_{\theta r}-(u_1)_\theta (u_1)_r  \\ (u_1)_{\theta r}-(u_1)_\theta (u_1)_r & \f12 (u_1)_{\theta\theta}-\f12 (u_1)_{rr}+\f 12 (u_1)_r^2-\f12(u_1)_\theta^2  \end{bmatrix}~.\]
 We now need to compute $\B{h_{\infty, 0}} {u_1}$ on $(\dot\gamma,\dot\gamma)$ with $\gamma(s)=(0,\f{2\pi}{\ell}s)$:
  \begin{align*}  \B{h_{\infty, 0}} {u_1}(\dot\gamma,\dot\gamma) &=\B{h_{\infty, 0}} {u}\left(\left(0,\f{2\pi}\ell\right),\left(0,\f{2\pi}\ell\right)\right)\\
 	&=\left(\f{2\pi}\ell\right) ^2\left(\f12 (u_1)_{\theta\theta}-\f12 (u_1)_{rr}+\f 12 (u_1)_r^2-\f12(u_1)_\theta^2\right)~,\end{align*}
 and, remembering that $e^{-m}\leq 1/16$, see Remark \ref{mbound}, we can bound the norm as follows:
 \begin{align*}
\norm{ \B{h_{\infty, 0}} {u}}_{\gamma}^{gr}
&\leq \f{\ell}{4}\left(\f{2\pi}\ell\right)^2
\left(
8 W\f{e^{m} (e^{m}+1)}{(e^{m}-1)^3}
+
80W^2\left( \f{e^m}{(e^m-1)^2}\right)^2
\right)
\\
&=
\f{\pi^2}{\ell}
\left(
8We^{-m}\f{1+e^{-m}}{(1-e^{-m})^3}
+
80W^2e^{-2m}\f{1}{(1-e^{-m})^4}
\right)
\\
&\leq
\f{\pi^2}{\ell}
\left(
8We^{-m}\f{17}{16}\left(\f{16}{15}\right)^3
+
80W^2e^{-2m}\left(\f{16}{15}\right)^4
\right)
\\
&=
\f{\pi^2}{\ell}We^{-m}
\left(
8\f{17}{16}\left(\f{16}{15}\right)^3
+
80We^{-m}\left(\f{16}{15}\right)^4
\right)
\\
&\leq
127\pi^2\f{e^{-m}}{\ell},
\end{align*}
concluding the proof.
 \epf
 
  For the already defined $G(\ell)$ as in Remark \ref{Mell}, $\bar{G}(\ell)=\left( W G(\ell) \left( \f {16}{15}\right)^2\right)^2$ and $0<W\leq 3.7$ as in Lemma \ref{lem1}, we now bound the norm $\norm{\cdot}^{gr}_{\gamma}$ of the Osgood-Stowe tensor $\B{h_0} {\log(r)}$.
 \bprop \label{bounB3gr} 
 The following bound holds
 \[ \abs{\f12\int_0^\ell \B{h_0} {\log(r)}(\dot\gamma,\dot\gamma)ds-\f{\pi^2}{\ell} }\leq144\pi^4 \f{\bar{G} (\ell)}\ell  e^{-\pi^2/\ell}~.\] Thus, $\norm{\B{h_0} {\log(r)} }^{gr}_{\gamma}\sim \frac{\pi^2}\ell$ as $\ell\rar0$.\eprop
 
 \bpf
 The proof of Lemma \ref{bounB3} essentially goes through and we just have to use the norm $\norm{\cdot}^{gr}_{\gamma}$. Denoting by $q=-S(f)=\f{1}{2}dz^2$ we have: 
 \begin{align*}
 	\f 12\int_{\R/\ell\Z}\re(q(\dot\gamma(t),\dot\gamma(t))) dt&=\f 14\re \int_0^\ell\dot\gamma(s)^2ds\\
 	&=\f14\re \int_0^\ell \left(\f{2\pi}\ell i\right)^2ds+\f14\re \int_0^\ell(\dot\gamma(s)-\bar\gamma')^2ds\\
 	&=-\f{\pi^2}{\ell}+\f14\re \int_0^\ell(\dot\gamma(s)-\bar\gamma')^2ds
 	~.
 \end{align*}
 Then,
 \begin{align*}
 	\abs{\f12\int_0^\ell \B{h_0} {\log(r)}(\dot\gamma,\dot\gamma)ds- \f{\pi^2}{\ell}}&=\abs{\re\left(\f12\int_0^\ell S(f)(\dot\gamma(s),\dot\gamma(s))ds\right)-\f{\pi^2}{\ell}}\\
 	&\leq \f14 \int_0^\ell \norm{\dot\gamma(s)-\bar\gamma'}^2_{h_0}ds\\
 	& \overset{\eqref{normbound} }{\leq}   \f\ell 4  \left(\f{  24 W G(\ell) \pi^2}\ell\f{e^m}{(e^m-1)^2}\right)^2\\
 	& =144 \pi^4 W^2 \f{G(\ell)^2}\ell \f{e^{2m}}{(e^m-1)^4}~.
 \end{align*}
 By applying the estimate in Remark \ref{mbound} \[e^{-m}\leq e^{-\pi^2/(2\ell)},\  \quad (e^m-1)^{-1}<\f {16}{15}e^{-m}~,\]
 and looking at the definition of $\bar{G}(\ell)$ right above the statement, we obtain the required result. \epf

  We can now prove the main estimate of this section.
 
 \begin{customthm}{\ref{bgestimategr}}
 	 We have the following bound
 	\[\abs{
\frac12\int_\gamma
B(\rho,e^{2u}\rho)(\dot\gamma,\dot\gamma)
-\frac{\pi^2}{\ell}}\leq \frac{\ell}{4}+ 144 \pi^4 \f{\bar G(\ell)}\ell  e^{-\pi^2/\ell}    +127\pi^2 \f{e^{-\pi^2/(2\ell)}}\ell~,\]
 	where $\ell$ is the length of the core $\gamma$ with respect to $h_{\infty}$. In particular $ \norm{\B \rho u}_{\gamma}^{gr}\sim\frac{\pi^2}\ell$ as $\ell\rar 0$.
 \end{customthm}
 
 \begin{proof}
 	 We have that:\small
 	\[ \norm{\B \rho u}_{\gamma}^{gr}= \abs{\frac 1 2\int_{\gamma} \left( \B {h_{\infty}} {u_0}+\B{h_{\infty, 0}} {u_1}+\B{h_0}{\log(r)}\right) (\dot\gamma,\dot\gamma)}\]\normalsize
Moreover, keeping track of the sign before applying the absolute value in
Proposition~\ref{bounB3gr}, and using Propositions~\ref{bounB1gr} and~\ref{bounB2gr}, we obtain
\[
\left|
\frac12\int_\gamma
B(\rho,e^{2u}\rho)(\dot\gamma,\dot\gamma)
-\frac{\pi^2}{\ell} 
\right|
\leq \frac{\ell}{4}+ 144\pi^4\frac{\bar G(\ell)}{\ell}e^{-\pi^2/\ell}
+
127\pi^2\frac{e^{-\pi^2/(2\ell)}}{\ell}.
\]
In particular, $\norm{\B{\rho}u}^{gr}_{\gamma}\sim\frac{\pi^2}{\ell}$ as $\ell\to0$. \end{proof}

 We now restate and prove Theorem \ref{pairing} in a more expanded form, namely we give explicit expressions for the functions $F_e(\ell)$ and $F_{gr}(\ell)$. Note that these are the functions that bound, from above, a complex pairing which depends on the whole complex structure. However, the bounds are explicit functions only depending on $\ell$.

\begin{customthm}{\ref{pairing}}\label{mainthm}
	Let $Z$ be a complex projective surface containing a long tube of core $\gamma$, let $S(f)$ be its Schwarzian, and let $X=\pi(Z)$ its underlying Riemann surface. Let also $\mu$ and $\nu$ be, respectively, the infinitesimal earthquake and grafting on the simple closed curve $\gamma\subseteq X$ of hyperbolic length $\ell\leq \varepsilon_0$. Then 
	
	\[\abs{\re\langle S(f), \mu\rangle}\leq 2553 \pi^4 \f{ (G(\ell))^2}\ell  e^{-\pi^2/\ell}    +148\pi^2 \f{e^{-\pi^2/(2\ell)}}\ell~, \] 
	and 
	\[\left|\re\langle S(f), \nu\rangle+\f{\pi^2}{\ell}\right|\leq  \f\ell4+ 2553 \pi^4 \f{ (G(\ell))^2}\ell  e^{-\pi^2/\ell}    +127\pi^2 \f{e^{-\pi^2/(2\ell)}}\ell ~, \]
 where $G (\ell)\rar 1$ as $\ell\rar 0$ and is always bounded by $e^{3.3}$.
\end{customthm}

\begin{proof}
	Thanks to Lemma \ref{intbypart}, the norms of the  pairings between $S(f)$ and $\mu$ or $\nu$ coincide, respectively, with the norms of the Osgood-Stowe differential of Definitions \ref{definB} and \ref{definBgr}. Therefore, the result follows from Theorem~\ref{bgestimate} and from the signed
estimate established in the proof of Theorem~ \ref{bgestimategr}.\end{proof}

\subsection{Fourier Analysis}\label{ssc:fourier}
 
In this section we bound, on a neighborhood of the core $\gamma$ of the thin tube $\mathbb{T}(\ell)$, the function $u_1$ and its derivatives. Recall that $u_1$ is the logarithmic conformal factor between the flat metrics $h_{\infty,0}$ and $h_0=e^{2u_1}h_{\infty,0}$, and that $\ell\leq \epsilon_0$ denotes the length of the core $\gamma$ with respect to the hyperbolic metric $h_{\infty}$.

 The relationship between $\mathbb{T}(\ell)$ and $\mathcal A$ is explained at the beginning of Section \ref{sc:estimates}, but the relevant fact here is that $\mathbb{T}(\ell)$ and $(\mathcal A,h_\infty)$ are isometric. We will use both viewpoints as needed and specifically when concerned with hyperbolic geometry we will use $\mathbb{T}(\ell)$. For our study we need a new metric on $\mathcal{A}$, called the \textit{Thurston metric} (see Section \ref{ssc:thmetric}). Recall that $\mathcal{A}$ is a long complex projective tube of a complex projective surface $Z$ with associated developing map $f$. Therefore, by Definition \ref{defilongprojtube}, the restriction of $f$ to the wide neighborhood $\widetilde{\mathcal{A}}$ of $\gamma$ is a covering with deck group $\langle\gamma\rangle$. In particular, the push-forward of the Thurston metric on $\mathcal{A}$ is well-defined.

 As per Section \ref{ssc:thmetric}, we identify  $\mathbb{CP}^1$ with the conformal boundary of $\mathbb{H}^3$, and we denote by $CH(\mathbb {CP}^1\setminus f(\cup_{z\in\widetilde{\mathcal{A}}}D_z))$ the smallest closed convex subset of $\mathbb{H}^3$ whose closure in $\mathbb{CP}^1$ is $\mathbb{CP}^1\setminus f(\cup_{z\in\widetilde{\mathcal{A}}}D_z)$, and by $\partial CH(\mathbb {CP}^1\setminus f(\cup_{z\in\widetilde{\mathcal{A}}}D_z))$ its boundary. We also remark that $\partial CH(\mathbb {CP}^1\setminus f(\cup_{z\in\widetilde{\mathcal{A}}}D_z))$ is a locally convex pleated surface, and therefore it carries a well-defined induced hyperbolic metric. In what follows below, in order to lighten the notation, we set $D_{\mathcal{A}}\eqdef f(\cup_{z\in\widetilde{\mathcal{A}}}D_z)$.

We start by bounding $u_1$ on $\partial \mathcal{A}$, and, for technical reasons, on a shorter curve of length $\varepsilon_0$ inside the tube. The second estimate is used to have a sharper bound on $u_1$ at the center of the tube, i.e. on $\gamma$. In what follows we will use the Thurston metric $h_{Th}$ and we define:
\begin{itemize}
\item $h_{Th}=e^{2u_2} h_0$ on $\mathcal{A}$;
\item $h_\infty=e^{2u_3} h_{Th}$ on $\pi(Z)$, and, in particular, on the restriction to the annulus $\mathbb{T}(\ell)$.
\end{itemize} 

Let $\alpha_{d}$ be the closed curve in the thin tube $\mathbb{T}(\ell)$ given by one connected component of the set of the points at distance $d$ from $\partial \mathbb{T}(\ell)$. View $\alpha_d$ in $\mathcal{A}$, let $\widetilde{\alpha}_d\subset \partial CH(\mathbb {CP}^1\setminus D_{\mathcal{A}})$ be its pull-back via the normal projection from the boundary $\partial CH(\mathbb {CP}^1\setminus D_{\mathcal{A}})$ to $D_{\mathcal{A}}$, and $\ell_{\partial CH(\mathbb {CP}^1\setminus D_{\mathcal{A}})}(\widetilde{\alpha}_d)$ be the length of $\widetilde{\alpha}_d$ with respect to the induced hyperbolic metric on $\partial CH(\mathbb {CP}^1\setminus D_{\mathcal{A}})$. We recall that $D_{\mathcal{A}}\eqdef f(\cup_{z\in\widetilde{\mathcal{A}}}D_z)$ is the union of the images of the maximal disks of $\widetilde{\mathcal A}$ and that by Definition \ref{defilongprojtube} we have that $D_{\mathcal A}\subset\mathbb C^*$.

\blem\label{horosphere}
For $\alpha_d$ and $u_2$ as above, the restriction of $u_2$ to $\alpha_d$ satisfies 
\[ 0\leq u_2|_{\alpha_d}\leq \f{\ell_{\partial CH(\mathbb{CP}^1\setminus D_{\mathcal{A}})}(\tilde{\alpha}_d)}{2}\leq b(\ell_{\infty}(\alpha_d))~,\]
with 
\[b(x)=2\pi e^{0.502\pi}e^{-\f{\pi^2}{\sqrt{e}x}} \qquad \ell_{\infty}(\alpha_d)=\ell \cosh(L-d)\qquad L=\operatorname{arcsinh}\left(\f{1}{\sinh(\ell/2)}\right)~.\]

\elem

\begin{proof}
  By definition, the Thurston metric at some point $z$ is given by $\dfrac{1}{r^2_{Th}}\abs{dz}^2$ where $r_{Th}$ is the radius of the horosphere in $\mathbb{H}^3$ at $z$ tangent to $\partial CH(\mathbb{CP}^1\setminus D_{\mathcal{A}})$ (see Section \ref{ssc:thmetric}). Note that, since the domain $D_{\mathcal{A}}$ is contained in $\mathbb{C}^*$, the geodesic $\mathcal L\subset\mathbb H^3$ connecting $0$ and $\infty$ belongs to $CH(\mathbb{CP}^1\setminus D_{\mathcal{A}})$, and, as $\mathcal{A}$ is essential in $\mathbb{C}^*$, it goes through a compression disk\footnote{This is any disk with boundary isotopic to $\Tilde{\alpha}_d$.} for $\Tilde{\alpha}_d$. Then, $r_{Th}$ is less than the radius $r_{0}$ of the horosphere at $z$ tangent to $\mathcal{L}$, therefore $u_2=\log{\dfrac{r_{0}}{r_{Th}}}$ is positive.
    
    We now prove the other inequality. Let us consider the pull-back of $\alpha_d$ through the normal projection $\Tilde{\alpha}_d$. We notice that any point $p\in \Tilde{\alpha}_d\subseteq \partial CH(\mathbb{CP}^1\setminus D_{\mathcal{A}})$ stays at distance from $\mathcal L$ less than $\ell_{\partial CH(\mathbb{CP}^1\setminus D_{\mathcal{A}})}(\Tilde{\alpha}_d)/2$. For any point $z\in \alpha_d\subseteq \mathbb{C}$ we call $p_z\in \partial CH(\mathbb{CP}^1\setminus D_{\mathcal{A}})$ the point on the horosphere $H_{z}$ centered at $z$ and tangent to $\partial CH(\mathbb{CP}^1\setminus D_{\mathcal{A}})$. Then we have
 \[\ell_{\partial CH(\mathbb{CP}^1\setminus D_{\mathcal{A}})}(\Tilde{\alpha}_d)/2\geq d_{\mathbb{H}^3}(p_z,\mathcal L)\geq d_{\mathbb{H}^3}(H_z,\mathcal L)=\log(r_{0})-\log(r
    _{Th})=u_2~.\] 
 Since $L$ is half the length of the thin tube $\mathbb{T}(\ell)$ and $\alpha_d$ stays at distance $L-d$ from the core of the tube, the length of $\alpha_d$ with respect to $h_{\infty}$ is 
 \[ \ell_\infty(\alpha_d)=\ell \cosh(L-d)~.\]
 The bound of $\ell_{\partial CH(\mathbb{CP}^1\setminus D_{\mathcal{A}})}(\tilde{\alpha}_d)$ by a function of $\ell_{\infty}(\alpha_d)$ follows from \cite[Theorem 5.1]{C2001}\footnote{Note that in \cite[Theorem 5.1]{C2001} the author requires the length to be less than one, however that is only needed to get a linear bound, see page 11 of \cite[Theorem 5.1]{C2001} the second line of equations.}. \end{proof} 

\brem \label{alphazero} In the notation of the previous lemma. if $\alpha_0$ coincides with a component of $\partial \mathbb{T}(\ell)$, its length with respect to $h_{\infty}$ is 

\[\ell_{\infty}(\alpha_0)=\ell \cosh{\left(\operatorname{arcsinh}{\left(\f{1}{\sinh(\ell/2)}\right)}\right)}~.\] 
For $\ell\in (0,\varepsilon_0]$ the length $\ell_{\infty}(\alpha_0)$ is increasing, it is maximized at $\epsilon_0$ by
 \[ \varepsilon_0 \cosh{\left(\operatorname{arcsinh}{\left(\f{1}{\sinh(\varepsilon_0/2)}\right)}\right)}< 2.5~, \]
  and it is greater than its limit at $0$, which is $2$, therefore
  \[2\leq \ell_{\infty}(\alpha_0)\leq 2.5~.\]
\erem

In what follows we let $\alpha=\alpha_{d_{\epsilon_0}}$ be the simple closed curve on $\mathbb{T}(\ell)$ given by one connected component of the set of points staying at distance $d_{\varepsilon_0}$ from $\partial\mathbb{T}(\ell)$, with $d_{\varepsilon_0}\geq 0$ such that $\ell_{\infty}(\alpha)=\ell\cosh(L-d_{\varepsilon_0})=\varepsilon_0$. 
      
\blem \label{lem1}
There exists a constant $W=W(\ell)>0$ such that if $h_0=e^{2u_1}h_{\infty,0}$, then $\abs{ u_{1\vert_{ \partial \mathcal{A}}}}\leq W\leq 3.7$. Moreover, with the notation above there exists a constant $W_0=W_0(\ell)>0$ such that $\abs{u_{1\vert _{\alpha}}}\leq W_0\leq 2.3$.\elem

\bpf
Recall that by $\alpha_d$ we mean the simple closed curve given by one connected component of the set of points at distance $d$ from $\partial \mathbb{T}(\ell)$. Note that, in particular, $\alpha_0$ is one connected component of $\partial \mathcal{A}$, and $\alpha_{d_{\varepsilon_0}}=\alpha$. We also recall that $L=:\operatorname{arcsinh}\left(1/\sinh(\ell/2)\right)$ is half of the length of the thin tube $\mathbb{T}(\ell)$. From the proof one can get explicit formulas $W(\ell)$, however for simplicity we just give upper bounds at every step.

We have $u_1=-( u_2+u_3+u_0)$ with $u_{0\vert _{\alpha_d}}=\log\left(\f{2\pi}{\ell_{\infty}(\alpha_d)}\right)$, and $0\leq u_{2\vert_{ \alpha_d}} \leq b(\ell_{\infty}(\alpha_d))$, as in Lemma \ref{horosphere}, and
\[\ell_\infty(\alpha_d)=\ell \cosh{(L-d)}~.\] 
In particular
\[\ell_\infty(\partial\mathcal A)=\ell \cosh{\left(\operatorname{arcsinh}{\left(\f{1}{\sinh(\ell/2)}\right)}\right)}~,\] 
and also
\[\ell_{\infty}(\alpha)=\varepsilon_0=2\operatorname{arcsinh}(1)~.\]
We now consider $u_3$. By the Schwarz Lemma we have $h_{\infty}\leq h_{Th}$ and so $u_3\leq 0$ everywhere. For a point $p\in \alpha_d$, we denote by $\operatorname{inj}(p)$ the injectivity radius at $p$ with respect to $h_{\infty}$. Then, for any $\epsilon\leq \operatorname{inj}(p)$, by \cite[Theorem 2.8]{BBB2018} and \cite[Corollary 2.12]{BBB2018}, we have:
\begin{equation}\label{eq:thurstoneq}
h_{Th}\leq h_{\infty}(1+3\coth^2(\epsilon/2))\coth^2(R^{\varepsilon}_d/2)~,\end{equation}
 where $R^{\varepsilon}_d$ is such that the injectivity radius inside a ball of radius $R^{\varepsilon}_d$ centered at any $p\in\alpha_d$ is at least $\epsilon$. The injectivity radius, see \cite[Thm 4.1.6]{Bu1992}, at a
 point $p$ at distance $d$ from $\partial \mathbb{T}(\ell)$ satisfies
 \[\sinh(\operatorname{inj}(p))=\sinh(\ell/2)\cosh(L-d)~,\]
 then, the condition on $R^{\varepsilon}_d$, being $\sinh(x)$ increasing for $x>0$, becomes
 \begin{equation}\label{eq:thurstoneq2} \sinh(\varepsilon)\leq \sinh(\ell/2)\cosh(L-d-R^{\varepsilon}_d)~.\end{equation}
Since $\coth^2(x/2)$ is always greater than one, with vertical asymptote at $0$  and decreasing for $x>0$, we need to find suitable $\epsilon$, $d$, and $R^{\epsilon}_d$, as in equation \eqref{eq:thurstoneq} and \eqref{eq:thurstoneq2}, and bound it from below. In order to prove the estimate for $\abs{u_1}$ on $\partial \mathcal{A}$ and $\alpha$, we are interested in the two cases $d=0$ and $d=d_{\epsilon_0}$ respectively.

\textbf{Case $d=0$.} The injectivity radius at any point of the boundary $\partial \mathbb{T}(\ell)$ is at least $\varepsilon_0/2=\operatorname{arcsinh}(1)$, so we take $\varepsilon=\varepsilon_0/4$ and define our radius to be
\[R_0=\inf_{\ell\in(0,\epsilon_0]} \left(\operatorname{arcsinh}\left(\f{1}{\sinh(\ell/2)}\right)-\operatorname{arcosh}\left(\f{\sinh(\varepsilon_0/4)}{\sinh(\ell/2)}\right)\right)~,\]
as the function we are taking the infimum of is increasing in $\ell$, the expression is given by its limit as $\ell\to0$ which is greater than $\pi/4$. Then, $\coth^2(R_0/2)$ is bounded above by $7.2$. Therefore, since

\[\abs{u_{3\vert\partial \mathcal{A}}}\leq \frac{1}{2}\log\left(\left(1+3\coth^2(\varepsilon_0/8)\right)\coth^2(R_0/2)\right),\]
we get
\[-3.08\leq u_{3\vert\partial \mathcal{A}}\leq 0~.\]
As $u_1=-( u_2+u_3+u_0)$, by merging the estimates:
\[ -b(\ell_{\infty}(\alpha_0))-\log\left(\dfrac{2\pi}{\ell_{\infty}(\alpha_0)}\right)\leq u_{1\vert \partial \mathcal{A}}\leq 3.08-\log\left(\dfrac{2\pi}{\ell_{\infty}(\alpha_0)}\right)~.\]
Recalling the definition of the function $b(x)$ from Lemma \ref{horosphere}, and since $\ell_{\infty}(\alpha_0)\in [2, 2.5]$, see Remark \ref{alphazero}, the worst bounds for the previous equation are obtained by evaluating the terms at $2.5$, hence
\[-3.7 \leq u_{1\vert \partial \mathcal{A}}\leq 2.16~.\]
\textbf{Case $d=d_{\epsilon_0}$.} Since $\cosh(L-d_{\varepsilon_0})=\varepsilon_0/\ell$ the injectivity radius of the points $p$ in $\alpha$ is \[\operatorname{inj}(p)=\operatorname{arcsinh}(\sinh(\ell/2)\varepsilon_0/\ell)~,\]
then, by choosing $\varepsilon=\operatorname{inj}(p)/2$, we take
 \[R_{d_{\varepsilon_0}}=\inf_{\ell\in(0,\epsilon_0]}\left(  \operatorname{arcosh}(\varepsilon_0/\ell)-\operatorname{arcosh}\left(\f{\sinh(\frac{1}{2}\operatorname{arcsinh}(\sinh(\ell/2)\varepsilon_0/\ell))}{\sinh(\ell/2)}\right)\right)~,\]
 which is greater than $0.77$. Then, $\coth^2(R_{d_{\varepsilon_0}}/2)$ is bounded above by $7.423$. Thus, \[\abs{u_{3\vert\alpha}}\leq \frac{1}{2}\log\left(\left(1+3\coth^2(\operatorname{inj}(p)/4)\right)\coth^2(R_{d_{\varepsilon_0}}/2)\right)~,\]
and \[\operatorname{inj}(p)=\operatorname{arcsinh}(\sinh(\ell/2)\varepsilon_0/\ell)\geq \operatorname{arcsinh}(\varepsilon_0/2)~.\]
Using these estimates and the fact that $u_3\leq 0$ we obtain:
\[-3.187\leq u_{3\vert_{\alpha}}\leq 0~.\]

To bound $u_1\vert_{\alpha}$ we argue as in the previous case and we obtain:
\[ -b(\varepsilon_0)-\log\left(\dfrac{2\pi}{\varepsilon_0}\right)\leq u_{1\vert_{ \alpha}}\leq 3.187-\log\left(\dfrac{2\pi}{\varepsilon_0}\right)\] from which 
\[-2.3 \leq u_{1\vert_{ \alpha}}\leq 1.92~.\]
Defining $W$ and $W_0$ as the evaluation of the function $b(x)+\log(2\pi/x)$ at respectively $\ell_{\infty}(\alpha_0)$ and $\varepsilon_0$ we get the statement.\epf

Now we want to solve the following PDE on $\mathcal{A}\cong  [-m,m]\times\mathbb S^1$:
\[ \Delta u=0 \qquad \abs{ u\vert_{\partial \mathcal{A}}}\leq W\leq 3.7~,\]
where $m$ is half of the length of the thin tube as in Proposition \ref{bounB2} and $W$ is as in Lemma \ref{lem1}.

Using Fourier analysis, the solution is of the form
\[u(r,\theta)=\sum_{k\in\Z} u_k(r) e^{ik\theta}~.\]
Then, we get that the coefficients need to satisfy:
\[ \ddot u_k(r)=k^2 u_k(r),\qquad \abs{u_k(\pm m)}\leq W~.\]
This yields solutions of the form:
\[k=0: u_0(r)=C_0+C_1r~,\]
and:
\[k\neq0: u_k(r) = v_1(k) \cosh(kr) + v_2(k) \sinh (kr)~.\]
Then, we get the following linear system of equations:
\be
\begin{cases} 
&v_1(k)\cosh (km)+v_2(k)\sinh(km)=u_k(m)\\
& v_1(k)\cosh (km)-v_2(k)\sinh(km)=u_k(-m)\\
\end{cases} 
\ee 
which can be rewritten as:
\be
\begin{cases} 
&v_1(k)=\f{u_k(m)+u_k(-m)}2 (\cosh(km))^{-1}\\
& v_2(k)=\f{u_k(m)-u_k(-m)}2 (\sinh(km))^{-1}\\
\end{cases} 
\ee 
which yields for $k\neq0$, since $e^{-2\abs km}<e^{-m}\leq\f 1{16} < \f 12$ and $W\leq 3.7$:
\begin{equation}\label{eq:four1} \abs{v_1(k)} \leq2 W e^{-\abs km}\leq 7.4 e^{-\abs km} \qquad \abs{v_2(k)} \leq 4 W e^{-\abs km}\leq 14.8 e^{-\abs km}~.\end{equation}
 We now show that $C_1=0$ with a separate geometric lemma:

\blem In the setting above $C_1=0$.\elem
\bpf Let $\Phi\colon (-m,m)\times(\mathbb{R}/2\pi\mathbb{Z})
\longrightarrow \mathcal{A}\subset\mathbb{C}^{*}$ be our conformal diffeomorphism, 
unwrap the annulus to its universal covering and let $p$ be the covering map. Let $\widetilde{\Phi}\eqdef \Phi\circ p$. Thus $\widetilde{\Phi}$ is holomorphic and $2\pi i$-periodic. Since $\Phi$ is a conformal diffeomorphism and its image is contained in $\mathbb{C}^{*}$, both $\widetilde{\Phi}$ and
$\widetilde{\Phi}'$ are nonzero. Pulling back $h_{0}$ gives $\widetilde{\Phi}^{*}h_{0}
=
\frac{|\widetilde{\Phi}'(w)|^{2}}
     {|\widetilde{\Phi}(w)|^{2}}
|dw|^{2}$. Consequently,
\[
u(r,\theta)
=
\log\left|\widetilde{\Phi}'(r+i\theta)\right|
-
\log\left|\widetilde{\Phi}(r+i\theta)\right|.
\]
Fix $r\in(-m,m)$ and consider the closed curve  $
\alpha_{r}(\theta)
\eqdef
\widetilde{\Phi}(r+i\theta)$, with $
\theta\in[0,2\pi]$. Differentiating the zero Fourier coefficient, we obtain
\begin{align*}
C_{1}
=
\frac{1}{2\pi}
\int_{0}^{2\pi}u_{r}(r,\theta)\,d\theta
=
\frac{1}{2\pi}
\operatorname{Re}
\int_{0}^{2\pi}
\left(
\frac{\widetilde{\Phi}''}{\widetilde{\Phi}'}
-
\frac{\widetilde{\Phi}'}{\widetilde{\Phi}}
\right)(r+i\theta)\,d\theta.
\end{align*}
Since $\dot{\alpha}_{r}(\theta)
=
i\widetilde{\Phi}'(r+i\theta)$ and
$
\ddot{\alpha}_{r}(\theta)
=
-\widetilde{\Phi}''(r+i\theta),
$
the winding-number formula gives
\begin{align*}
\operatorname{wind}(\alpha_{r},0)
&=
\frac{1}{2\pi i}
\int_{0}^{2\pi}
\frac{\dot{\alpha}_{r}(\theta)}
     {\alpha_{r}(\theta)}\,d\theta
=
\frac{1}{2\pi}
\int_{0}^{2\pi}
\frac{\widetilde{\Phi}'}
     {\widetilde{\Phi}}(r+i\theta)\,d\theta,
\\
\operatorname{wind}(\dot{\alpha}_{r},0)
&=
\frac{1}{2\pi i}
\int_{0}^{2\pi}
\frac{\ddot{\alpha}_{r}(\theta)}
     {\dot{\alpha}_{r}(\theta)}\,d\theta
=
\frac{1}{2\pi}
\int_{0}^{2\pi}
\frac{\widetilde{\Phi}''}
     {\widetilde{\Phi}'}(r+i\theta)\,d\theta.
\end{align*}
Both integrals are integers and hence real. Therefore, $C_{1}= \operatorname{wind}(\dot{\alpha}_{r},0)
- \operatorname{wind}(\alpha_{r},0)$. The curve $\alpha_{r}$ is a regular simple essential closed curve in
$\mathbb{C}^{*}$. Thus $0$ lies in its bounded complementary component,
and $\operatorname{wind}(\alpha_{r},0)=\pm1$. By the turning-number theorem for regular simple closed plane curves, the tangent vector makes one full turn with the same orientation: $
\operatorname{wind}(\dot{\alpha}_{r},0)=\operatorname{wind}(\alpha_{r},0)$. Thus $C_1=0$, as claimed.\epf

Thus, from Lemma \ref{lem1} we also get $\abs{v_1(0)}=\abs{C_0}\leq W _0\leq 2.3$. Similar to equation \eqref{eq:four1} one can obtain bounds on the $n$-derivatives of $u(r,\theta)$ of the form $C\abs k^n e^{-\abs k(m-\abs r)}$, as:
\begin{align*}n\text{ even:}\quad &\abs {\difrac {^n u_k (r)}{r^n}}= \abs{v_1(k)k^n\cosh(kr)+v_2(k) k^n\sinh(kr)},  \\
n\text{ odd:}\quad &\abs {\difrac {^n u_k (r)}{r^n}}= \abs{v_1(k)k^n\sinh(kr)+v_2(k) k^n\cosh(kr)},   \end{align*}
and the same computation gives bounds of the order of $\abs k^n e^{-\abs k(m-\abs r)}$ which still decays exponentially. Thus, we get that the decay of the coefficients of $u(r,\theta)$ and its $n$-th derivatives is of the order $\abs k^ne^{-\abs k(m-\abs r)}$.\

To sum-up we have
\[ u(r,\theta)= C_0 +\sum_{\abs{k}>0}u_k(r)e^{ik\theta}~,\]
with coefficients:
 \[ u_k(r) = v_1(k) \cosh(kr) + v_2(k) \sinh (kr)~\quad u_0=C_0,\]
satisfying
 \[ \abs{v_1(k)} \leq 2W e^{-\abs  km} \qquad \abs{v_2(k)} \leq 4W e^{-\abs km}\quad \abs {C_0}\leq W_0~,\]
where $W$ and $W_0$ are as in Lemma \ref{lem1}.

 \subsubsection{Estimates}\label{est}
 We can now compute and bound the first and second derivatives of $u$, which, when $u=u_1$, allow to bound the Osgood-Stowe tensor $\B{h_{\infty,0}} {u_1}$ for the two flat metrics of the previous sections. As we have already seen, it is also necessary to have a good bound on $u_1$ at the core curve. We will obtain it at the end of this subsection. From now on, to simplify the notation, we will drop the subscript $1$. Then, using the previous discussion we have:
\begin{itemize}
\item $u_r(r,\theta)=\sum_{\abs{k}>0} k\left(v_1(k)\sinh(kr)+v_2(k)\cosh(kr)\right)e^{ik\theta}$,
\item $u_\theta(r,\theta)=\sum_{\abs{k}>0}iku_k(r)e^{ik\theta}$,
\item $u_{\theta r}(r,\theta)=\sum_{\abs{k}>0} ik^2\left(v_1(k)\sinh(kr)+v_2(k)\cosh(kr)\right)e^{ik\theta}$,
\item $u_{rr}(r,\theta)=\sum_{\abs{k}>0} k^2\left(v_1(k)\cosh(kr)+v_2(k)\sinh(kr)\right)e^{ik\theta}$,
\item $u_{\theta \theta}(r,\theta)=-\sum_{\abs{k}>0}k^2 u_{k}(r)e^{ik\theta}$.
\end{itemize}

Using the bounds on $v_1(k)$ and $v_2(k)$ we obtain the following estimates:
\begin{itemize}
\item $\abs{u_r(r,\theta)}\leq 2\sum_{k>0}k\left( 2W e^{-km}\abs{\sinh(kr)}+ 4 W e^{-km}\cosh(kr)\right)$,
\item $\abs{u_\theta(r,\theta)}\leq 2\sum_{k>0} k \left( 2 W e^{-km}\cosh(kr)+4 W e^{-km}\abs{\sinh(kr)}\right)$,
\item $\abs{u_{\theta r}(r,\theta)}\leq2\sum_{k>0} k^2\left( 2 W e^{-km}\abs{\sinh(kr)}+4 W e^{-km}\cosh(kr)\right)$,
\item $\abs{u_{rr}(r,\theta)}\leq 2\sum_{k>0} k^2\left( 2 W e^{-km}\cosh(kr)+4 W e^{-km}\abs{\sinh(kr)})\right),$
\item $\abs{u_{\theta\theta}(r,\theta)}\leq2\sum_{k>0} k^2\left( 2We^{-km}\cosh(kr)+4We^{-km}\abs{\sinh(kr)}\right)$.
\end{itemize}
If we just care about the geodesic at $r=0$ (which is the image of $\gamma$ in the flat coordinates) we have the following estimates:
\begin{itemize}
\item $\abs{u_r(0,\theta)}\leq   8 W \sum_{k>0} k e^{-k m}=8 W\f{e^{m}}{(e^{m}-1)^2}\leq 30\f{e^{m}}{(e^{m}-1)^2} $,
\item $\abs{u_\theta(0,\theta)}\leq 4 W \sum_{k>0} k e^{-km}=4W\f{e^m}{(e^m-1)^2}\leq 15\f{e^m}{(e^m-1)^2}$,
\item $\abs{u_{\theta r}(0,\theta)}\leq 8 W \sum_{k>0} k^2  e^{-km}= 8W\f{e^{m} (e^{m}+1)}{(e^{m}-1)^3}\leq 30\f{e^{m} (e^{m}+1)}{(e^{m}-1)^3}$,
\item $\abs{u_{rr}(0,\theta)}\leq 4 W \sum_{k>0} k^2  e^{-km}=4 W\f{e^{m} (e^{m}+1)}{(e^{m}-1)^3}\leq 15\f{e^{m} (e^{m}+1)}{(e^{m}-1)^3}$,
\item $\abs{u_{\theta\theta}(0,\theta)}\leq 4W \sum_{k>0} k^2  e^{-km}=4W\f{e^{m} (e^{m}+1)}{(e^{m}-1)^3}\leq 15\f{e^{m} (e^{m}+1)}{(e^{m}-1)^3}$.
\end{itemize}
In the more general case:
\begin{itemize}
\item $\abs{u_r(r,\theta)}\leq   10 W \sum_{k>0} k e^{-k(m-\abs r)}=10 W\f{e^{(m-\abs r)}}{(e^{(m-\abs r)}-1)^2} $,
\item $\abs{u_\theta(r,\theta)}\leq 8 W \sum_{{k}>0} k e^{-k(m-\abs r)}=8 W\f{e^{(m-\abs r)}}{(e^{(m-\abs r)}-1)^2}$,
\item $\abs{u_{\theta r}(r,\theta)}\leq 10 W \sum_{{k}>0} k^2  e^{-k(m-\abs r)}=10W\f{e^{(m-\abs r)} (e^{(m-\abs r)}+1)}{(e^{(m-\abs r)}-1)^3}$,
\item $\abs{u_{rr}(r,\theta)}\leq 8 W \sum_{{k}>0} k^2  e^{-k(m-\abs r)}=8 W\f{e^{(m-\abs r)} (e^{(m-\abs r)}+1)}{(e^{(m-\abs r)}-1)^3}$,
\item $\abs{u_{\theta\theta}(r,\theta)}\leq 8W\sum_{{k}>0} k^2  e^{-k(m-\abs r)}=8W\f{e^{(m-\abs r)} (e^{(m-\abs r)}+1)}{(e^{(m-\abs r)}-1)^3}$.
\end{itemize}
Note that, in the notation of the previous sections, for any fixed $r$, when $2m$ is the length of the Euclidean tube $(\mathcal{A}, h_{\infty, 0})$, the last two blocks of estimates are in $O(e^{-m})\leq O(e^{-\pi^2/(2\ell)})$. 

\brem\label{cstcorebound}
At the core $\gamma$:
\begin{align*}
\abs{u(\gamma)}=\abs{u(0,\theta)}&=\abs{C_0+\sum_{\abs k>0} v_1(k)e^{ik\theta}}\leq\abs {C_0}+ \sum_{\abs k>0}  \abs{v_1(k)}\leq W_0 +4W\sum_{k=1}^\infty  e^{-km}\\
&\leq W_0+\f{4W}{e^m-1}\leq W_0+\f{4W}{\exp\left(\f{2\pi\arctan(\sinh(L))}\ell\right)-1}\leq 2.3+1=3.3~,\\
\end{align*}
where the last inequality is obtained by substituting $L=\argsinh(1/\sinh(\varepsilon_0/2))$, $W\leq 3.7$ and $W_0\leq 2.3$. Also, we see that for $\ell\rar 0$ the bound goes to $W_0$ as
\begin{align*}
\abs{u(\gamma)}&\leq W_0+ \f{4We^{-\pi^2/(2\ell)}}{1-e^{-\pi^2/(2\ell)}}\\
&\leq W_0+4.4 W e^{-\pi^2/(2\ell)}~.
\end{align*}
\erem

We now want to improve the estimate for $u$ on the core (and on a neighborhood) also by using the bound of Lemma \ref{almostround}:
\[ \max_{s\in[0,\ell]} \norm{ D^0_{\dot\gamma}\dot\gamma(s)}_{h_0} \leq 12We^{3.3}\left( \f{2\pi}\ell\right)^2 \f{e^m}{(e^m-1)^2},\quad W\leq 3.7~,\]

\blem \label{estcore}
On the core $\gamma$ the function $u$ satisfies:
\[ e^u\leq  1+4W\pi(3e^{3.3}+ 2)\f{e^m}{(e^m-1)^2}  +  96 W^2e^{3.3}\pi^2\left(  \f{e^m}{(e^m-1)^2}\right)^2.\]
In particular, for small $\ell$ we have $e^u=1+O\left(e^{-\pi^2/(2\ell)}\right)$.\elem

\bpf
The core is parameterized as $\gamma(s)=\left(0,\f {2\pi}\ell s\right)$ with $s\in[0,\ell]$ with respect to the coordinates $(r, \theta)$ of $h_{\infty,0}$, and $h_0=e^{2u}h_{\infty,0}$. Then, by using that $u(s)=u\left(0,\f {2\pi}\ell s\right)$ we have:
\[ 2\pi\leq\ell_{h_0}(\gamma)=\int_0^\ell\norm{\dot\gamma}_0ds=\int_0^\ell e^{u(s)}\norm{\dot\gamma}_{h_{\infty,0}}ds=\f{2\pi}\ell \int_0^\ell e^{u(s)}ds~.\]

We now estimate $\ell_{h_0}(\gamma)$ by using the mean $\bar\gamma'=(0,\f{2\pi}\ell)$ of $\dot\gamma$ as seen with respect to $h_0$, as we did in Lemma \ref{bounB3}.
\begin{align*}
\int_0^\ell \norm{\dot\gamma}_{h_0}ds&\leq\int_0^\ell  \norm{\dot\gamma-\bar\gamma'}_{h_0}ds+\int_0^\ell  \norm{\bar\gamma'}_{h_0}ds\\
&\leq \int_0^\ell\f{  24 We^{3.3}\pi^2}\ell\f{e^m}{(e^m-1)^2}ds+2\pi\\
&\leq 2\pi+ 24 We^{3.3}\pi^2\f{e^m}{(e^m-1)^2}~.
\end{align*}
Here the second inequality follows from Lemma \ref{almostround} and $ \norm{\dot\gamma-\bar\gamma'}_{h_0}\leq\f \ell 2\max_{s\in[0,\ell]} \norm{ D^0_{\dot\gamma}\dot\gamma(s)}_{h_0}$. Thus,
\[2\pi\leq \f{2\pi}\ell \int_0^\ell e^{u(s)}ds\leq  2\pi+ 24 We^{3.3}\pi^2\f{e^m}{(e^m-1)^2}~,\]
which gives:
\[ \ell\leq  \int_0^\ell e^{u(s)}ds\leq \ell+ 12 We^{3.3}\pi \ell \f{e^m}{(e^m-1)^2}~.\]
Therefore, there must be $\xi\in[0,\ell]$ such that 
\[e^{u(\xi)}\leq 1+ 12 We^{3.3}\pi\f{e^m}{(e^m-1)^2}~.\] 
For all $s\in[0,\ell]$ we have:
\begin{align*}
\abs{e^{u(s)}-e^{u(\xi)}}&=\abs{\int_\xi^s e^{u(t)}\dot u(t)dt}\\
&\leq \int_\xi^s e^{u(t)}\max_{s\in[0,\ell]}\abs{\dot u(s)} ds\\
&\leq 4W\left(\f{2\pi}\ell\right)\f{e^m}{(e^m-1)^2} \int_\xi^s e^{u(t)}dt\\
&\leq 4W\left(\f{2\pi}\ell\right)\f{e^m}{(e^m-1)^2} \int_0^\ell e^{u(t)}dt\\
&\leq 4W\left(\f{2\pi}\ell\right)\f{e^m}{(e^m-1)^2}\left(  \ell+ 12 We^{3.3}\pi \ell \f{e^m}{(e^m-1)^2}\right)\\
&=8\pi W\f{e^m}{(e^m-1)^2}\left( 1+ 12 We^{3.3}\pi  \f{e^m}{(e^m-1)^2}\right),
\end{align*}
and by using $e^{u(\xi)}\leq 1+12 We^{3.3}\pi\f{e^m}{(e^m-1)^2}$ and re-arranging the terms we obtain the required statement. Analogously, there exists $\xi'\in[0,\ell]$ such that $e^{u(\xi')}\geq 1$, and we can apply the same argument.
The fact that $e^u=1+O(e^{-\pi^2/(2\ell)})$ follows directly by Equation \eqref{mest} in Remark \ref{mbound}.
 \epf

\brem \label{Mell} The bound of Lemma \ref{estcore} beats the one of Remark \ref{cstcorebound} when $\ell$ is small enough. Then, we can write:
\[ \max_{s\in[0,\ell]} \norm{ D^0_{\dot\gamma}\dot\gamma(s)}_{h_0} \leq 12W G(\ell) \left( \f{2\pi}\ell\right)^2 \f{e^m}{(e^m-1)^2},\quad W\leq 3.7,\]
with:
\[ G(\ell)\eqdef \min \set{e^{3.3},   1+4W\pi(3e^{3.3}+ 2)\f{e^m}{(e^m-1)^2}  +  96 W^2e^{3.3}\pi^2\left(  \f{e^m}{(e^m-1)^2}\right)^2}.\]
By Equation \eqref{mest} in Remark \ref{mbound}, $G(\ell)$ is dominated by:
\[1+4W\pi(3e^{3.3}+ 2) \left( 16/15\right)^2 e^{-\pi^2/(2\ell)}+  96W^2e^{3.3}\pi^2  \left( 16/15\right)^4e^{-\pi^2/\ell}\]
which is bounded and of the form  $1+O\left(e^{-\pi^2/(2\ell)}\right)$ for small $\ell$.

\erem

 \brem \label{utube}
Lemma \ref{estcore} together with the last block of estimates for the partial derivatives of $u$ above, imply that the same bound for $e^{u}$ still holds in a neighborhood of the core $\gamma$ in $\mathcal{A}$. The estimates themselves will depend on the chosen neighborhood $U_r=(-r,r)\times\mathbb S^1$ and will be of the form $1+C e^{r-m}$, for $C$ a universal constant. However, as we will not need it, for concreteness, we work in $U=U_1$.

Specifically, by using the Taylor expansion at $(0, \theta_0)$, to first order with the formula of the Lagrange error in $U$, we obtain 
\[
e^{u(r,\theta)}
=e^{u(0,\theta_0)}
+e^{u(0,\theta_0)}u_r(0,\theta_0)\,r
+e^{u(0,\theta_0)}u_\theta(0,\theta_0)\,(\theta-\theta_0) + R_\xi(r,\theta)\]
with
\[ R_\xi(r,\theta)=
\frac12 e^{u(\xi)}\Big[
\big(u_r(\xi)^2+u_{rr}(\xi)\big)r^2
+2\big(u_r(\xi)u_\theta(\xi)+u_{r\theta}(\xi)\big) r(\theta-\theta_0)
+\big(u_\theta(\xi)^2+u_{\theta\theta}(\xi)\big)(\theta-\theta_0)^2
\Big],
\]
for a certain point $\xi\in U$.
We then use the bounds on the derivatives of $u$ coming from the Fourier analysis. Assuming $\abs{r}\leq 1$, since also $e^x/(e^x-1)^2\leq 4e^{-x}$ for any $x>\log(2)$ and $m-\abs r\geq m-1\geq \log(16)-1>\log(2)$, we have the following bounds:

\begin{itemize}
\item $\abs{u_r}$, $\abs{u_\theta}\leq 40e  W e^{-m}$,
\item  $\abs{u_{rr}}$, $\abs{u_{r\theta}},\abs{u_{\theta\theta}} \leq 64e W  e^{-m}$.
\end{itemize}

Thus, we can bound the remainder term as
\begin{align*}
\abs{R_\xi(r,\theta)}
&\leq
\frac12 e^{u(\xi)}
\left(
1600e^2W^2e^{-2m}+64eWe^{-m}
\right)
\left(\abs r+\abs{\theta-\theta_0}\right)^2\\
&\leq
\frac12 e^W(1+2\pi)^2
\left(
1600e^2W^2e^{-2m}+64eWe^{-m}
\right),
\end{align*}
where in the second inequality we used $\abs r\leq1$,
$\abs{\theta-\theta_0}\leq2\pi$, and the maximum principle,
which gives $\abs{u}\leq W$ on $\mathcal A$. Using $e^{-m}\leq e^{-\pi^2/(2\ell)}$, we obtain
\[
\abs{R_\xi(r,\theta)}
\leq
e^{-\pi^2/(2\ell)}D(\ell),
\]
where
\[
D(\ell)
\eqdef
\frac12 e^W(1+2\pi)^2
\left(
64eW
+
1600e^2W^2e^{-\pi^2/(2\ell)}
\right).
\]
Since $W\leq3.7$, the function $D(\ell)$ is uniformly bounded. Therefore, from the Taylor expansion,
\begin{align*}
\abs{e^{u(r,\theta)}-1}
&\leq \abs{e^{u(0,\theta_0)}-1}+e^{u(0,\theta_0)}
\left( 8W\abs r + 4W\abs{\theta-\theta_0}
\right)\frac{e^m}{(e^m-1)^2} +\abs{R_\xi(r,\theta)}\\
&\leq
\abs{e^{u(0,\theta_0)}-1}
+
e^{u(0,\theta_0)}
8W(1+\pi)\frac{e^m}{(e^m-1)^2}+
e^{-\pi^2/(2\ell)}D(\ell).
\end{align*}

By Lemma \ref{estcore},
\[
\abs{e^{u(0,\theta_0)}-1}
\leq
12\pi We^{3.3}\frac{e^m}{(e^m-1)^2}
+
8\pi W \frac{e^m}{(e^m-1)^2}+
96\pi^2 W^2e^{3.3} \frac{e^{2m}}{(e^m-1)^4}.
\]
Moreover, by Remark \ref{mbound}, $
\frac{e^m}{(e^m-1)^2}
\leq \left(\frac{16}{15}\right)^2e^{-m}
\leq \left(\frac{16}{15}\right)^2
e^{-\pi^2/(2\ell)}$. It follows that, uniformly on
$U=(-1,1)\times\mathbb S^1$,
\[
e^{u(r,\theta)}
=
1+O\left(e^{-\pi^2/(2\ell)}\right).
\]\erem

  \section{Applications to Renormalized Volume}
\label{sc:appl}

\subsection{Renormalized volume bounds under earthquakes on short compressible curves} In this section we prove the version of Theorem \ref{mainVR} with an explicit formula for $F(\ell)$. Note that the complex projective structures naturally arising on the conformal boundary of convex co-compact hyperbolic $3$-manifolds $M$ are realized as quotients of invariant proper domains $\Omega$ of $\mathbb{CP}^1$ by Kleinian groups. In particular, the associated developing maps are coverings, and, when the boundary of $M$ is compressible, their image is a non-simply-connected hyperbolic domain in $\mathbb{CP}^1$. Then, Theorems \ref{Schwarzian} and \ref{pairing} hold for any short compressible simple closed curve in the conformal boundary at infinity.

\begin{customthm}{\ref{mainVR}}
	Let $M$ be a convex co-compact hyperbolic $3$-manifold. Let $X_0\in \mathcal{T}(\partial\overline{M})$, and let $X_t\in \mathcal{T}(\partial\overline{M})$ be the Riemann surface obtained by a parameter $t\in\mathbb{R_+}$ earthquake on $X_0$ along a simple closed compressible geodesic $\gamma$ in $\partial\bar{M}$ of length $\ell_X(\gamma)\leq\epsilon_0$. Then, the following estimate for the renormalized volume of the associated convex co-compact manifolds $M_0$ and $M_t$ holds:
	\[\abs{V_R(M_t)-V_R(M_0)}\leq 
148\pi^2 \f{e^{-\pi^2/(2\ell)}}\ell t + 2553 \pi^4 \f{G (\ell)^2}\ell  e^{-\pi^2/\ell}t~,\]
	with $G (\ell)\rar 1$ as $\ell\rar 0$ and always bounded by $e^{3.3}$.
\end{customthm}

\begin{proof}
	We first recall that the differential of the renormalized volume $d V_R$ at a point $X_s$ coincides with $\re(S(f_{M_s}))$, with $M_s$ the convex co-compact hyperbolic $3$-manifold associated to $X_s$, and $f_{M_s}$ the developing map of its associated complex projective structure. Let us define the earthquake path $X_s=eq_s(X_0)$, with $s\in[0,t]$, and write the difference of renormalized volumes as \[V_R(M_t)-V_R(M_0)=\int_0^t (dV_R)_{X_s}(\mu_s)ds~,\]
	with $\mu_s$ the Beltrami differential obtained by differentiating the earthquaking path at time $s$, that is the infinitesimal earthquake at $X_s$.
	Then, the result is a direct corollary of Theorem \ref{mainthm} in the version of Section \ref{ssc:estimatesearth}.
\end{proof}

\subsection{Asymptotic behaviour of $V_R$ along pinching a compressible curve}

In this section we prove Theorem \ref{asymptotic}.

Differently from the previous section, in which to prove Theorem \ref{mainVR} we just integrated the differential given by the first equality in Theorem \ref{pairing}, we now have to take care of two facts. Firstly, the length of the simple closed geodesic on which we are grafting changes its length along the deformation path. Secondly, while earthquaking on the same geodesic produces a flow, grafting does not.

We will consider the path $X_s=gr_{s\gamma}(X_0)$ obtained by grafting on the short compressible simple closed curve $\gamma$ in $X_0=\partial_{\infty}M_0$, and denote by $(M_s)_{s\in[0, \infty)}$ the path of convex co-compact hyperbolic $3$-manifolds associated to $X_s$ through the Uniformization Theorem (see Section \ref{hypback}). We will consider here the composition $\ell_{\cdot}(\gamma)$ of the length function 
\[\ell_\gamma \colon \mathcal{T}(\partial\overline{M})\rightarrow \mathbb{R}^{+}~,\] 
with the path $(X_s)_{s\in[0,\infty)}$ which associates to $s$ the length $\ell_{s}(\gamma)$ of $\gamma$ with respect to the hyperbolic representative in $X_s$. For the sake of notation, whenever the dependence on $\gamma$ is clear, we will just write $\ell_s=\ell_s(\gamma)$. We recall that as $s$ goes to infinity $\ell_{s}(\gamma)$ tends to $0$.

The following lemma was obtained by Diaz and Kim in \cite[Proposition $3.4$]{graftingrays} (or also in \cite[Lemma $4.1$]{henselgraftingrays}). Our contribution is just to explicitly write down the constant appearing in the lower bound.

\blem\label{graftlength}
Let $X_0\in \mathcal{T}(\partial\overline{M})$, and let $X_s\in \mathcal{T}(\partial\overline{M})$, with $s>0$, be the Riemann surface obtained by a parameter $s$ grafting along a short geodesic $\gamma$ of length $\ell_{0}\leq \varepsilon_0$. Then,
 \[\f{ \pi}{\pi+s} \ell_{0}\geq \ell_{s}\geq\f{\pi}{2(\pi+s)}\ell_{0}~.\] 
\elem

\bpf
The first inequality was proven in \cite[Proposition $3.4$]{graftingrays}, by a direct computation applying the definition of grafting at the universal cover (see Subsection \ref{complexearth}). In the same statement, we can also find the lower bound 
\[\ell_{s}\geq \f{2\theta_0}{2\theta_0+s}\ell_0~,\]
where $\theta_0$ is such that the thin tube around $\gamma$ (see Definition \ref{margulitubedefi}) in $X_0$ can be isometrically lifted in $\mathbb{H}^2$ to a fundamental domain \[\{\rho e^{i\theta}\in\mathbb{H}^2\ |\ \rho\in[1, e^{\ell}],\ \theta\in[\pi/2-\theta_0, \pi/2+\theta_0]\}~.\] 
By elementary hyperbolic geometry (see for example \cite[Lemma $5.2.7$]{Mar16}), the angle $\theta_0$ satisfies: \[\theta_0=\argtg(\sinh(L_0))=\argtg\left(\f{1}{\sinh(\ell_0/2)}\right)~,\]
where the last equality follows from the fact that the width $L$ of the thin tube of core length $\ell$ is \[L=\argsinh\left(\f1{\sinh(\ell/2)}\right)~.\] 
Since now $\ell_0\in[0, \varepsilon_0]$, with $\varepsilon_0=2\argsinh(1)$, we can estimate $\theta_0$ \[\pi/4\leq\theta_0\leq\pi/2~,\]
from which the right inequality of our statement follows.\epf

\brem
The upper bound in Lemma \ref{graftlength} is asymptotically sharp, see \cite[Theorem $6.6$]{DumWolf}.
\erem

We define:   
\[\nu_0^s=\left.\frac{d}{dt}\right|_{t=0}gr_{t\gamma}(X_s)~.\]
and 
\[\nu_s^0=\f{d}{dt}\bigg\vert_{t=s}gr_{t\gamma}(X_0)~.\]
Note that both Beltrami differentials belong to the tangent space $T_{X_s}\mathcal{T}(\partial\bar{M})$, and that the first one corresponds to the infinitesimal grafting at $X_s$, while the second one is the derivative of the grafting path $gr_{t\gamma}(X_0)$ at $t=s$, and they need not coincide for $s>0$ but they are equal at $s=0$. For suitable bounded Beltrami representatives, by explicit computation, we have that their $L^\infty$ norm is bounded as:

\[
\|\nu_s^0\|_{\infty}
\leq \frac{1}{2s+\pi},
\qquad
\|\nu_0^s\|_{\infty}
\leq \frac{1}{\pi},
\qquad s\geq 0.
\]

\blem\label{Gardiner}
In the notations above:
\[(d V_R)_{X_s}(\nu_s^0)=\left(\f{\pi}{4}+\f{\pi^3}{\ell_s^2}\right)d\ell_{\gamma}(\nu_s^0)+E(s)~,\]
where $E(s)$ stands for a bounded real function such that $E(s)\in O\left(\f{e^{-\pi^2/(2\ell_s)}}{\ell_s^2}\right)$ for large $s$.
\elem

\begin{proof}
	Let us denote by $z=x+iy$ the complex coordinate of the half-space model $\mathbb{H}^2$. Note that the holomorphic quadratic differential $dz^2/z^2$ is invariant under the action of the hyperbolic isometry $\varphi_{\gamma, s}(z)=e^{\ell_s}z$ relative to the simple closed geodesic $\gamma$ in $X_s$. By Gardiner's formula for the differential of the length function of $\gamma$ (see \cite{Gard}):
	\begin{equation}\label{dl} d\ell_{\gamma}(\mu)= \f{2}{\pi}\re\left\langle \f{dz^2}{z^2}, \widetilde{\mu }\right\rangle_{\mathcal{A}_{\gamma,s}}\end{equation}
	where $\mu$ is a harmonic Beltrami differential in $T_{X_s}\mathcal{T}(\partial\bar{M})$, and $\widetilde{\mu}$ denotes its lift to the annulus $\mathcal{A}_{\gamma,s}=\mathbb{H}^2/\langle \varphi_{\gamma,s}\rangle$. If we take $\mu=\nu_s^0$, being this equal to zero out of a tubular neighborhood of $\gamma$, the pairing is such that 
	\begin{equation}\label{null} \re\left\langle \f{dz^2}{z^2}, \widetilde{\nu_s^0 }\right\rangle_{\mathcal{A}_{\gamma,s}}=\re\left\langle  \f{dz^2}{z^2}, \nu_s^0\right\rangle ~,\end{equation}
	where the right hand side stands for the coupling between holomorphic quadratic and Beltrami differentials on the universal covering of whole surface $X_s$. Now, recall, see Theorem \ref{dVR}, that the differential of the renormalized volume function at $X_s$ is the Schwarzian derivative of the developing map of the boundary at infinity of $M_s$: 
	\[(d V_R)_{X_s}(\mu)=\re\langle S(f), \mu\rangle~.\] 
	Moreover, on a smaller tube $U_s$ contained in the collar (note that we can assume $\nu_s^0$ to be zero outside it), by Theorem \ref{Schwarzian} and in the notations of Section \ref{sec:symmsch}, this can be expressed as a sum of a symmetric part and a correction term: 
	\[ S(f)=S\left(f_{\ell_{s}}\right)+q_{\ell_s}(z)dz^2~,\] 
	with 
	\begin{equation}\label{sym}
		S\left(f_{\ell_{s}}\right)=\f{1}{2z^2}\left(1+\f{4\pi^2}{\ell_{s}^2}\right)dz^2~,
	\end{equation} 
	and $\norm{q_{\ell_s}}_{\infty}$ bounded, such that
	\begin{equation}\label{O} \norm{q_{\ell_s}(z)}_{\infty}\in O\left(\f{e^{-\pi^2/(2\ell_s)}}{\ell_s^2}\right)~.\end{equation}
	The identities \eqref{dl}, \eqref{null} and \eqref{sym} give the first term in the right hand side of the equality in the statement:
	\[\re\left\langle  S\left(f_{\ell_{s}}\right), \nu_s^{0}\right\rangle = \re\left\langle \left(\f12+\f{2\pi^2}{\ell_{s}^2}\right)\f{dz^2}{z^2}, \nu_s^{0}\right\rangle =\left(\f{\pi}{4}+\f{\pi^3}{\ell_s^2}\right)d\ell_{\gamma}(\nu_s^0)~. \]
Then, it remains the term 
\[  E(s)\eqdef\re\left\langle q_{\ell_s}(z)dz^2, \nu_s^0\right\rangle~,\]
that, since $\nu_s^0=\nu_s \f{\partial}{\partial z} \otimes d\bar{z}$ with $\nu_s\in L^{\infty}(X_s)$ uniformly in $s$, by \eqref{O}, we estimate as
	 \[\abs{\re\left\langle q_{\ell_s}(z)dz^2, \nu_s^0\right\rangle} = \abs{\re\int_{\mathcal{A}_{\gamma, s}} q_{\ell_s}(z)\nu_s(z) dxdy} \leq \left(\sup_{z\in \mathcal{A}_{\gamma, s}}\abs{\nu_s(z)} \right)\norm{q_{\ell_s}(z)}_\infty Area_{hyp}(U_s )~,\]
	 as $\norm{q_{\ell_s}(z)}_\infty\in  O\left(\f{e^{-\pi^2/(2\ell_s)}}{\ell_s^2}\right)$ and the rest is uniformly bounded this concludes the proof\footnote{Note that the area term $Area_{hyp}(U_s)$ is uniformly bounded by $4$.}.

\end{proof}
	
	We can now prove the main result of this section.

\begin{customthm}{\ref{asymptotic}}
		Let $M_0$ be a convex co-compact hyperbolic $3$-manifold, and let $\gamma\subseteq\partial\bar{M}$ be a simple closed compressible curve of length $\ell_0(\gamma)\leq\varepsilon_0$ in its boundary. The composition of the renormalized volume with the grafting path $(M_s)_{s\in[0,\infty)}$ satisfies
	\[V_R(M_s)-V_R(M_0)= -\f{\pi^3}{\ell_{s}(\gamma)}+\f{\pi^3}{\ell_0(\gamma)}+\f{\pi}{4}(\ell_s(\gamma)-\ell_0(\gamma))+C(s)~,\]
		with $C(s)$ bounded and such that $C(s)\rar 0$ as $s\rar 0$. In particular, when $s\rightarrow \infty$ the renormalized volume diverges as $ -\f{\pi^3}{\ell_s(\gamma)}$.
\end{customthm}

\begin{proof}
First note that, in our notations, for any $t\in[0,s]$:
	\[\ell_t(\gamma)=\ell_{\gamma}(gr_{t\gamma}(X_0))~,\] so that the differential of $u(t)=\ell_{t}(\gamma)$ is 
	\[du=(d\ell_{\gamma})_{gr_{t\gamma}(X_0)}gr'_{t\gamma}(X_0)dt=(d\ell_{\gamma})_{X_t}\nu_t^0dt~.\]
	
	By integrating Lemma \ref{Gardiner} along the path we have:
	\begin{align*}
		V_R(M_s)&=V_R(M_0)+\int_0^s (d V_R)_{X_t}(\nu_t^0)dt \\
		&=V_R(M_0)+\int_0^s \left(\f{\pi}{4}+\f{\pi^3}{\ell_t^2}\right)(d\ell_{\gamma})_{X_t}(\nu_t^0)dt+\int_0^s E(t) dt \\
		&=V_R(M_0)+\int_{\ell_0}^{\ell_s} \left(\f{\pi}{4}+\f{\pi^3}{u^2}\right)du+\int_0^s E(t)dt\\
		&=V_R(M_0)+(\ell_s-\ell_0)\f{\pi}{4}-\f{\pi^3}{\ell_s}+\f{\pi^3}{\ell_0}+\int_0^sE(t)dt~,
	\end{align*}
	
	where the third equality is obtained by the change of variables $u(t)=\ell_t(\gamma)$ for the second term.
	
	We now define $C(s)=\int_0^sE(t)dt$. As $E(t)$ is bounded by Lemma \ref{Gardiner} we get that as $s\rar 0$ $C(s)\rar0$. Moreover, since again by Lemma \ref{Gardiner} $E(t)\in O\left(\f{e^{-\pi^2/(2\ell_t)}}{\ell_t^2}\right)$, for large enough $t>\bar t$:
	
	\[\abs{E(t)}\leq K_1 \f{e^{-\pi^2/(2\ell_t)}}{\ell_t^2}~,\]
	for some positive constant $K_1$. 
	
	Now, by using Lemma \ref{graftlength} we have
	
	\[ \f{e^{-\pi^2/(2\ell_t)}}{\ell_t^2}\leq \f{4(\pi+t)^2}{\pi^2\ell_0^2}e^{-\pi^2(\pi+t)/(2\pi\ell_0)}~.\] 
	Thus, for large enough $t>\bar t$ one obtains:
	\[\abs{E(t)}\leq K_2 e^{-\pi t/(2\ell_0)}t^2~,\]
	for some positive constant $K_2$.
By integrating the above by parts, this gives the boundedness of $\abs{C(s)}$, concluding the proof of the statement.  \end{proof}

\brem[Sharper Error in Theorem \ref{asymptotic}] \label{sharperboundschwarzian}By merging more carefully the estimate coming from the Fourier analysis of Section \ref{est}, and in particular Remark \ref{utube}, into the proof of Theorem \ref{Schwarzian}, one can observe that the coefficients of the tensors of type $Q_2(f(\ell))$ are not merely functions in $O\left(f(\ell)\right)$, but are in fact bounded by a universal constant times $f(\ell)$. Accordingly, Theorem \ref{Schwarzian} could be improved by stating that the norm of $q_{\ell}$ is bounded by a universal constant times $\f{e^{-\pi^2/(2\ell)}}{\ell^2}$. As a consequence, in our same notations of Theorem \ref{Schwarzian} and Lemma \ref{Gardiner}, the function $E(s)=\re\left\langle q_{\ell_s}(z)dz^2, \nu_s^0\right\rangle$ is bounded for sufficiently large $s$ by $K\f{e^{-\pi^2/(2\ell_s)}}{\ell_s^2}\leq K'e^{-\pi s/(2\ell_0)}s^2$, for some universal constants $K$ and $K'$, which could be explicitly computed. Therefore, for the function $C(s)$ of Theorem \ref{asymptotic}, one could extract the more precise bound \[\abs{C(s)}\leq\abs{ C(s_0)}+\int_{s_0}^s K'e^{-\pi t/(2\ell_0)}t^2~dt~.\] 
\erem

\newpage
\appendix

\section{The Osgood-Stowe tensor and the space of horospheres}
\label{standardmetrics}

In this section we provide a proof of the following elementary but useful statement.

\bthm \label{tm:standard}
  Let $f:\Omega\to \Omega'$ be a locally univalent holomorphic map, where $\Omega,\Omega'\subset \mathbb{CP}^1$. Let $h$, resp. $h'$ be a standard metric on $\Omega$, resp. $\Omega'$. Then
  $$ B(h,f^*h') = \re\left(S(f)\right)~. $$ 
\ethm

Here by {\em standard metric} we mean (the restriction of) a metric compatible with the conformal structure of $\mathbb{CP}^1$, which can be:
\begin{itemize}
\item a spherical metric on $\mathbb{CP}^1$,
\item a Euclidean metric on the complement of a point in $\mathbb{CP}^1$,
\item the hyperbolic metric in a round disk in $\mathbb{CP}^1$. 
\end{itemize}

 A direct proof of Theorem \ref{tm:standard} can be found in \cite[Proposition $4.2$]{Epstein-wvolume-OS}.
 We prefer here to provide a geometric argument using the space of horospheres, which moreover yields additional information. The argument presented here is based on \cite{horo}, and the main ideas are already in \cite{schwarzian}, although presented in a more straightforward and direct way here. 

\subsection{The ``full''  Osgood-Stowe tensor}

Let $\Omega\subset \R^d, d\geq 2$, and let $h, h'$ be two metrics on $\Omega$, with $h'=e^{2u}h$, for a function $u:\Omega\to \R$. Recall that the Osgood-Stowe tensor is defined in \cite{OS1992} as:
$$ B(h,h') = \hess_h(u) - du\otimes du - \frac 12\left((\Delta_h u)-\| du\|_h^2\right) h~. $$
It is traceless (with respect to $h$) by construction.

We introduce a ``non-traceless'' version of this tensor as follows:
$$ \bar B(h,h') =  \hess_h(u) - du\otimes du + \frac 12\| du\|_h^2 h~. $$
Clearly, the traceless part (with respect to $h$) of $\bar B$ is $B$.

This ``full'' Osgood-Stowe tensor, including a trace part, has a remarkable composition property, extending that of the usual, trace-free Osgood-Stowe tensor.

\blem\label{Bbaradditivity}
  Let $\Omega\subset \R^d$, and let $h_1, h_2, h_3$ be three conformal metrics on $\Omega$. Then
  $$ \bar B(h_1, h_3) = \bar B(h_1, h_2) + \bar B(h_2, h_3)~. $$
\elem

\begin{proof}
  Suppose that $h_2=e^{2u_1}h_1, h_3=e^{2u_2}h_2$, so that $h_3=e^{2u_1+2u_2}h_1$.
  By definition,
  \begin{eqnarray*}
    \bar B(h_1, h_3) & = & \hess_{h_1}(u_1+u_2) - d(u_1+u_2)\otimes d(u_1+u_2) + \frac 12\| d(u_1+u_2)\|_{h_1}^2 h_1 \\
                     & = & \hess_{h_1}(u_1)+\hess_{h_1}(u_2) - du_1\otimes du_1 - du_2\otimes du_2 \\
    & - & du_1\otimes du_2- du_2\otimes du_1 
    + \frac 12(\| du_1\|_{h_1}^2+\|du_2\|_{h_1}^2 +2h_1(du_1, du_2))h_1~.
  \end{eqnarray*}
  Moreover,
  \begin{equation}
    \label{eq:h1h2}
    \bar B(h_1,h_2) =  \hess_{h_1}(u_1) - du_1\otimes du_1 + \frac 12\| du_1\|_{h_1}^2 h_1~.
  \end{equation}
  To compute $\bar B(h_2, h_3)$, recall (see \cite{Besse}) that the Levi-Civita connection of $h_2$ is expressed as
  $$ D^2_xy = D^1_xy + du_1(x)y + du_2(y)x-h_1(x,y)D^1u~, $$
  where $D^1$ denotes both the Levi-Civita connection of $h_1$ and the gradient for $h_1$.
  As a consequence,\small
  $$ \hess_{h_2}(u_2)(x,y) = \hess_{h_1}(u_2)(x,y)   - du_1(x) du_2(y) - du_1(y)du_2(x) + h_1(x,y) h_1(du_1,du_2)~, $$\normalsize
  so that
    $$ \hess_{h_2}(u_2)= \hess_{h_1}(u_2) - du_1\otimes du_2 - du_1\otimes du_2 + h_1(du_1,du_2)h_1~. $$
  Therefore
  \begin{eqnarray*}
    \bar B(h_2, h_3) & = & \hess_{h_2}(u_2) - du_2\otimes du_2 + \frac 12\| du_2\|_{h_2}^2 h_2 \\
                     & = & \hess_{h_1}(u_2) - du_1\otimes du_2 - du_2\otimes du_1 + h_1(du_1,du_2)h_1 \\
    & - & du_2\otimes du_2 + \frac 12\| du_2\|_{h_1}^2 h_1~.  
  \end{eqnarray*}
  Adding this last equation to \eqref{eq:h1h2} shows the result.
\end{proof}

\subsection{The space of horospheres}

Following notations from \cite{horo}, we denote by $C^{d+1}_+$ the space of horospheres in $\HH^{d+1}, d\geq 2$. It can be identified with $S^d\times \R$, equipped with the degenerate metric $e^{2t}m_1$, where $m_1$ is the standard metric on $S^d$. This space $C^{d+1}_+$ is also equipped with a foliation by ``vertical lines'', which are defined as the kernel of the differential of the projection to $S^d$, or, equivalently, the kernel of the (degenerate) metric. Moreover, each leaf of this foliation is equipped with a distance and orientation through $t$. (Note that each vertical line corresponds to the set of horospheres centered at a given point at infinity, and the distance between two points is equal to the (constant) distance between the two corresponding horospheres.)

The isometric action of $O(d+1,1)$ on $\HH^{d+1}$ extends to an isometric action $C^{d+1}_+$, and for $d=2$ the action of $\mathrm{PSL}(2,\mathbb C)$ on $\HH^3$ extends to an isometric action on $C^3_+$. In addition, $C^{d+1}_+$ has a family of ``totally geodesic'' hypersurfaces, defined as the duals of points in $\HH^{d+1}$: for each $x\in \HH^{d+1}$, the set $H_x$ of points in $C^{d+1}_+$ corresponding to horospheres containing $x$ is one of those ``totally geodesic'' hypersurfaces. The group $O(d+1,1)$ acts transitively on the space of those hypersurfaces. More precisely, for each $x^*\in C^{d+1}_+$ and for all hyperplane $P\subset T_{x^*}C^{d+1}_+$ which is transverse to the vertical direction, there is a unique $x\in \HH^{d+1}$ such that $H_x$ is tangent along $P$ at $x$. 

There is a natural duality between hypersurfaces in $\HH^{d+1}$ and in $C^{d+1}_+$. If $H$ is an oriented hypersurface in $\HH^{d+1}$ which is {\em horospherically convex} -- its principal curvatures are strictly bigger than $-1$ at every point -- then the set of points in $C^{d+1}_+$ corresponding to the horospheres tangent to $H$ on the positive side is a smooth surface in $C^{d+1}_+$ denoted by $H^*$, see \cite{horo}. 

Let $H^*$ be a smooth surface in $C^{d+1}_+$, everywhere transverse to the vertical direction. The ambient metric on $C^{d+1}_+$ induces a Riemannian (non-degenerate) metric on $H^*$, which we denote by $I^*$. 

A crucial fact here is that any metric $h$ conformal to the spherical metric on $S^{d}$ can be realized as the first fundamental form of an embedded hypersurface $\Sigma^*$ in $C_+^{d+1}$, which is unique up to isometries of $C_+^{d+1}$ \cite[Theorem $6.3$]{horo}. We say that $\Sigma^*$ is associated to $h$.
 We can also define a notion of second fundamental form of $H^*$, denoted by $\II^*$, as follows. Given $x^*\in H^*$, let $H_0$ be the unique totally geodesic hyperplane tangent to $H^*$ at $x^*$. Let $u$ be the function defined over $H^*$ as the oriented distance along the vertical lines to $H_0$. Then
$$ \II^*_x = \hess_x(u)~. $$

The forms $I^*$ and $\II^*$ satisfy a modified version of the Gauss-Codazzi equation, see \cite{horo}. First, $\II^*$ is Codazzi with respect to $I^*$:
$$ d^{D^*}B^*=0~, $$
where $D^*$ is the Levi-Civita connection of $I^*$, and $B^*$ is defined by the condition that $\II^*=I^*(B^*\cdot, \cdot)=I^*(\cdot, B^*\cdot)$. Moreover, $\II^*$ satisfies a modified form of the Gauss equation with respect to $I^*$, in particular, for $d=2$, the curvature of $I^*$ is
$$ K^* = 1 - \tr(B^*)~. $$

We are particularly interested in three types of hypersurfaces, namely those associated to the standard metrics and which are
 {\em umbilic} in the sense that $\II^*$ is proportional to $I^*$:
\begin{itemize}
\item Totally geodesic surfaces, with $I^*$ of constant curvature $1$, already discussed.
\item The set of points dual to the horospheres tangent to a given horosphere (but with opposite orientation). For those surfaces $B^*=Id$ so that $K^*=0$.
\item The set of points dual to the horospheres tangent to a totally geodesic plane in $\HH^{d+1}$, for which $I^*$ is hyperbolic. 
\end{itemize}

\subsection{The Osgood-Stowe tensor and hypersurfaces}

The relation between the full Osgood-Stowe tensor and the space of horospheres is based on the following property.

\bthm \label{tm:difference}
  Let $\Omega\subset \mathbb{CP}^1$, let $h$ and $\bar h$ be two conformal metrics on $\Omega$, and let $\Sigma^*$ and $\bar \Sigma^*$ be the associated hypersurfaces in $C^{d+1}_+$, identified through the projection along the vertical lines. Let $\II^*$ and $\bar\II^*$ be their second fundamental forms. Then
  $$ \bar B(h,\bar h)=(\II^*+(1/2)I^*) - (\bar\II^*+(1/2)\bar I^*)~. $$
\ethm

The proof uses the following simple computation.

\blem \label{lm:spherical}
  Let $h_0$ be the Euclidean metric on $\R^{d}$, and let $h_1$ be the spherical metric defined on $\R^d$ as
  $$ h_1 = \frac{4h_0}{(1+r^2)^2}~, $$
  where $r$ denotes the Euclidean distance to $0$ in $\R^d$. Then
  $$ \bar B(h_0,h_1)= - \frac 12 h_1~. $$
\elem

\begin{proof}
  Let $u:\R^d\to \R$ be such that $h_1=e^{2u}h_0$, then $u=\log(2)-\log(1+r^2)$.
  Then
  $$ du = -\frac{d(r^2)}{1+r^2}~, $$
  so that
  $$ \hess(u) = -\frac{\hess(r^2)}{1+r^2} + \frac{d(r^2)\otimes d(r^2)}{(1+r^2)^2}~. $$
  However $\hess(r^2) = 2h_0$, while
  $$ \| du\|_{h_0}^2 = \frac{4r^2}{(1+r^2)^2}~. $$
  So we obtain that
  \begin{eqnarray*}
    \bar B(h_0, h_1) & = & \hess(u) - du\otimes du + \frac 12 \| du\|^2h_0 \\
                     & = & \frac{-2h_0}{1+r^2} + 2\frac{r^2}{(1+r^2)^2}h_0 \\
                     & = & \frac{-2}{(1+r^2)^2}h_0 \\
    & = & -\frac {h_1}2~. 
  \end{eqnarray*}
\end{proof}

\bcor\label{cr:spsp}
  Let $h_1$ and $\bar h_1$ be two conformal spherical metrics on $S^d$. Then $\bar B(h_1, \bar h_1)=(1/2)(h_1-\bar h_1)$.  
\ecor

\begin{proof}
  Let $h_0$ be a Euclidean metric on $S^d\setminus \{ p\}$, for some $p\in S^d$. Since $h_1$ and $\bar h_1$ are isometric to the standard spherical metric above, it follows from the previous lemma and Lemma \ref{Bbaradditivity} that
  \begin{eqnarray*}
    \bar B(h_1, \bar h_1) & = & \bar B(h_1, h_0) + \bar B(h_0, \bar h_1) \\
    & = & \frac 12 h_1 - \frac 12\bar h_1~. 
  \end{eqnarray*}
\end{proof}

\begin{proof}[Proof of Theorem \ref{tm:difference}]
  Let $\Sigma^*$ and $\bar\Sigma^* \subset C^{d+1}_+$ be the hypersurfaces associated to $h$ and $\bar h$, respectively. 

   Let $x\in \Sigma^*$, and let $\bar x$ be the point of $\bar\Sigma^*$ on the same vertical line as $x$. Let $P_1$ be the totally geodesic sphere tangent to $\Sigma^*$ at $x$, and let $\bar P_1$ be the totally geodesic sphere tangent to $\bar \Sigma$ at $\bar x$.

  We denote by $u$ and $\bar u$ the functions with graphs $\Sigma^*$ and $\bar \Sigma^*$, so that $h=e^{2u}m_1, \bar h =e^{2\bar u}m_1$, and let $u_1$ and $\bar u_1$ be the functions with graphs $P_1$ and $\bar P_1$, for $m_1$ the standard spherical metric. Finally, let $h_1=e^{2u_1}m_1$, $\bar h_1=e^{2\bar u_1}m_1$, and note that these are conformal spherical metrics \cite[Proposition $2.1$]{horo}.

  By definition of $\II^*$, we have
  $$ \II_x^* =\hess_x(u-u_1)=\bar B(h_1, h)~, $$
  $$ \bar \II_x^* =\hess_x(\bar u-\bar u_1)=\bar B(\bar h_1, \bar h)~. $$
  So
  \begin{eqnarray*}
    \bar \II_x^* - \II_x^* & = & \bar B(h_1, h) - \bar B(\bar h_1, \bar h) \\
                           & = & - \bar B(h, h_1) - \bar B(\bar h_1, \bar h) \\
                           & = & - \bar B(h, \bar h) + \bar B(h_1,\bar h_1) \\
                           & = & - \bar B(h, \bar h) + (1/2)(h_1-\bar h_1) \\
    & = & - \bar B(h, \bar h) + (1/2)(I_x^*-\bar I_x^*)~, 
  \end{eqnarray*}
  and the result follows.
\end{proof}

\bcor \label{cr:standard}
  Let $\Omega\subset S^d$, and let $h,\bar h$ be two standard metrics on $\Omega$. Then $B(h,\bar h)=0$. 
\ecor

\begin{proof}
  The hypersurfaces corresponding to standard metrics are all umbilic, that is, for those hypersurfaces, $\II^*$ is proportional to $I^*$.  According to Theorem \ref{tm:difference}, $\bar B(h,\bar h)$ is conformal to $m_1$, so that $B(h,\bar h)$, the traceless part of $\bar B(h,\bar h)$, is zero.
\end{proof}

\subsection{Proof of Theorem \ref{tm:standard}}

We now consider the case $d=2$. Let $\Omega\subset \C$, and let $f:\Omega\to \Omega'\subset \C$ be a holomorphic function. We already know (see \cite{OS1992}) that
$$ \re\left(S(f)\right) = B(h_0, f^*h_0)~, $$
where $h_0$ denotes the Euclidean metric on $\C$.

Let $h$ be a standard metric on $\Omega$, and let $h'$ be a standard metric on $\Omega'$. Then the composition rule \eqref{comprule} shows that
\begin{eqnarray*}
  B(h, f^*h') & = & B(h, h_0)+B(h_0, f^*h_0)+B(f^*h_0, f^*h') \\
  & = &  B(h, h_0)+B(h_0, f^*h_0)+f^*B(h_0, h')~. 
\end{eqnarray*}
Corollary \ref{cr:standard} indicates that
$$ B(h_0,h)=B(h_0, h')=0~, $$
so that
$$ B(h, f^*h') = B(h_0, f^*h_0) = \re\left(S(f)\right)~, $$
as claimed.

\thispagestyle{empty}
{\small
\markboth{References}{References}
\bibliographystyle{alpha}
\bibliography{mybib}{}

\begin{thebibliography}{BBVP23}

\bibitem[Ago04]{AG2004}
Ian Agol.
\newblock {Tameness of hyperbolic 3-manifolds}.
\newblock {\url{http://arxiv.org/abs/math/0405568}}, 2004.

\bibitem[BB24]{Epstein-wvolume-OS}
Martin Bridgeman and Kenneth Bromberg.
\newblock Epstein surfaces, $w$-volume, and the osgood-stowe differential.
\newblock 2024.

\bibitem[BBB19]{BBB2018}
Martin Bridgeman, Jeffrey Brock, and Kenneth Bromberg.
\newblock Schwarzian derivatives, projective structures, and the
  weil--petersson gradient flow for renormalized volume.
\newblock {\em Duke Mathematical Journal}, 168(5):867 -- 896, 2019.

\bibitem[BBB23]{bridgeman-brock-bromberg:gradient}
Martin Bridgeman, Jeffrey Brock, and Kenneth Bromberg.
\newblock The {W}eil-{P}etersson gradient flow of renormalized volume and
  3-dimensional convex cores.
\newblock {\em Geom. Topol.}, 27(8):3183--3228, 2023.

\bibitem[BBVP23]{bridgeman-bromberg-pallete:convergence}
Martin Bridgeman, Kenneth Bromberg, and Franco Vargas~Pallete.
\newblock Convergence of the gradient flow of renormalized volume to convex
  cores with totally geodesic boundary.
\newblock {\em Compos. Math.}, 159(4):830--859, 2023.

\bibitem[BC17]{BC2017}
Martin Bridgeman and Richard~D. Canary.
\newblock Renormalized volume and the volume of the convex core.
\newblock {\em Annales de l'Institut Fourier}, 67(5):2083--2098, 2017.

\bibitem[Bes87]{Besse}
Arthur Besse.
\newblock {\em Einstein Manifolds}.
\newblock Springer, 1987.

\bibitem[Bus10]{Bu1992}
Peter Buser.
\newblock {\em Geometry and spectra of compact {R}iemann surfaces}.
\newblock Modern Birkh\"{a}user Classics. Birkh\"{a}user Boston, Ltd., Boston,
  MA, 2010.
\newblock Reprint of the 1992 edition.

\bibitem[Can01]{C2001}
R.~D. Canary.
\newblock {The conformal boundary and the boundary of the convex core}.
\newblock {\em Duke Mathematical Journal}, 106(1):193 -- 207, 2001.

\bibitem[CG06]{CG2006}
Danny Calegary and David Gabai.
\newblock {Shrinkwrapping and the taming of hyperbolic 3-manifolds}.
\newblock {\em {Journal of American Mathematical Society}}, 19(2):385--446,
  2006.

\bibitem[CGS25]{CGS2024}
Tommaso Cremaschi, Viola Giovannini, and Jean-Marc Schlenker.
\newblock Filling {Riemann} surfaces by hyperbolic {Schottky} manifolds of
  negative volume.
\newblock {\em J. Geom. Phys.}, 217:18, 2025.
\newblock Id/No 105628.

\bibitem[DK07]{graftingrays}
Raquel Diaz and Inkang Kim.
\newblock Asymptotic behavior of grafting rays.
\newblock 2007.

\bibitem[Dum08]{dumas-survey}
David Dumas.
\newblock Complex projective structures.
\newblock 13:455--508, 2008.

\bibitem[DW08]{DumWolf}
David Dumas and Michael Wolf.
\newblock {Projective structures, grafting and measured laminations}.
\newblock {\em Geometry \& Topology}, 12(1):351 -- 386, 2008.

\bibitem[Eps84]{epsteinSurf}
Charles~L. Epstein.
\newblock Envelopes of horospheres and weingarten surfaces in hyperbolic
  3-space.
\newblock {\url{https://arxiv.org/abs/2401.12115}}, 1984.

\bibitem[FM11]{FM2011}
Benson Farb and Dan Margalit.
\newblock {\em {Primer on Mapping Class Groups}}.
\newblock {Princeton Mathematical Series}, 2011.

\bibitem[Gar75]{Gard}
Frederick~P. Gardiner.
\newblock Schiffer's interior variation and quasiconformal mapping.
\newblock {\em Duke Math. J.}, 42(1):371--380, 1975.

\bibitem[Gar00]{QTT}
Nikola Gardiner, Frederick P.;~Lakic.
\newblock {\em Quasiconformal Teichm\"uller Theory}, volume~76 of {\em
  Mathematical Surveys and Monographs}.
\newblock American Mathematical Society, first edition, 2000.

\bibitem[Hen08]{henselgraftingrays}
Sebastian~W. Hensel.
\newblock Iterated grafting and holonomy lifts of teichmueller space, 2008.

\bibitem[Hub16]{Hubbard2016}
John~H Hubbard.
\newblock {\em Teichm{\"u}ller theory and applications to geometry, topology,
  and dynamics}, volume 1-2.
\newblock Matrix Editions, 2016.

\bibitem[KM]{KM}
I.~Kra and B.~Maskit.
\newblock {\em Remarks on srojective structures: Riemann Surfaces and Related
  Topics.}

\bibitem[KS08]{S2008}
Kirill Krasnov and Jean-Marc Schlenker.
\newblock On the renormalized volume of hyperbolic 3-manifolds.
\newblock {\em Communications in Mathematical Physics}, (279):637--668, 2008.

\bibitem[KT92]{KamishimaTan1992}
Y.~Kamishima and S.~P. Tan.
\newblock Deformation spaces on geometric structures.
\newblock {\em Aspects of Low-Dimensional Manifolds}, 20:263--299, 1992.

\bibitem[Mar16a]{Ma2016}
Albert Marden.
\newblock {\em {Hyperbolic Manifolds: An introduction in 2 and 3 dimensions}}.
\newblock {Cambridge University Press}, 2016.

\bibitem[Mar16b]{Mar16}
Bruno Martelli.
\newblock {An introduction to geometric topology}.
\newblock {\url{arXiv:1610.02592v1}}, 2016.

\bibitem[McM98]{McMcomplexearth}
Curt McMullen.
\newblock Complex earthquakes and teichm\"uller theory.
\newblock {\em Journal of the American Mathematical Society}, 11:283--320,
  1998.

\bibitem[MT98]{MT1998}
Katsuhiko Matsuzaki and Masahiko Taniguchi.
\newblock {\em {Hyperbolic Menifolds and Kleinian Groups}}.
\newblock {Oxford University Press}, 1998.

\bibitem[Neh49]{Ne1949}
Zeev Nehari.
\newblock The schwarzian derivative and schlicht functions.
\newblock {\em Bull. Amer. Math. Soc.}, 55:545--551, 1949.

\bibitem[OS92]{OS1992}
Brad Osgood and Dennis Stowe.
\newblock The schwarzian derivative and conformal mapping of riemannian
  manifolds.
\newblock {\em Duke Mathematical Journal}, 67(1):57 -- 99, 1992.

\bibitem[Pal16]{VP2016}
Franco~Vargas Pallete.
\newblock Continuity of the renormalized volume under geometric limits.
\newblock arXiv:1605.07986, 2016.

\bibitem[Sch02]{horo}
Jean-Marc Schlenker.
\newblock Hypersurfaces in {$H\sp n$} and the space of its horospheres.
\newblock {\em Geom. Funct. Anal.}, 12(2):395--435, 2002.

\bibitem[Sch13]{compare}
Jean-Marc Schlenker.
\newblock The renormalized volume and the volume of the convex core of
  quasifuchsian manifolds.
\newblock {\em Math. Res. Lett.}, 20(4):773--786, 2013.
\newblock Corrected version available as arXiv:1109.6663v4.

\bibitem[{Sch}17]{schwarzian}
Jean-Marc {Schlenker}.
\newblock {Notes on the Schwarzian tensor and measured foliations at infinity
  of quasifuchsian manifolds}.
\newblock {\em ArXiv e-prints}, August 2017.

\bibitem[SW22]{averages}
Jean-Marc Schlenker and Edward Witten.
\newblock No ensemble averaging below the black hole threshold.
\newblock {\em J. High Energy Phys.}, (7):Paper No. 143, 50, 2022.

\end{thebibliography}

}

\end{document}